\def\TwoCol{1}

		\documentclass[twocolumn,a4paper,11pt]{article}
\usepackage{newtxtext} 
\usepackage[a4paper, top=2.5cm, bottom=2.5cm, left=1.3cm, right=1.3cm]{geometry}
\usepackage{cite}
\usepackage{authblk}
\usepackage{graphicx}
\usepackage{amsmath}
\usepackage{fancyhdr}
\usepackage{hyperref}

\let\oldbibliography\thebibliography
\renewcommand{\thebibliography}[1]{%
  \oldbibliography{#1}%
  \setlength{\itemsep}{0.3ex}
}

\usepackage{etoolbox}
\AtBeginEnvironment{tabular}{\scriptsize}
\apptocmd{\thebibliography}{\small}{}{}

\pretocmd{\bibliography}
    {
      \bibliographystyle{unsrturl}
  	}{
    }{
    }

\usepackage{amsthm}

\hypersetup{
    colorlinks=true, 
    linkcolor=blue,
    urlcolor=blue, 
    citecolor=blue,
    pdfauthor={Herrmann-Wicklmayr & Flaßkamp},
}

\theoremstyle{plain}

\theoremstyle{definition}

\theoremstyle{remark}

\makeatletter

\newcommand{\titleblockwidth}{0.8\textwidth}

\newcommand{\articletype}[1]{}
\newcommand{\received}[1]{}
\newcommand{\revised}[1]{}
\newcommand{\accepted}[1]{}
\newcommand{\journal}[1]{}
\newcommand{\volume}[1]{}
\newcommand{\copyyear}[1]{}
\newcommand{\startpage}[1]{}
\newcommand{\jnlcitation}[1]{}
\newcommand{\cname}[1]{}
\newcommand{\ctitle}[1]{}
\newcommand{\cjournal}[1]{}
\newcommand{\cvol}[1]{}
\newcommand{\FirstPg}[1]{}
\newcommand{\LastPg}[1]{}
\newcommand{\fundinginfohead}[1]{}
\newcommand{\fundinginfoheadtext}[1]{}
\newcommand{\fundingInfo}[1]{}

\let\oldtitle\title
\renewcommand{\title}[1]{\oldtitle{#1}\gdef\@title{#1}}
\newcommand{\titlemark}[1]{\gdef\@titlemark{#1}}
\newcommand{\authormark}[1]{\gdef\@authormark{#1}}
\let\address\affil
\newcommand{\orgdiv}[1]{#1}
\newcommand{\orgname}[1]{#1}
\newcommand{\orgaddress}[1]{#1}

\newcommand{\city}[1]{#1}

\newcommand{\country}[1]{#1}
\newcommand{\corres}[1]{\gdef\@corres{#1}}
\newcommand{\presentaddress}[1]{\gdef\@presentaddress{#1}}
\newcommand{\email}[1]{Email: \href{mailto:#1}{#1}}
\newcommand{\keywords}[1]{\gdef\@keywordscontent{#1}}

\gdef\@abstracttitle{}
\gdef\@abstractcontent{}
\renewcommand{\abstract}[2][Abstract]{%
  \gdef\@abstracttitle{#1}%
  \gdef\@abstractcontent{#2}%
}

\renewcommand{\maketitle}{%
    \twocolumn[
      \begin{@twocolumnfalse}
        \pagestyle{fancy}%
        \hfill%
        \begin{minipage}{\titleblockwidth}

            \begin{center}%
                {\LARGE\bfseries \@title \par}
                \vspace{1.5em}
                {\large\@author \par}
            \end{center}
            \vspace{0.0em}
            
            \ifdefined\@corres
                \noindent\rule{\textwidth}{0.4pt}\par
                \small\@corres
                \par\vspace{0.25\baselineskip}
            \fi
            \ifdefined\@presentaddress
                \small\@presentaddress
                \par\vspace{-0.5em}
            \fi
            \ifdefined\@corres
                \noindent\rule{\textwidth}{0.4pt}\par
            \fi
            \vspace{0em}

            \if\relax\detokenize{\@abstractcontent}\relax\else
                \begin{center}\bfseries \@abstracttitle\end{center}\par\vspace{-0.5em}
                \small \@abstractcontent \par
                \vspace{1.5em}
            \fi
            
            \ifdefined\@keywordscontent
                \par\noindent\textbf{Keywords:}\enspace \@keywordscontent\par
            \fi

        \end{minipage}%
        \hfill%
      \end{@twocolumnfalse}
      \vspace{2em}%
    ]
    \thispagestyle{plain}%
}

\fancypagestyle{plain}{\fancyhf{}}

\makeatother

		\def\AddAckn{0}
	
\newlength{\ShortenLastPage}%
	\providecommand{\pgfsyspdfmark}[3]{} 

\label{packages}
	\usepackage{microtype}
	\usepackage[utf8]{inputenc}
	\usepackage{amsmath}
	\usepackage{amssymb}
	\usepackage{algorithm}
	\usepackage{mathtools}
	\usepackage{aligned-overset}
	\usepackage{blindtext}
	\usepackage{csquotes}
	\usepackage[dvipsnames,table]{xcolor}
	\usepackage{enumitem}
		\setlist[enumerate]{itemsep=0.1em,topsep=0.4em,leftmargin=*,labelindent=1em}		
		\setlist[itemize]{itemsep=0.1em,topsep=0.4em,leftmargin=*,labelindent=1em,label={$\bullet$}}		
	\usepackage{siunitx}
		\let\oldunit\unit
		\renewcommand{\unit}[1]{\,\oldunit{#1}}
	\usepackage{resizegather}
	\usepackage{subcaption} 
		
	\usepackage{standalone} 
	\usepackage{tabularray}
		\UseTblrLibrary{booktabs}
	\usepackage{lipsum}
	
	\usepackage{booktabs} 
	\usepackage{eurosym}
	\usepackage{adjustbox} 
	\usepackage{makecell} 
	
	\usepackage{matlab-prettifier}
	\usepackage{graphicx}  
		\definecolor{udsblue}    {RGB}{0,72,119}
		\definecolor{udsdarkblue}{RGB}{1,40,63}
		\definecolor{slightgray}{RGB}{215,215,215}
		\definecolor{whitegray}{rgb}{0.985,0.985,0.985}
	\usepackage{tikz}

\label{boolean variables, penalties and whitespace}
	\newif\ifAppendix
	\newif\ifTodos   
	\newif\ifFormReview        
	\newif\ifIFACClassFile
	\newif\ifOldIEEEClassFile 
	\newif\ifAcronymsVisible
	\newif\ifAnimation
	\newif\ifAcronymCount
	
	\newcommand{\locDefPenalites}[1]{%
		\begingroup%
			\binoppenalty = 700%
			\relpenalty   = 500%
			#1%
		\endgroup%
	}

\label{standard macros/user-defined commands}
	
	\usepackage{listings}
		\newcommand{\inv}[1][1]{^{\raisebox{-0.0ex}{$\scriptscriptstyle - #1$}}}
		\newcommand{\tr}{^{\raisebox{0.37ex}{$\scriptscriptstyle\intercal$}}}

	\newcommand{\diag}{\operatorname{diag}}
	\DeclareMathOperator*{\argmin}{arg\,min}
	\renewcommand{\(}{\left(}
	\renewcommand{\)}{\right)}
	\newcommand{\numberthis}{\refstepcounter{equation}\tag{\theequation}}

	\newcommand{\tagref}[1]{\ref{#1}}
	
	\newcommand{\concSep}[1][]{ \, #1| \, }

	\newcommand{\ceil}[1]{\left\lceil #1 \right\rceil}

	\newcommand{\tableFontSizeWiley}{\fontsize{8}{10}\selectfont}

\label{define list}

	\newcounter{NonFloatListCtr}
	\newenvironment{NonFloatList}[1][]{%
	    \refstepcounter{NonFloatListCtr}%
	    \par\medskip
	    \hrule height 1.2pt \smallskip\nobreak      
	    \noindent\textbf{List \theNonFloatListCtr}%
	    \ifx#1\empty\else\ #1\fi%
	    \vspace{0.1\baselineskip}%
	    \par\hrule height 1.2pt \medskip\nobreak  
	    \noindent
	}{%
	    \par\medskip
	    \hrule height 1.2pt \par\medskip            
	}

\label{uncommented}

	\AcronymsVisiblefalse
	\AcronymCounttrue
	\FormReviewfalse

	\usepackage{hyperref}

	\hypersetup{
		draft=false,
		colorlinks=true,
		hypertexnames=false,
		linkcolor = black,
		citecolor = black,
		urlcolor  = black,
	}

	\usepackage{autonum}

\usepackage[nameinlink,noabbrev,capitalise]{cleveref}
	
	\crefname{theorem}{Theorem}{Theorems}
	\crefname{assumption}{Assumption}{Assumptions}
	\crefname{definition}{Definition}{Definitions}
	\crefname{convention}{Convention}{Conventions}
	\crefname{eexample}{Example}{Examples}
	\crefname{subsection}{Subsection}{Subsections}
	\crefname{table}{Table}{Tables}
	\crefname{NonFloatListCtr}{List}{Lists}
	\crefname{stepsi}{Step}{Steps} 

	\newenvironment{bThm}[2][]{%
	\def\envPH{#2}%
	\begin{\envPH}[#1]
}{%
	\end{\envPH}
}

	\label{prepare acronyms (acc, accp)}
\usepackage[nolist,nohyperlinks,printonlyused]{acronym}

	\makeatletter
	\AtBeginDocument{%
		\renewcommand*\AC@verridelabel[1]{%
			\@bsphack
			\protected@write\@auxout{}{\string\AC@undonewlabel{#1}}%
			\protected@write\@auxout{}{\string\AC@undonewlabel{#1@cref}}%
			\ifdefined\ltx@label\ltx@label{#1}\else\label{#1}\fi
			\AC@overriddenmessage rs{#1}%
			\AC@overriddenmessage rs{#1@cref}%
			\@esphack}}
	\makeatother

\definecolor{colAcro}{RGB}{0, 110, 14} 	
\ifAcronymsVisible
		\newcommand{\highlight}[1]{\textcolor{colAcro}{#1}}

		\let\oldacl\acl	  \renewcommand*\acl[1]{\highlight{\oldacl{#1}}}
		\let\oldacp\acp	  \renewcommand*\acp[1]{\highlight{\oldacp{#1}}}
		\let\oldaclp\aclp \renewcommand*\aclp[1]{\highlight{\oldaclp{#1}}}
		\let\oldAcl\Acl	  \renewcommand*\Acl[1]{\highlight{\oldAcl{#1}}}
		\let\oldAcp\Acp	  \renewcommand*\Acp[1]{\highlight{\oldAcp{#1}}}
		\let\oldAclp\Aclp \renewcommand*\Aclp[1]{\highlight{\oldAclp{#1}}}
\else
	\newcommand{\highlight}[1]{#1}
\fi

\newcommand{\contentAcronymUsageTable}{}

	\newcommand{\incrementacronym}[1]{%
	    \stepcounter{acronymcount#1}%
	}

		\makeatletter
		\newcommand{\noLongAcronym}[2]{%
		    \if\relax\detokenize{#1}\relax
		        #2
		    \else
		        \@noLongAcronym#1\@nil
		    \fi
		}
		\def\@noLongAcronym#1#2\@nil{%
		    \if l#1%
		    \else
		        #2
		    \fi
		}
		\makeatother

	\makeatletter
	\newcommand{\myacro}[3][]{%
		\edef\@ctrname{acronymcount#2}%
		\expandafter\newcounter\expandafter{\@ctrname}%
		\expandafter\gdef\csname increment#2\endcsname{%
			\incrementacronym{#2}%
		}%
		\acro{#2}{#3}%
		\gappto{\contentAcronymUsageTable}{
		 	\displayacronymcount[#1]{#2}%
		}%
	}
	
	\newcommand{\myacroWO}[3][]{
		\edef\@ctrname{acronymcount#2}%
		\expandafter\newcounter\expandafter{\@ctrname}%
		\expandafter\setcounter\expandafter{\@ctrname}{-1}
		\expandafter\gdef\csname increment#2\endcsname{%
		}%
		\acro{#2}{#3}%
		\gappto{\contentAcronymUsageTable}{
		 	\displayacronymcount[#1]{#2}%
		}%
	}
	
	\newcommand{\condWO}[4]{%
		\ifnum\csname c@acronymcount#2\endcsname=-1\relax
		   	\csname#3l#4\endcsname{#2}%
		\else
			\csname#3#1#4\endcsname{#2}%
		\fi%
		\noLongAcronym{#1}{\csname increment#2\endcsname}
	}
	
	\newcommand{\acc}[2][]{%
		\condWO{#1}{#2}{ac}{}
	}
	\newcommand{\accp}[2][]{%
		\condWO{#1}{#2}{ac}{p}
	}	
	\newcommand{\Acc}[2][]{%
		\condWO{#1}{#2}{Ac}{}%
	}
	\newcommand{\Accp}[2][]{%
		\condWO{#1}{#2}{Ac}{p}%
	}
	\makeatother

	\newcommand{\displayacronymcount}[2][]{%
		#2
		&\edef\temp{\noexpand\arabic{acronymcount#2}}\temp
	   	& #1
   		\\
	}

	\definecolor{burntorange}{rgb}{0.8, 0.33, 0.0}
\ifFormReview
	\newcommand{\form}[1]{\textcolor{burntorange}{#1}}

\else
	\newcommand{\form}[1]{#1}

\fi

\let\oldLipsum\lipsum
	\renewcommand{\lipsum}[1][]{\textcolor{lightgray}{\oldLipsum[#1]}}

	\usepackage[titles]{tocloft}
	\makeatletter
		\renewcommand{\@pnumwidth}{4em}
		\renewcommand{\@tocrmarg}{4em}
	\makeatother
	\AtEndPreamble{
    \ifdefined\DIFaddbegin

    \fi
}

	\usepackage{twoopt}

\definecolor{PFgreen}{RGB}{0,150,0}
\definecolor{PFlightgreen}{RGB}{0,210,0}

\newcommand{\Rset}{\mathbb{R}}
\newcommand{\task}[1]{\fbox{#1}}
\newcommand{\uc}[1]{\expandafter\MakeUppercase\expandafter{#1}}

\newcommand{\id}[1]{\mathbb{I}_{#1}} 
\newcommand{\opt}{^*}

\newcommand{\n}{k} 

	\newcommand{\FunitBox}{\mathbb{F}_{\mathrm{UB}}}
	\newcommand{\FunitSimplex}{\mathbb{F}_{\mathrm{US}}}

	\newcommand{\Seq}[1]{\boldsymbol{#1}}
	\newcommand{\Set}[1]{\mathbf{#1}}

	\newcommand{\xSys}{z} 
	\newcommand{\uSys}{v} 
	
	\newcommand{\uSeq}{\Seq{\uSys}}
	\newcommand{\uSet}{\Set{\uc{\uSys}}}

	\newcommand{\xSet}{\Set{\uc{\xSys}}}

	\newcommand{\term}{_0} 

	\newcommand{\SO}{\mathrm{SO}}
	
	\newcommand{\NP}{\mathrm{NP}}
	\newcommand{\UP}{\mathrm{UP}}

	\newcommand{\DM}{\mathrm{DM}}

	\newcommand{\DC}{\mathrm{DC}}
	\newcommand{\lb}{\mathrm{lb}}
	\newcommand{\ub}{\mathrm{ub}}
	
	\newcommand{\FHP}{\mathrm{FHP}}
	\newcommand{\meas}{}
	
	\newcommand{\cNPUP}{C_{\NP,\UP}}
	\newcommand{\WS}{\mathrm{WS}}

	\newcommand{\hpAcc}{}
	\newcommand{\wFHP}{\hpAcc{w}} 
	\newcommand{\bFHP}{\hpAcc{b}} 
	\newcommand{\dFHP}{\hpAcc{d}}
	
	\newcommand{\fScal}{f_{\mathrm{scale}}}
	\newcommand{\fNormal}{f_{\perp}}

	\newcommand{\J}{\mathcal{J}}
	\newcommand{\JMO}{\J}
	
	\newcommand{\xMO}{x}
	\newcommand{\JStarTime}{\J^*(t_k)}
	\newcommand{\lInR}[1]{} 
	\newcommand{\betaPref}[1][]{\beta_{#1}}
	
	\newcommand{\conePS}[1][\realize]{#1{\mathsf{K}}}
	\newcommand{\baryCO}{\lambda}
	
		\newcommand{\param}{\hat}
		\newcommand{\realize}{\check}
		\newcommand{\paramAccent}{hat\xspace}
		\newcommand{\realizeAccent}{check\xspace}

	\newcommand{\signD}{q}
	
	\newcommand{\wChim}[1][]{#1{\eta}}
	\newcommand{\bChim}{\rho}
	\newcommand{\wChimQuasi}{\wChim[\tilde]}
	\newcommand{\wChimVisual}{\wChim_{\mathrm{vis}}}

		\newcommand{\xOP}{x}
		\newcommand{\PO}{\mathrm{P}}
		\newcommand{\costFun}{J}
		\newcommand{\costFunVec}{\mathcal{J}}
		\newcommand{\nObj}{{n_J}}
		\newcommand{\decSetFeas}{\Set{\uc{\xOP}}}
		\newcommand{\decSetFeasMPC}{\decSetFeas}
		\newcommand{\decSetPareto}{\Set{\uc{\xOP}}_{\PO}}
		\newcommand{\decSetFeasMOMPC}{\tilde{\Set{\uc{\xOP}}}}
		\newcommand{\imSetFeas}{\Set{J}}
		\newcommandtwoopt{\imSetPareto}[2][][]{#1{\Set{J}}_{\PO#2}}
		\newcommand{\feasInitialStates}{_{\mathrm{FI}}}

	\newcommand{\fCT}[1][]{#1{f}_{\mathrm{CT}}}
	\newcommand{\fDT}[1][]{#1{f}}
	\newcommand{\hMicro}{h_{\mathrm{micro}}}
	\newcommand{\hMacro}{h_{\mathrm{macro}}}
	\newcommand{\scale}{s}
	\newcommand{\scaleTime}{\scale_{\mathrm{time}}} 
	\newcommand{\dT}{\Delta T_{\mathrm{cond}}}
	\newcommand{\CS}{\breve} 
	\newcommand{\PS}{} 
	\newcommand{\airInside}{\alpha_{\mathrm{mix}}}

\newcommand{\pSmoothAbs}{L}
\newcommand{\cSS}[1][]{_{*#1}} 

	\newcommandtwoopt{\Qstar}[2][][]{#1{Q}_{#2}^*}
	\newcommandtwoopt{\Q}[2][][]{#1{W}_{\xSys#2}}
	\newcommandtwoopt{\R}[2][][]{#1{W}_{\uSys#2}}
	\newcommand{\DeltaQ}[1][]{\Delta #1{Q}}

	\newcommand{\termCostFun}{F} 
	\newcommand{\termContr}{\kappa}
	\newcommand{\expkLoc}{}	
	\newcommand{\kMax}{K} 

\newcommand{\circled}[1]{[#1]}
\newcommand{\testSetup}[2]{[#1,#2]}

\if\TwoCol1
	\newcommand{\proofColDepAlignBefore}{ \\ } 
	\newcommand{\proofColDepAfter}{ & }        
\else
	\newcommand{\proofColDepAlignBefore}{  }   
	\newcommand{\proofColDepAfter}{}           
\fi

\AddToHook{env/samepage/begin}{%
	\AddToHook{para/before}{\nopagebreak}%
	\AddToHook{para/after}{\nopagebreak}%
}

\newenvironment{mat}{
  \begin{bmatrix}
}{
  \end{bmatrix}
}

\NewDocumentCommand{\matInline}{O{} O{} m}{%
	#1[#3#2]
}

\newsavebox{\tempbox}%
\newcommand{\InfoBoxFontSize}{\normalsize}
\newcommand{\InfoBox}[1]{%
	\begin{minipage}[b][\ht\tempbox][c]{0.10\textwidth}
		\setlength{\fboxsep}{3pt}
		\fbox{%
			\begin{minipage}{0.9\linewidth}
				\centering%
				\InfoBoxFontSize%
				#1
			\end{minipage}%
		}%
	\end{minipage}
	\hfill
}%

\newcommand{\hyperanchor}[1]{%
	\phantomsection%
	\label{#1}
}

\newlist{steps}{enumerate}{1} 

	\usepackage{environ}
\usepackage{anyfontsize}
\usepackage{fp}

\newcommand*\newfontsize{}
\newcommand*\newbaselineskip{}
\newsavebox{\alignatResizeBox}
\NewEnviron{alignatResize}[3]{%
    \sbox{\alignatResizeBox}{%
        $%
        \begin{alignedat}{#1}
            \BODY
        \end{alignedat}%
        $%
    }

    \FPdiv\scalefactor{\number\wd\alignatResizeBox}{\number\linewidth}%

    \FPifgt{\scalefactor}{1.001}%
        \FPmul\scalefactor{\scalefactor}{1.27}
    
        \begingroup
            \FPeval{\newfontsize}{#2/scalefactor}
            \FPeval{\newbaselineskip}{newfontsize * 1.2}%
            \fontsize{\newfontsize pt}{\newbaselineskip pt}\selectfont%
            \begin{alignat}{#1}
                \BODY
                #3
            \end{alignat}%
        \endgroup
    \else
        \begin{alignat}{#1}
            \BODY
            #3
        \end{alignat}%
    \fi
}

	\begin{acronym}
	
	\myacroWO{MO}{multi-objective}
	\myacroWO{OP}{optimization problem}
	\myacroWO{MOO}{multi-objective optimization}
	\myacroWO{MOOP}{multi-objective optimization problem}
	\myacro{pOP}{parameterized optimization problem}	
	\myacro[its MPC ;)]{MPC}{model predictive control}
	\myacro{MOMPC}{multi-objective model predictive control}
	
	\myacro{PF}{Pareto front}
	\myacro{IM}{individual minimum}
		\acroplural{IM}[IM]{individual minima}
	\myacro{CHIM}{convex hull of individual minima}
	\myacro{POS}{Pareto optimal solution}
	\myacroWO{ES}{efficient solution}
	\myacroWO{POP}{Pareto optimal point}
	\myacro{DM}{decision-making}
	\myacroWO{NIS}{normalized image space}
	
	\myacro[used as index]{NP}{nadir point}
	\myacro[used as index]{UP}{utopia point}
	\myacro[used as index]{SO}{shooting origin}
	\myacroWO{SDV}{shooting direction vector}
	
	\myacroWO{KP}{knee-point}
	\myacroWO{QKP}{quasi-knee-point}
	\myacro[used as index]{FHP}{furthest-to-hyperplane-point}
	\myacroWO[used as index]{NCHIM}{nadir-CHIM}
	
	\myacro{WS}{weighted-sum}
	\myacro{NBI}{normal boundary intersection}
	\myacro{SRI}{spherical ray intersection}
	\myacro{PS}{Pascoletti-Serafini}
	
	\myacroWO{CT}{continuous-time}
	\myacroWO{DT}{discrete-time}
\end{acronym}

	\graphicspath{
		{./figures/}
		{./figures/sim_cases/}
		{./figures/DM_comparison/}
	}
	\newcommand{\FigureFileExt}{pdf}
	
	\makeatletter
		\@ifundefined{medsize}{}{}
	\makeatother
	
	\newlength{\WidthBase}

\begin{document}


\articletype{Article Type}

\received{Date Month Year}
\revised{Date Month Year}
\accepted{Date Month Year}
\journal{Journal}
\volume{00}
\copyyear{2023}
\startpage{1}

\raggedbottom

\title{%
    Individual Minima-Informed Multi-Objective Model Predictive Control for Fixed Point Stabilization
}
\titlemark{Individual Minima-Informed sMOMPC}

\author[1]{Markus Herrmann-Wicklmayr}
\author[1]{Kathrin Flaßkamp}
\authormark{Herrmann-Wicklmayr \& Flaßkamp}

\address[1]{%
    \orgdiv{Systems Modeling and Simulation},\ %
    \orgname{Saarland University},\ %
    \orgaddress{%
	    \city{Saarbrücken},\ %
	    \country{Germany}%
    }%
}

\corres{%
    Correspondence:\ %
    Markus Herrmann-Wicklmayr,\ %
    Systems Modeling and Simulation,\ %
    Saarland University,\ %
    Saarbrücken,\ %
    Germany.\ %
    \email{markus.herrmannwicklmayr@uni-saarland.de}%
}

\presentaddress{%
	Present address:\ %
	Systems Modeling and Simulation,\ %
	Saarland University,\ %
	Saarbrücken,\ %
	Germany.%
}

	\fundinginfohead{\textbf{Acknowledgments/Funding}}
	\fundinginfoheadtext{}
	\fundingInfo{The authors received no specific funding for this work.}

\jnlcitation{%
    \cname{%
        \author{Herrmann-Wicklmayr M.},\ %
        \author{Flaßkamp K.}
    }.
    \ctitle{Individual Minima-Informed Multi-Objective Model Predictive Control for Fixed Point Stabilization.}
    \cjournal{\it Optimal Control Applications and Methods.}
    \cvol{2025;00(00):1--18}.
}

\abstract[Abstract]{%
In multi-objective model predictive control (MOMPC), a single Pareto-optimal solution has to be selected at each sampling instant by an automated a priori decision-making (DM) mechanism, since real-time constraints do not permit a human decision-maker to be involved.
This paper considers MOMPC for fixed point stabilization, where the DM translates a high-level preference into the solution to be implemented, while closed-loop stability has to be preserved.
To this end, an approach called individual minima-informed DM is introduced.
These methods leverage the individual minima, which are characteristic Pareto-optimal points, to facilitate the translation of a high-level preference to a point on the Pareto front.
They can be implemented through two sequential optimizations, regardless of the number of objectives, thereby facilitating a real-time implementation of MOMPC.
Two novel methods are proposed and compared with four established ones.
These methods are embedded into a stabilizing MOMPC scheme that builds on a quasi-infinite horizon framework from the literature.
Closed-loop stability is proven under a descent condition that is less restrictive than the one in the original framework and retains the largest DM space among all descent conditions that prevent any objective from increasing.
The required terminal ingredients can be constructed with a practical procedure.
In a numerical case study, the resulting scheme stabilizes the fixed point and illustrates the potential benefit of adapting the preference online.
}


\maketitle

\acresetall

\if\TwoCol1
	\setlength{\WidthBase}{\linewidth}
	\newcommand{\ColDepLinebreak}{\linebreak}
\else
	\setlength{\WidthBase}{0.5\linewidth}
	\newcommand{\ColDepLinebreak}{}
\fi

\section{Introduction}
\label{sec:Intro}
\Acc{MO} optimization is appealing because it yields a set of optimal trade-offs instead of a single optimal solution.
These trade-offs arise because several conflicting objectives are optimized simultaneously, so that, in general, no solution is optimal with respect to all objectives at once \cite{miettinen_nonlinear_1998,ehrgott_multicriteria_2005}.
The solutions of interest are the so-called \accp{ES}, also referred to as Pareto-optimal solutions, for which no objective can be improved without worsening at least one other objective, and their images form the \acc{PF} (see \cref{def:Nondom}).
Hence, a single solution has to be chosen from the \acc{PF}, which requires a \acc{DM} process.
\Acc{MOMPC}\acused{MPC} transfers this concept to the control of dynamical systems: at each sampling instant, an \acc{MO} optimal control problem is solved, and one of its \accp{ES} is applied to the system \cite{peitz_survey_2018,grune_multiobjective_2019}.
Consequently, both the computation of (parts of) the \acc{PF} and the \acc{DM} have to be carried out under real-time constraints, which neither allow for a dense sampling of the \acc{PF} nor permit a human decision-maker to be involved.

\acc{DM} methods can be classified by the stage at which preference information enters the solution process.
Following \cite{miettinen_nonlinear_1998}, it is distinguished between no-preference, \emph{a priori}, \emph{a posteriori}, and interactive methods.
The latter are not applicable in \acc{MOMPC}, though, since no human decision-maker can be involved due to the real-time constraints.
Most \acc{MOMPC} approaches rely on \emph{a priori} \acc{DM}.
In industrial practice and large parts of the literature, the objectives are even collapsed into a single-objective \acc{MPC} by fixed weighting factors, which lack an intuitive meaning \cite{marler_weighted_2010}, require extensive manual tuning, and cannot adapt to changing operating conditions.
A widely adopted strategy minimizes the distance to a reference point in the image space.
This point is typically computed in advance and kept fixed \cite{zavala_stability_2012}, which prevents an adaptation to changing operating conditions.
Alternatively, it is recomputed as the \acc{UP} at each step \cite{flores-tlacuahuac_multiobjective_2012}.
Other approaches employ the \acc{WS} method with time-varying weights \cite{bemporad_multiobjective_2009,shen_path-following_2019}, whose adaptation laws are mostly application-specific, or lexicographic optimization \cite{he_stability_2015}, which requires several sequential optimizations.
\emph{A posteriori} approaches approximate the \acc{PF}, typically by \acc{MO} evolutionary algorithms, and select a point afterwards by a predefined rule, e.g., fuzzy rules \cite{valera_garcia_intelligent_2012}, the minimum distance to a reference point \cite{laabidi_multi-criteria_2008}, or knee-point criteria \cite{schmitt_incorporating_2022}.
Due to the computational effort of the approximation, however, real-time requirements are often difficult to meet \cite{schmitt_incorporating_2022}.
Hence, the strict real-time constraints induce a trade-off between the amount of information about the \acc{PF} that can be computed and the quality of the resulting decision.
Furthermore, the \acc{MOMPC} approaches mentioned above rarely examine how accurately a specified preference is translated into the selected \acc{ES}, let alone quantify this translation.
\\
Between \emph{a priori} and \emph{a posteriori} methods lies a class of approaches, called \textit{informed \acc[l]{DM}} in this paper, where a few characteristic points on the \acc{PF} are computed and leveraged to determine the \acc{ES} ultimately used.
A well-known representative of this class is (distance-based) \acc[l]{KP}-guided optimization
\cite{das_characterizing_1999,li_knee_2020,chiu_minimum_2016}%
.
The characteristic points used in this approach are the \accp{IM}, i.e.,~the objective vectors obtained by optimizing each objective individually.
Altogether, this method requires one more optimization than the number of objectives, i.e.~it does not suffer from the curse of dimensionality, as a dense sampling of the entire \acc{PF} would.

For fixed point stabilization, the \acc{DM} method has to be embedded into a stabilizing \acc{MOMPC} scheme.
While stability proofs exist
	for linear systems \cite{bemporad_multiobjective_2009},
	for reference point-based (utopia-tracking) formulations \cite{zavala_stability_2012},
	for lexicographic formulations \cite{he_stability_2015},
	for general nonlinear systems \cite{grune_multiobjective_2019},
	and for nonlinear systems under additional assumptions on dissipativity \cite{eichfelder_relaxed_2023,krugel_dissipativity-based_2023}, 
it remains to show whether the \accp{IM}-informed \acc{DM} methods developed in this paper can be integrated into a general nonlinear \acc{MOMPC} framework in a way that preserves closed-loop stability.
To this end, this work builds on the nonlinear \acc{MOMPC} stability framework of \cite{grune_multiobjective_2019}, adapted to the proposed \acc{DM} schemes.
The resulting stability argument relies on suitable terminal ingredients and a descent condition that is less restrictive than the one used in \cite{grune_multiobjective_2019}.

The main contribution of this paper is an \emph{a priori} \acc{MOMPC} strategy for fixed point stabilization, in which scalarization-based \acc{DM} methods informed by the \accp{IM} select the \acc{ES} to be applied, and for which closed-loop asymptotic stability is guaranteed.
The strategy consists of two parts that can be considered independently of each other:
\begin{enumerate}
	\item
		\emph{Decision-making:} Six \acc{DM} methods, four established and two novel ones, are analyzed, which aim to improve the translation of a high-level preference vector to a point on the \acc{PF}, especially with respect to differing objective ranges and different shapes of the \acc{PF}.
		They require only two sequential optimization stages, which facilitates a real-time implementation of \acc{MOMPC}, and share a common preference vector that can be adapted online.
	\item
		\emph{Stability:} The \acc{DM} methods are integrated into the quasi-infinite horizon \acc{MOMPC} framework of \cite{grune_multiobjective_2019}.
		It is proven that, given suitable terminal ingredients, the closed-loop system is asymptotically stable for the \acc{DM} methods, using a descent condition that is less restrictive than in \cite{grune_multiobjective_2019} and thus retains a larger \acc{DM} space.
		The \accp{IM}-informed constructions only determine which \acc{ES} is selected, whereas the stability guarantee relies on the terminal ingredients and the descent condition.
\end{enumerate}

The remainder of this paper is structured as follows.
\Cref{sec:MOO} lays the theoretical groundwork by introducing \acc{MOO} problems and common scalarization techniques.
Building upon this, \cref{sec:IMIDM} addresses the first part, where the six \accp{IM}-informed \acc{DM} methods are developed and analyzed, together with metrics for their assessment.
The close interconnection between scalarization methods and \acc{DM} is shown, which allows for a systematic derivation of both well-known and novel methods.
\Cref{sec:IMI_MOMPC} then addresses the second part by embedding these \acc{DM} methods into the \acc{MOMPC} framework, providing the corresponding stability proof, and presenting a practical construction of the terminal ingredients.
The methods and results from \cref{sec:IMIDM} and \cref{sec:IMI_MOMPC} are then applied to a numerical case study in \cref{sec:num_ex}, where the resulting \acc{MO} schemes are compared with fixed-weight \acc[l]{WS} controllers to illustrate the potential benefit of adapting the preference online.
Finally, \cref{sec:conclusion} concludes the paper and outlines potential directions for future research.

\goodbreak
\textbf{Notation.} 
The vertical concatenation of two scalars, vectors, or matrices $X$ and $Y$ with matching column dimensions is denoted by $[X \concSep Y] := (X\tr, Y\tr)\tr$.  
Bold symbols $\boldsymbol{0}$, $\boldsymbol{1}$, and $\boldsymbol{\infty}$ denote vectors (or matrices) of zeros, ones, and infinities, respectively, with context-dependent dimensions.
For a vector $v \in \Rset^n$ (in)equalities like $v\geq0$, should be understood component-wise, i.e.~$v_i \geq 0$ for all $i \in \{1,\ldots,n\}$.
The following operators are used:
$\operatorname{sum}(v) := \sum_i v_i$, which computes the sum over all elements of a vector $v$; 
$\diag(v)$, which returns a diagonal matrix with $v$ on its diagonal; 
and the sign function $\operatorname{sign}(x)$, defined as $\operatorname{sign}(x) = 1$ if $x > 0$, $\operatorname{sign}(x) = 0$ if $x = 0$, and $\operatorname{sign}(x) = -1$ if $x < 0$.
The angle between two vectors $u$ and $v$ in degrees is computed by $\measuredangle(u,v)$.
The operation $\mathbb{A}\oplus\mathbb{B}$ denotes the Minkowski sum of the two sets $\mathbb{A}$ and $\mathbb{B}$.

\section{Multi-Objective Optimization}
\label{sec:MOO}
This section covers the theoretical basics of multi-objective optimization needed to develop the \acc{MOMPC} scheme, providing the necessary background to keep this work self-contained.
Hence, this section is mostly a review of the literature.

\subsection{Problem Statement}
\label{ssec:Intro_MOO}

Consider the \acc{MOOP} 
\begin{equation}
	\min_{\xMO \in \decSetFeas} \
	\costFunVec(\xMO)
	\tag{$\mathcal{P}$}
	\label{eq:MOOP}
\end{equation}
with the vector $\costFunVec(\xMO) := \matInline{\costFun_1(\xMO), \ldots, \costFun_\nObj(\xMO)}\tr$
of $n_J$ objectives $\costFun_i:\decSetFeas \rightarrow \mathbb{R}$, $i=1,\ldots,n_J$.
The set $\decSetFeas \subseteq \mathbb{R}^{n_{\xMO}}$ denotes the feasible set, i.e.~a point $x\in \mathbb{R}^{n_{\xMO}}$ is called feasible, if $x\in \decSetFeas$.
The meaning of the vector-valued minimization in \eqref{eq:MOOP} is clarified by the following definitions and conventions adopted from \cite{stieler_performance_2018}.
\goodbreak
\begin{bThm}[Efficiency, nondominance]{definition}
	\label{def:Nondom}
	A point $\xMO^{\star} \in \decSetFeas$ is called an \emph{\acc{ES}}
	to the \acc{MO} \acc[l]{OP} \eqref{eq:MOOP} if there is no feasible $\xMO \in \decSetFeas$ such that
	\begin{align}
		& \costFun_i(\xMO) \leq \costFun_i\left(\xMO^{\star}\right)
		\text{ for all }  i \in\{1, \ldots, \nObj\} 
		\text{ and }
		\\
		& \costFun_k(\xMO)<\costFun_k\left(\xMO^{\star}\right) 
		\text { for at least one } k \in\{1, \ldots, \nObj\}.
	\end{align}
	The respective image value 
	$
	\costFunVec\left(\xMO^{\star}\right)
	$ 
	is called \emph{nondominated}.
	The set of all nondominated points is called the \emph{nondominated set} or \emph{\acc[l]{PF}}.
	\\
	A point $\xMO^{\star} \in \decSetFeas$ is called a \emph{weakly \acc{ES}} to the \acc{MO} \acc[l]{OP} \eqref{eq:MOOP} if there is no feasible $\xMO \in \decSetFeas$ such that $\costFun_i(\xMO)<\costFun_i\left(\xMO^{\star}\right)$ holds for all $i \in\{1, \ldots, \nObj\}$. 
	The image value of a weakly \acc{ES} is called \emph{weakly nondominated}.
\end{bThm}
\begin{bThm}{convention} 
	Throughout the paper, the $\min$- and $\argmin$-operator  in the context of \acc{MO} optimization are defined as
	\begin{align}
		\argmin_{\xMO \in \decSetFeas} \,
		\costFunVec(\xMO)	
		&=: \decSetPareto
		,
		\qquad
		\min_{\xMO \in \decSetFeas} \,
		\costFunVec(\xMO)
		=: \imSetPareto
		.
	\end{align}
\end{bThm}
Equivalently, the set of \accp{ES}, also called efficient set or Pareto set, is defined by
$\decSetPareto:= \left\{\xMO \in \decSetFeas \ \middle| \ \xMO \text{ is an \acc{ES} to \eqref{eq:MOOP}} \right\}$
and the \acc{PF} as
$\imSetPareto:=\left\{\costFunVec(\xMO)  \ \middle| \ \xMO \in \decSetPareto\right\}$.
Furthermore, the set of feasible image points of the objective is defined by
\begin{equation}
	\imSetFeas:=\left\{\costFunVec(\xMO)  \ \middle| \ \xMO \in \decSetFeas\right\}.
	\label{eq:feasImage}
\end{equation}
Now, in \acc{MOO}, the set $\imSetPareto$, the \acc{PF}, typically consists of infinitely many points.
Therefore, one must address the task of how to select a single point from this set.
This task is called \acc[f]{DM}.

\subsection{Characteristic Quantities}
\label{ssec:Characteristic_Quantities}
In the following, characteristic quantities of \acc{MOOP} are derived, which are key to the development of the \acc{DM} methods presented below.

The \accp[f]{IM} are denoted by $\xMO_{i}\opt$, $i = 1,\ldots,n_J$, i.e.~the solutions to the single-objective \acc{OP}
\begin{equation}
	\min_{\xMO \in \decSetFeas} \ J_i(\xMO), \quad i=1,\dots,n_J.
	\label{eq:P_IM_i}
\end{equation}
Evaluating the full objective vector $\JMO$ at all \acc{IM} defines the so-called pay-off matrix (see, e.g.~\cite{ehrgott_multicriteria_2005})
\begin{equation}
	\Phi =
	\begin{mat}
		\JMO\big(\xMO_{1}\opt\big), &\ldots, &\JMO\big(\xMO_{n_J}\opt\big)
	\end{mat}
	.
	\label{eq:Phi}
\end{equation}
To derive further quantities based on the pay-off matrix $\Phi$, the following two functions are required:%
\begin{bThm}{definition}
	\label{defi:scaling_normal}
	For $v \in \Rset^n \setminus \{0\}$, let $i_{\max}(v) := \min\big(\arg\max_i |v_i|\big)$ be the smallest index of a component with maximal absolute value.
	The scaling function $\fScal: \Rset^n \setminus \{0\} \to \Rset^n$ is defined by
	\begin{equation}
		\fScal(v) := -\operatorname{sign}\big(v_{i_{\max}(v)}\big) \, \frac{v}{\|v\|_1}.
	\end{equation}
	Hence, $\fScal(v)$ is the multiple of $v$ with $\|\fScal(v)\|_1 = 1$ whose component of largest absolute value is negative.
	It is used to fix the length and orientation of normal vectors, which are otherwise only determined up to a non-zero factor.
	In particular, $\fScal(cv) = \fScal(v)$ for all $c \in \Rset \setminus \{0\}$.
\end{bThm}
\begin{bThm}{definition}
	\label{defi:normal_vector}
	Let the columns of $M = \matInline{m_1,\ldots,m_n} \in \Rset^{n \times n}$ be $n$ affinely independent points.
	The hyperplane through these points is defined by
		$
		\mathcal{H}(M) := \left\{ M\gamma \ \middle| \ \gamma \in \Rset^{n}, \ \operatorname{sum}(\gamma) = 1 \right\}
		$
		.
	A vector $w \in \Rset^n \setminus \{0\}$ with $w\tr (m_i - m_1) = 0$ for $i = 2,\ldots,n$ is called a normal vector of $\mathcal{H}(M)$.
	The scaled normal vector is defined by
	\begin{equation}
		\fNormal(M) := \fScal(w),
	\end{equation}
	where $w$ is any normal vector of $\mathcal{H}(M)$.
\end{bThm}
Since the $n$ points are affinely independent, $\mathcal{H}(M)$ is a unique $(n-1)$-dimensional hyperplane in $\Rset^n$, and its normal vectors are unique up to a non-zero factor.
As $\fScal(cv) = \fScal(v)$ for all $c \neq 0$, $\fNormal(M)$ does not depend on the choice of $w$.
The quantities derived from $\Phi$ are collected in \cref{list:objects_from_Phi}
and visualized in \cref{fig:quants_from_Phi}.
The front shown there is not a sketch, but the \acc{PF} of the \acc{MOOP} introduced in \cref{example2}.
\goodbreak
\begin{NonFloatList}[Quantities derived from $\Phi$]
	\begin{enumerate}[itemsep=0em,parsep=0em]
		\item 
			\label{item:quants_from_Phi_1}
			\textbf{Nadir point.}
			The (\textit{practical}) \acc{NP} $\J_\NP$ is defined as the row-wise maximum of $\Phi$%
			\footnote{%
				Although, generally, for $n_J > 2$, the actual or real \acc{NP} may be larger than $\J_{\NP}$, this definition provides an easy-to-compute approximation. According to \cite[Chapter 2.2]{ehrgott_multicriteria_2005}, finding the \textit{real} \acc{NP} is a \enquote{very difficult problem}, for which \enquote{no efficient method} is known to address arbitrary problems. Hence, \enquote{heuristics are often used}, e.g.\ the estimation via \enquote{pay-off tables} such as the matrix $\Phi$ introduced above.%
			}.
		\item \textbf{Utopia point.}
		The \acc{UP} $\J_\UP$ is defined as the row-wise minimum of $\Phi$.
		\item \textbf{\Acc[s]{CHIM}.}
		As in \cite{das_normal-boundary_1998}, the \acc{CHIM} is introduced as the simplex formed by the columns of the pay-off matrix $\Phi$, i.e.
		$\mathbb{P} := \{\Phi\baryCO \mid \baryCO\in\FunitSimplex\}$,
		where $\FunitSimplex:= \{\gamma \in [0,1]^{n_J} \mid \operatorname{sum}(\gamma)=1\}$ denotes the unit simplex.
		$\Phi$ is assumed to be full rank in the following.
		\item \textbf{Normal vector.}
		The \acc{CHIM} is contained in the $(n_J-1)$-dimensional hyperplane $\mathcal{H}(\Phi)$ (cf.~\cref{defi:normal_vector}).
		Its scaled normal vector is denoted by $\wChim := \fNormal(\Phi)$, such that all $\chi \in \mathcal{H}(\Phi)$ satisfy $\wChim\tr \chi +\bChim=0$ for some $\bChim \in \Rset$.
		\item 
			\label{item:quants_from_Phi_5}
			\textbf{Quasi-normal vector.}
			A quasi-normal vector of the \acc{CHIM}, $\wChimQuasi$, is the (scaled) vector that connects a point on the \acc{CHIM} with the \acc{UP} (cf.~\cite{das_normal-boundary_1998}), 
			i.e.
			\begin{equation}
				\label{eq:wChimQuasi}
				\wChimQuasi = \fScal\(- (\Phi \baryCO - \J_{\UP})\),
			\end{equation}
			with $\baryCO \in \FunitSimplex$, such that $\Phi \baryCO$ lies on the \acc{CHIM}.
			Note that $\wChimQuasi$, by definition, only has non-positive components.
		\item 
			\label{item:quants_from_Phi_6}
			\textbf{Barycentric center.}
			A special point on the \acc{CHIM} is the barycentric center that is obtained via $\Phi \baryCO_{\mathrm{c}}$ with $\baryCO_{\mathrm{c}} = 1/n_J \cdot \boldsymbol{1}$.
			The quasi-normal corresponding to $\baryCO_{\mathrm{c}}$ is denoted by $\wChimQuasi_{\mathrm{c}}$.
		\item \textbf{Normalization.}
		Let $\cNPUP$ denote the positive definite, diagonal matrix
  		\begin{equation}
			\cNPUP:= 
			\diag(\J_{\NP}-\J_{\UP})\inv.
			\label{eq:cNPUP}
		\end{equation}
		Then, the \acc{UP} and \acc{NP} can be scaled, such that they form a hypercube or unit box of edge length one, i.e.~$\cNPUP (\J_{\NP}-\J_{\UP})=\boldsymbol{1}$.
		Using $\cNPUP \J$ as the objective in problem \eqref{eq:MOOP} then generates the \emph{\acc{NIS}} \cite{rao_modified_1991}. 
		$\cNPUP\J$ is called a normalized objective vector.
		\item \label{item:eta_vis}
		\textbf{Visual normal vector.}
		The \emph{visual} \acc{CHIM} normal $\wChimVisual$ is introduced.
		In two or three dimensions, this vector appears to be perpendicular to the \acc{CHIM} hyperplane if the box spanned by \acc{NP} and \acc{UP} appears to be a cube.
		This effect can be seen in \cref{fig:quants_from_Phi}, where the \emph{true} normal $\wChim$ appears to be quite far from being perpendicular to the \acc{CHIM} hyperplane, whereas the \emph{visual} normal $\wChimVisual$ does not.
		Hence, $\wChimVisual$ is a direction rather than a normal vector in the original image space; it is normal to the \acc{CHIM} only in the \acc{NIS}.
		It can be computed by computing a normal in normalized image space and transforming it back to the original image space:
		\begin{equation}
			\begin{aligned}
				\wChimVisual :=& \
				\fScal \(\cNPUP\inv \ \fNormal \(\cNPUP\Phi\)\)
				.
			\end{aligned}
			\label{eq:eta_visual}
		\end{equation}
	\end{enumerate}
	\label{list:objects_from_Phi}
\end{NonFloatList}
\begin{figure}[h]
	\centering
		\adjincludegraphics[
			width=0.95\WidthBase, clip, 
			trim={{0.0\width} {0.01\height} {0.0\width} {0.0\height}}
		]{quants_from_Phi.\FigureFileExt}
	\caption[%
		Convex Pareto front with some of the quantities mentioned in \cref{list:objects_from_Phi}.
	]
	{%
		Convex Pareto front with
		some of the
		quantities mentioned in \cref{list:objects_from_Phi}.
		The vectors
			\if\TwoCol1
				\par\vspace{4pt}
				{%
					\centering%
					$
					\begin{aligned}
						\wChim       &\approx -\matInline{0.979, -0.011, 0.009}\tr, \\
						\wChimQuasi  &\approx -\matInline{0.011,  0.074, 0.914}\tr \qquad \text{and} \\
						\wChimVisual &\approx -\matInline{0.010, -0.004, 0.986}\tr 	
					\end{aligned}
					$
					\par
				}
				\vspace{4pt}\noindent
			\else
				$\wChim       \approx -\matInline{0.979, -0.011, 0.009}\tr$,
				$\wChimQuasi  \approx -\matInline{0.011,  0.074, 0.914}\tr$ and 
				$\wChimVisual \approx -\matInline{0.010, -0.004, 0.986}\tr$ 
			\fi
		are scaled according to \cref{defi:scaling_normal}.
		Note that all components of $\wChimQuasi$ are negative, whereas $\wChim$ and $\wChimVisual$ have a positive second component.
	}
	\label{fig:quants_from_Phi}
\end{figure}
Let $e_i \in \Rset^{n_J}$ denote the $i$-th unit vector.
In addition to $\Phi$, it can be of interest to compute the \emph{regularized} pay-off matrix $\Phi_{\mathrm{reg}}$ and the quantities derived from it (cf.~\cref{list:objects_from_Phi}), which are used below to characterize the \acc{PF} near the \accp{IM}.
This can be done by computing the \emph{regularized} or \emph{non-extreme} \accp{IM}%
\footnote{
	The following \acc{OP} is a \acc[l]{WS} scalarization of \eqref{eq:MOOP}; 
	this method is formally introduced later in \cref{sssec:WS}.
}
\makeatletter%
	\AC@reset{WS}%
\makeatother%
\begin{subequations}
\begin{align}
	\xMO_{i,\mathrm{reg}}\opt
	&=
	\argmin_{\xMO \in \decSetFeas} \ 
	\tilde{e}_i\tr \, f_{\mathrm{NIS}}\(\J(\xMO)\)
	\quad
	i = 1,\ldots,n_J,
	\notag \\ & 
	\label{eq:NEIM}
\intertext{with}
	\tilde{e}_i
	&:= e_i + w_{\mathrm{reg},i}
	,
	\label{eq:unit_vec_reg}
	\\
	w_{\mathrm{reg},i} 
	&:= 
	\diag(\boldsymbol{1}-e_i) \boldsymbol{\delta} - e_i(\boldsymbol{1}-e_i)\tr\boldsymbol{\delta}
	,
	\label{eq:w_reg_i}
	\\
	\boldsymbol{\delta} &\in 
	\Big\{
		\tilde{\boldsymbol{\delta}} \in \Rset_{\geq 0}^{n_J}
	\, \Big| \,
		\textstyle
		\sum_{k \neq \, i} \tilde{\boldsymbol{\delta}}_k < 1
		\ 
		\if\TwoCol1
			\\&\hspace{0.37\linewidth}
		\fi
		\forall  \, i \in \{1,\ldots,n_J\}
	\Big\}
	,
	\label{eq:bs_delta_restriction}
	\\
	f_{\mathrm{NIS}}(\J_{0}) &:= \cNPUP(\J_{0}-\J_{\UP})
	.
	\label{eq:f_NIS}
	\numberthis
\end{align}
Thus, $\tilde{e}_i$ has the entries $\boldsymbol{\delta}_k$ for $k \neq i$ and $1-\sum_{k \neq i}\boldsymbol{\delta}_k$ for $k = i$, i.e.~$\operatorname{sum}(\tilde{e}_i)=1$.
From $\boldsymbol{\delta}>0$ and \eqref{eq:bs_delta_restriction}, it follows that $\tilde{e}_i > 0$, which guarantees the solution of \eqref{eq:NEIM} to be efficient \cite[Theorem~3.1.2]{miettinen_nonlinear_1998}.
Note that the normalization via $f_{\mathrm{NIS}}$ is essential to compensate for potentially vastly different objective ranges. 
Without such a normalization, the regularization may, in extreme cases, have either a negligible effect or lead to an excessively strong distortion.
The normalization via $f_{\mathrm{NIS}}$ also simplifies the choice of a suitable regularization magnitude, since comparable values of $\boldsymbol{\delta}$ then induce comparable effects across objectives.
In practice, the isotropic choice $\boldsymbol{\delta} = \delta \cdot \boldsymbol{1}$ with $\delta \in [10^{-3},10^{-1}]$ is therefore used as a guideline in the normalized setting.
Smaller values keep the regularized \accp{IM} close to the original \accp{IM}, whereas larger values move them more strongly toward non-extreme trade-off points.
Naturally, $\delta$ must still satisfy \eqref{eq:bs_delta_restriction}.
As can be seen from \eqref{eq:w_reg_i}, the regularization explicitly depends on $\boldsymbol{\delta}$, and consequently 
$\Phi_{\mathrm{reg}}$
becomes a function of $\boldsymbol{\delta}$,
\begin{equation}
	\Phi_{\mathrm{reg}}(\boldsymbol{\delta})
	= 
	\begin{mat}
		\JMO\big(\xMO_{1,\mathrm{reg}}\opt(\boldsymbol{\delta})\big), 
		&\ldots, 
		&\JMO\big(\xMO_{n_J,\mathrm{reg}}\opt(\boldsymbol{\delta})\big)
	\end{mat}.
	\numberthis
\end{equation}
\label{eq:Phi_reg}%
\end{subequations}
This derivation gives rise to the following definition, which further characterizes certain properties of the \acc{PF}.
\begin{bThm}[Pathological Pareto front]{definition}
	\label{defi:pathological_PF}
	Assume that the regularized pay-off matrix $\Phi_{\mathrm{reg}}(\boldsymbol{\delta})$ (cf.\eqref{eq:Phi_reg}) was computed with 
	$\boldsymbol{\delta}=\delta \cdot \boldsymbol{1}$.
	The corresponding \acc{NP} and \acc{UP} are denoted by $\J_{\NP,\mathrm{reg}}(\delta)$ and $\J_{\UP,\mathrm{reg}}(\delta)$.
	The volume of the hyperbox spanned by them in the \acc{NIS} is
	\begin{equation}
		V_{\mathrm{box}}(\delta)
			:= \prod_{i=1}^{n_J} \Big[ \cNPUP\big(\J_{\NP,\mathrm{reg}}(\delta)-\J_{\UP,\mathrm{reg}}(\delta)\big) \Big]_i .
	\end{equation}
	Note that $\J_{\NP,\mathrm{reg}}(0)=\J_{\NP}$, $\J_{\UP,\mathrm{reg}}(0)=\J_{\UP}$, and hence $V_{\mathrm{box}}(0) = 1$.
	\\
	Now, a \acc{PF} is defined to be \emph{pathological} if, under a uniform increase of the regularization parameter $\delta$, the ratio $V_{\mathrm{box}}(\delta) / V_{\mathrm{box}}(0)$ exhibits an approximately exponential decay.
	Visually speaking, the hyperbox spanned by $\J_{\NP,\mathrm{reg}}(\delta)$, $\J_{\UP,\mathrm{reg}}(\delta)$ contracts rapidly.
	This is an indication that the \acc{PF} is extremely steep or extremely flat near the non-regularized \accp{IM}.
\end{bThm}
\begin{bThm}{remark}
	\label{rem:pathological_PF_rule_of_thumb}
	A quick practical assessment of whether a \acc{PF} is pathological can be carried out by computing $\Phi_{\mathrm{reg}}(\delta)$ using $n_J$ additional optimizations.
    Choose a single small $\delta>0$ and compute the contraction of the hyperbox. 
    As a rule of thumb, a \acc{PF} is considered pathological if 
    $V_{\mathrm{box}}(\delta) / V_{\mathrm{box}}(0) < 50\%$ for $\delta = 0.005$.
\end{bThm}

\subsection{Sampling Points on the Pareto Front}
In \cref{ssec:Intro_MOO} it is shown that the solution of a \acc{MOOP} of type \eqref{eq:MOOP} is a set $\imSetPareto$, the \acc{PF}. 
For most applications involving conflicting objectives, this set contains an infinite number of points.
In order to approximate the \acc{PF}, one can apply a scalarization method \cite[Chapter 2]{eichfelder_adaptive_2008}.
A scalarization method transforms \eqref{eq:MOOP} to a \emph{parameterized} single-objective \acc{OP}.
By varying the parameter(s), different points on the \acc{PF} are obtained.
This process is also called \emph{sampling the \acc[l]{PF}}.
Two scalarization methods are detailed in the following.

\subsubsection{Weighted-Sum Method}
\label{sssec:WS}
Arguably the most popular scalarization method is \acc{WS}.
Applying the \acc{WS} method to \eqref{eq:MOOP} using the parameter $\param{w}_{\mathrm{WS}} \in \FunitSimplex$ yields the parameterized single-objective \acc{OP}
\begin{equation}
	\begin{aligned}
		\min_{\xMO \in \decSetFeas} \ &\param{w}_{\mathrm{WS}}\tr \J(\xMO).
	\end{aligned}
	\tag{$P_{\mathrm{WS}}(\param{w}_{\mathrm{WS}})$}
	\label{eq:SM_WS}
\end{equation}
The \paramAccent accent is used here and throughout the remainder of the paper to indicate that a quantity is a parameter.
Choosing non-negative weight components guarantees obtaining weakly \accp{ES} \cite[Theorem 3.1.1.]{miettinen_nonlinear_1998}.
However, this method cannot sample points in non-convex regions of the \acc{PF} \cite[Chapter 2.5.7]{eichfelder_adaptive_2008}.

\subsubsection{Method by Pascoletti and Serafini}
\label{ssec:PS}
Similarly to \eqref{eq:SM_WS}, the scalarization of \eqref{eq:MOOP} using the \acc{PS} method 
\cite[Chapter 2.1]{eichfelder_adaptive_2008}, \cite{pascoletti_scalarizing_1984} 
yields
\begin{equation}
	\begin{aligned}
		\min_{\xMO \in \decSetFeas,\, l \in \mathbb{R}}
		\ & - l
		\\
		\text{s.t.} \quad 
		& \param{\J}_{\SO} + l\param{d} - \J(\xMO) \in \conePS[\param].
	\end{aligned}
	\tag{$P_{\mathrm{PS}}(\param{\J}_{\SO}, \param{d}, \conePS[\param])$}
	\label{eq:SM_PS}
\end{equation}
The quantity $\param{\J}_{\SO}$ is called \acc{SO} and $\param{d}$ \acc{SDV}.
Furthermore, $\conePS[\param]$ is a nonempty closed pointed convex cone.
This method can be interpreted as shooting a ray from $\param{\J}_{\SO}$ in the direction $\param{d}$, where $l$ is used to regulate the length of the ray.
For admissibility, $\J(\xMO)$ must remain within the cone originating at the ray's endpoint.

To improve the readability of the paper, a strict notation convention is adhered to.
While the \textit{tilde} ($\tilde{\cdot}$) is used flexibly for general modifications and transformations, the following two accents are reserved for specific contexts: 
the \textit{\paramAccent} ($\param{\cdot}$) for parameter variables
and
the \textit{\realizeAccent} ($\realize{\cdot}$) for their specific realizations.
\Cref{tab:notation_accents} provides a summary of these conventions.
\begin{table}[ht]
    \centering
    \tableFontSizeWiley
    \begin{tblr}{colspec={r|cl}, rowsep = 0.1em, row{2} = {abovesep = 0.4em}}
        \toprule
        Accent & Example & Reserved meaning / use case \\
        \midrule
        \textit{\paramAccent} & $\param{d}$ & parameter of an \acc{OP} \\
        \textit{\realizeAccent} & $\realize{d}$ & numerical value of a parameter \\ 
        \textit{tilde}  & $\tilde{\Phi}$ & general modifications of a specific quantity \\
        \bottomrule
    \end{tblr}
    \caption{Summary of mathematical notation conventions and reserved accents.}%
    \label{tab:notation_accents}%
\end{table}

\paragraph*{Normal Boundary Intersection}
In \cite{das_normal-boundary_1998}, \acc[f]{NBI}, was introduced and further discussed in \cite{das_characterizing_1999} in the context of knee-points.
Although the naming of the method suggests otherwise, the \acc{SDV} is not necessarily set to the normal vector of the \acc{CHIM}, but to $\realize{d}=-\tilde{\Phi}\boldsymbol{1}$ with the \textit{shifted} pay-off matrix%
\footnote{%
	Note that in \cite{das_normal-boundary_1998,das_characterizing_1999} the \textit{shifted} pay-off matrix is denoted by $\Phi$, i.e.~no tilde accent is used.
}%
\begin{equation}
	\tilde{\Phi}=\Phi - (\J_{\UP},\ldots,\J_{\UP}).
	\label{eq:Phi_shifted}
\end{equation}
The method parameterizes the \acc{SO} as $\realize{\J}_{\SO}(\baryCO) = \tilde{\Phi}\baryCO + \J_{\UP}$ 
and sets $\conePS=\{\boldsymbol{0}\}$.
The parameter $\baryCO$ is used to perform a convex combination of the (shifted) \accp{IM} which means that if $\baryCO\in \FunitSimplex$, then $\Phi\baryCO$ yields a point on the \acc{CHIM}.
\\
Equivalently, this can be represented as
\begin{equation}
	\realize{\J}_{\SO}(\baryCO) = \Phi\baryCO
	\ \text{and} \
	\realize{d}(\baryCO_{\mathrm{d}})=-n_J\cdot(\Phi \baryCO_{\mathrm{d}} - \J_{\UP}),
	\
	\baryCO_{\mathrm{d}} = \baryCO_{\mathrm{c}}
\end{equation}
(recall that $\Phi \baryCO_{\mathrm{c}}$ gives the barycentric center)
which has a more intuitive interpretation: 
The \acc{SO} is a convex combination of the \accp{IM} and the \acc{SDV} points from the barycentric center of the \acc{CHIM}%
\footnote{Other feasible $\baryCO \in \FunitSimplex$ are equally valid choices.}
to the \acc{UP}, i.e.~it is the quasi-normal vector.
\\
Choosing $\conePS=\Rset_{\geq 0}^{n_J}$, where $\Rset_{\geq 0}^{n_J}$ is the set of all non-negative real numbers in $\Rset^{n_J}$, enables the determination of (weakly) \accp{ES} even when shooting towards dominated points.
Then, the constraint in \eqref{eq:SM_PS} is replaced by
\begin{equation}
	0 \leq \param{\J}_{\SO} + l\param{d} - \J(\xMO).
	\label{eq:PS_constraint}
\end{equation}
An illustration of the method is given in \cref{fig:NBI_rays_coverage}.

\begin{figure*}[htb]
	\centering
	\begin{subfigure}[t]{0.49\textwidth}
		\centering
				\adjincludegraphics[width=\linewidth,trim={{0.0\width} 0 {0.0\width} {0.0\height}},clip]{NBI_rays_coverage.\FigureFileExt}
		\caption{Points obtained using the \acc{NBI} method.}
		\label{fig:NBI_rays_coverage}
	\end{subfigure}
	\hfil
	\begin{subfigure}[t]{0.49\textwidth}
		\centering
			\adjincludegraphics[width=\linewidth,trim={{0.0\width} 0 {0.0\width} {0.0\height}},clip]{SRI_rays_coverage.\FigureFileExt}
		\caption{Points obtained using the \acc{SRI} method.}
		\label{fig:SRI_rays_coverage}
	\end{subfigure}
	\caption{
		A \acc{PF} with $\J_{\NP}=\matInline{0,0,0}\tr$ and $\J_{\UP}=-\matInline{1,1,1}\tr$.
		The \acc{NBI} method samples only parts of the \accp{PF}, 
		whereas the \acc{SRI} sampling covers the entire front.
	}
	\label{fig:mainfig}
\end{figure*}

\paragraph*{Spherical Ray Intersection}
The second subform, called \acc[f]{SRI}, was introduced in \cite{herrmann-wicklmayr_accelerating_2025}.
The main idea of this method is to shoot rays from a point, which can be imagined to be the center of a (stretched) hypersphere, in the direction of points on (some parts of) the surface of the hypersphere. 
An illustration is given in \cref{fig:SRI_rays_coverage}.
Here, the \acc{SDV} is given by $\realize{d}(\tau) = \cNPUP\inv \, x_{\mathrm{seg}}(\tau)$, where $x_{\mathrm{seg}}(\tau) \leq 0$, $\tau \in [0,1]^{n_J-1}$, parameterizes a hypersphere segment in the nonpositive orthant \cite{herrmann-wicklmayr_accelerating_2025}.
The \acc{SRI} method uses the cone $\conePS = \Rset_{\geq 0}^{n_J}$, which yields the same constraint as in \eqref{eq:PS_constraint}.

A detailed comparison of the \acc{NBI} and \acc{SRI} method is given in \cref{sec:Comp_NBI_SRI}.

\section{Individual Minima-Informed Decision-Making}
\label{sec:IMIDM}
While knowledge of entire \accp{PF} is advantageous in design and process optimization, in the end, only a single solution can be produced or implemented.
This means, for any application, exactly one nondominated solution has to be chosen---a process called \acc[f]{DM}.
Human expert decision-makers typically chooses based on a visual examination of the \acc{PF}. 
This already becomes practically infeasible if the number of objectives is larger than three.
\\
Furthermore, in the context of \acc{MOMPC}, where a \acc{MOOP} is repeatedly solved at controller frequency, such a procedure would not be feasible.
Hence, an automated \acc{DM} strategy is needed. 
Due to the ever-present challenge of real-time capability of \acc{MPC}, the number of optimization problems to be solved should be as small as possible.

\subsection{Specifying Relative Importance of Objectives}
\label{ssec:Specifying_Relative_Importance}
One way to keep the number of optimizations small 
while yet robustly translate a-priorily given, high-level preferences to a single solution
is to compute a few characteristic Pareto optimal points as a basis for \acc{DM}.
The information encoded in the \accp{IM} is proposed to be used in automated \acc{DM} methods to achieve such a translation.
As a prerequisite, $n_J$ \accp{OP} have to be solved to compute all \accp{IM}.

Throughout this paper, the term \emph{\acc[l]{DM}} refers to the algorithmic determination of an \acc{ES} through a scalarization method parameterized by a preference vector $\betaPref \in \FunitSimplex \subset \Rset^{n_J}$, which is intended to encode the relative importance of the objectives.
Here, $\betaPref$ serves as a common parameter through which the considered methods are addressed.
Its technical role depends on the selected scalarization.
For \acc{WS}-based methods, it influences the weight vector, whereas for \acc{NBI}-type methods it acts as a barycentric coordinate on the \acc{CHIM} and thereby determines the shooting origin $\realize{\J}_{\SO}=\Phi\betaPref$.
For the \acc{NCHIM} method, it determines the point $\Phi\betaPref$ on the \acc{CHIM} through which the ray starting at the \acc{NP} is shot.
Accordingly, $\betaPref$ should not be interpreted as a direct objective weighting in general, but as a parameter through which the selected scalarization maps a high-level preference to a point on the \acc{PF}.
How such a translation can be assessed in a systematic and method-independent way is discussed later in this section.
Moreover, the resulting translation \form{need not be equally meaningful} for arbitrary \acc{PF} geometries; 
in particular, pathological \accp{PF} (cf.~\cref{defi:pathological_PF}) may lead to a deterioration of the translation quality, as discussed later in \cref{ssec:Challenges_of_Problem}.

Leveraging the information encoded in the \accp{IM}, such a scalarization-based \acc{DM} criterion can then be realized by applying a suitable scalarization of \eqref{eq:MOOP} with appropriately chosen parameters.
Thus, only one additional optimization is needed to determine the nondominated solution to be realized.
Hence, in total, $n_J+1$ optimizations have to be performed for \accp{IM}-informed \acc{DM}.
If the computation of the \accp{IM} is done in parallel, then only two sequential (batches of) optimizations have to be executed.
Otherwise, the number of sequential optimizations scales linearly with the number of objectives $n_J$.
Both approaches increase the likelihood of \acc{MOMPC} having real-time capabilities.

The presumably best-\emph{known} parameterization of a scalarized \acc{MOOP} with the goal of implementing a \acc{DM} criterion is the \acc{WS} method (cf.~\cref{sssec:WS}).
One might try to set the preference \emph{directly} by setting $\param{w}_{\mathrm{WS}}=\betaPref$.
However, the authors of \cite{marler_weighted_2010} \enquote{investigate the fundamental significance of the weights in terms of preferences}
and conclude that expressing preferences via \acc{WS} weights is difficult, because the objectives' ranges and magnitudes are not considered.
Moreover, in the context of \acc{MPC}, the objectives' ranges and magnitudes might be time-varying.
In \cite{marler_weighted_2010}, it is proposed to transform \enquote{all objective functions [...] such that they have similar ranges}. 
Equivalently, the weight vector can be transformed using a positive scaling vector $c_{\J}>0$:
\locDefPenalites{$
	\realize{w}_{\mathrm{WS}}\tr \diag(c_{\J}) \J
	= (\diag(c_{\J})\tr \realize{w}_{\mathrm{WS}})\tr \J
	= (\diag(c_{\J}) \realize{w}_{\mathrm{WS}})\tr \J
	= (\realize{w}_{\mathrm{WS}}^{\prime})\tr \J
$}
with $\realize{w}_{\mathrm{WS}}^{\prime} := \diag(c_{\J}) \realize{w}_{\mathrm{WS}}$.

\subsection{Novel Scalarization Methods Inspired by Decision-Making}
\label{ssec:Novel_SM_DM}
In the following, two (\acc{IM})-informed \acc{DM} methods are proposed, which have not previously been described in the literature, to the best of the authors' knowledge.
These novel methods aim to produce a more robust translation of preferences with respect to objective ranges than that achieved by established methods.
Both methods structurally belong to the \acc{PS} method.
Hence, their computational complexity is similar to that of other \acc{PS} subforms and the \acc{WS} method, which itself requires one optimization variable and $n_J$ inequality constraints less than \acc{PS} subforms.

\paragraph*{NBI Using the Visual CHIM Normal}
This method can be viewed as a subform of \acc{NBI} with $\realize{d}=\wChimVisual$. 
A closer examination of \eqref{eq:eta_visual} reveals that $\wChimVisual$ is the normal vector of the \acc{CHIM} in a \acc{NIS} transformed back into the original image space.
Its effect on robustness with respect to differing objective ranges is assessed numerically in \cref{ssec:Exemplary_Problems,fig:DM_all_methods,tab:metrics_DM}.%
\hyperanchor{anchor:NBI_vis}{}%

\paragraph*{Nadir-CHIM Method}
The \acc{NCHIM} constitutes a novel realization of the \acc{PS} method, given by the realization of the parameters
\begin{equation}
	\realize{\J}_{\SO} = \J_{\NP}, \quad
	\realize{d} = -(\J_{\NP}-\Phi\baryCO), \quad
	\conePS = \Rset_{\geq 0}^{n_J},
\end{equation}
where $\baryCO \in \FunitSimplex$.
It can be interpreted as a perspective projection from points on the CHIM onto the \acc{PF} using projection rays originating from $\J_{\NP}$.
Here, the use of points on the \acc{CHIM} provides robustness against different objective ranges.
Furthermore, the \acc{NCHIM} method can be viewed as a restriction of the \acc{SRI} method to those rays which shoot through the \acc{CHIM}, because this allows to replace the non-intuitive parameterization via $\tau$ by the $\Phi \baryCO$ sampling of the \acc{CHIM}, as it was derived for the \acc{NBI}-based \acc{DM}.

\subsection{Decision-Making Methods}
\label{ssec:DM_methods}
Six scalarization-based \acc{DM} methods are compared, summarized in \cref{list:DM_methods}.
All are based on the \acc[f]{WS} \eqref{eq:SM_WS} or \acc[f]{PS} \eqref{eq:SM_PS} scalarization. 
For \acc{PS}, the cone is fixed to $\conePS=\Rset_{\geq 0}^{n_J}$.
Consequently, the preference $\betaPref$ only influences the specific parameters of the chosen scalarization:
$\param{w}_{\mathrm{WS}}$ (for \acc{WS}) or $\param{\J}_{\SO}, \param{d}$ (for \acc{PS}).
\goodbreak
\begin{NonFloatList}[\acc{DM} methods (novel ones are highlighted in bold)]
	\begin{enumerate}[label={WS(\arabic*)},parsep=0em,itemsep=0.4\baselineskip,align=left,
	labelwidth=\widthof{WS(3)\hspace{0mm}},leftmargin=!]
		\item \label{enum:WS_1} 
			\acc{WS} with the preference vector, scaled using \accp{IM} information, as weight vector: $\realize{w}_{\mathrm{WS}} =
			\cNPUP \betaPref$
		
		\item \label{enum:WS_2} 
			\acc{KP}: $\realize{w}_{\mathrm{WS}} = -\wChim$; note that this approach only yields a single point
	\end{enumerate}
	\vspace{-0.3\baselineskip}%
	\begin{enumerate}[label={PS(\arabic*)},itemsep=0.4\baselineskip,align=left,
	labelwidth=\widthof{WS(3)\hspace{0mm}},leftmargin=!]
		\item \label{enum:PS_1} 
			\acc{NBI} with 
			$\realize{\J}_{\SO}=\Phi\betaPref$,
			$\realize{d}=\wChim$
		
		\item \label{enum:PS_2} 
			\acc{NBI} with
			$\realize{\J}_{\SO}=\Phi\betaPref$,	
			$\realize{d}=\wChimQuasi_{\mathrm{c}}$ (cf.~\cref{item:quants_from_Phi_6} in \cref{list:objects_from_Phi})
			
		\refstepcounter{enumi}
		\item[\textbf{PS(3)}] \label{enum:PS_3}
			\acc{NBI} with 
			$\realize{\J}_{\SO}=\Phi\betaPref$,
			$\realize{d}=\wChimVisual$
			
		\refstepcounter{enumi}
		\item[\textbf{PS(4)}] \label{enum:PS_4}
			\acc{NCHIM}: $\realize{\J}_{\SO}=\J_{\NP}$ 
			~~and
			\newline
			$\realize{d}=-(\J_{\NP}-\Phi\betaPref)$ (scaling is not required)
	\end{enumerate}
	\label{list:DM_methods}
\end{NonFloatList}
\Cref{list:DM_methods} contains the method \ref{enum:WS_2} that makes use of the well-known \acc{KP} concept (cf.~\cite{das_characterizing_1999}), which is discussed further in \cref{sec:KP}.
\\
All variants of the \acc{NBI} method in \cref{list:DM_methods} set the \acc[l]{SO} to $\realize{\J}_{\SO}=\Phi\betaPref$.
Since there is no \emph{unique} choice to set $\realize{d}$, two common choices (cf.~\cite{das_normal-boundary_1998}) and one novel variant (cf.~\ref{enum:PS_3}) are compared.

\subsection{Unified Problem Formulation}
\label{ssec:IMG_Prob_Formulation}
So far, the \accp{IM} required for all \acc{DM} methods shown in \cref{ssec:DM_methods} were assumed to be known.
In general, and specifically for \acc{MOMPC}, this is not the case. 
For the \acc{WS} and \acc{PS} scalarization, a unified \acc{OP} is formulated that, depending on the choice of its parameters, can compute the \accp{IM} or implement a specific \accp{IM}-informed \acc{DM} method.

The two formulations presented below are parameterized by a vector $p$, which is set according to the desired task.
These tasks are:
\begin{itemize}[label={\fbox{val.}-},align=parleft,labelindent=2ex]
	\item[\task{IM}] determining the \accp{IM} solutions,
	\label{task_IM}
	\item[\task{DM}] determining the solution to be implemented via \accp{IM}-informed \acc{DM}.
	\label{task_DM}
\end{itemize}
Executing task \hyperref[task_IM]{\task{IM}} enables the computation of the quantities discussed in \cref{list:objects_from_Phi}.
These, in turn, are required to compute the \accp{IM}-informed \acc{DM} solution in task \hyperref[task_DM]{\task{DM}}.

Applying the \acc{WS} scalarization to \eqref{eq:MOMPC} yields the single-objective problem
\begin{equation}
	\begin{aligned}
		\min_{\xMO \in \decSetFeas} \ &
		\param{w}_{\mathrm{WS}}\tr\JMO(\xMO),
	\end{aligned}
	\tag{$P_{\mathrm{WS}}(p,\decSetFeas)$}
	\label{eq:pOP_WS}
\end{equation}
which depends on the parameter vector $p$.
Furthermore, the \acc{PS} scalarization is augmented by introducing an additional parameter $\param{T}_{\J}$:
\begin{equation}
	\begin{aligned}
		\min_{\xMO \in \decSetFeas, \, l \in \Rset} \ & -l
		\\
		\text{s.t.} \quad
		&0 \leq \param{\J}_{\SO} + l \param{d} - \param{T}_{\J}\J(\xMO).		
	\end{aligned}
	\tag{$P_{\mathrm{PS}}(p,\decSetFeas)$}
	\label{eq:pOP_PS}
\end{equation}
\Cref{tab:p_WS_PS} details the composition of the parameter vector $\realize{p}$ for both problems. 
\newcommand{\pFontsize}{\normalsize}%
\begin{table}[tb]
	\pFontsize%
	\vspace*{-1\baselineskip}
	\begin{minipage}{1\linewidth}
		\hfil
		\begin{minipage}[t]{0.32\WidthBase}
			\begin{subequations}
				\begin{alignat}{1}
					p
					&=\
					\param{w}_{\mathrm{WS}}
					\label{eq:p_WS}
					\numberthis
					\\
					\overset{\text{\hyperref[task_IM]{\task{IM}}}}&{\equiv}
					e_i
					\label{eq:p_WS_IM}
					\numberthis
					\\
					\overset{\text{\hyperref[task_DM]{\task{DM}}}}&{\equiv}
					\realize{w}_{\mathrm{WS}}
					\label{eq:p_WS_FHP}
					\numberthis
				\end{alignat}
			\end{subequations}
		\end{minipage}%
		\hspace{1.1em}
		\begin{minipage}[t]{0.61\WidthBase}
			\begin{subequations}
				\begin{alignat}{7}
					p
					&=\
					& \big\{
						& \param{\J}_{\SO}, \
						& \param{d}, \
						& \param{T}_{\J}
					& \big\}
					\label{eq:p_PS}
					\numberthis
					\\
					\overset{\text{\hyperref[task_IM]{\task{IM}}}}&{\equiv}
					& \big\{
						& \boldsymbol{0}, \
						& -e_1, \
						& \matInline{e_i, \boldsymbol{0}}\tr
					& \big\}
					\label{eq:p_PS_IM}
					\numberthis
					\\
					\overset{\text{\hyperref[task_DM]{\task{DM}}}}&{\equiv}
					& \big\{
						& \realize{\J}_{\SO}, \
						& \realize{d}, \
						& I_{n_J}
					& \big\}
					\label{eq:p_PS_FHP}
					\numberthis
				\end{alignat}
			\end{subequations}
		\end{minipage}%
		\hfil
	\end{minipage}
	\\[-0\baselineskip]
	\caption[Task-dependent parameterization for the WS and PS formulations.]{
		Compositions of the parameter $p$ and its task-dependent realizations for the \acc{WS} formulation \eqref{eq:pOP_WS} (left) and the \acc{PS} formulation \eqref{eq:pOP_PS} (right).
		The index $i \in \{1,\ldots,n_J\}$ determines which $i$-th \acc{IM} is computed.
		The \emph{\realizeAccent} accent indicates that specific numerical values are assigned (cf. \cref{tab:notation_accents}). 
		For information on how these numerical values are computed for the various \acc{DM} methods and how they depend on $\betaPref$, see \cref{list:DM_methods}.
	}
	\label{tab:p_WS_PS}
\end{table}%
It also lists the specific parameter realizations for the tasks of computing the \accp{IM} and performing the \acc{DM}.
Note that with the parameter realization depicted in \eqref{eq:p_PS_IM}
one can effectively implement the \acc{WS} method using \acc{PS}:
\begin{equation}
	\begin{aligned}
		&\min_{x\in \mathbb{X}} \ -l
		\quad \text{s.t.} \quad
		\begin{mat}
			0 \\
			0_{n_J-1 \times n_J}
		\end{mat} 
		\leq 
		\begin{mat}
			0 - l - e_i\tr \J\\
			0_{n_J-1 \times n_J}
		\end{mat}
		\\[0.5\baselineskip]
		\equiv \quad
		&\max_{x\in \mathbb{X}} \ l
		\quad \text{s.t.} \quad
		e_i\tr \J \leq -l.
	\end{aligned}
\end{equation}
To facilitate readability and serve as a quick reference, \cref{tab:important_quantities} summarizes the most important quantities and mathematical notation used throughout this paper.
\begin{table}[htb]
    \centering
    \tableFontSizeWiley 
	\begin{tblr}{colspec={r|ll}, rowsep = 0.15em}
	    \toprule
	    Variable & Description & Reference \\
	    \midrule
	    $\param{w}_{\WS}$ & weights of \acc{WS} method & \eqref{eq:pOP_WS} \\
	    $\param{\J}_{\SO}$ & \acc[l]{SO} (\acc{PS} method)& \eqref{eq:pOP_PS} \\
	    $\param{d}$ & \acc[l]{SDV} (\acc{PS} method)& \eqref{eq:pOP_PS} \\
	    $\param{T}_{\J}$ & {transformation matrix for \\ objective vector (\acc{PS} method)} & \eqref{eq:pOP_PS} \\
	    $\J$ & objective vector & \eqref{eq:MOOP} \\
	    $\betaPref$ & preference vector & \cref{ssec:Specifying_Relative_Importance} \\
	    $\Phi$ & pay-off matrix &  \eqref{eq:Phi} \\
	    $\J_{\NP}$, $\J_{\UP}$ & \acc[l]{NP} and \acc[l]{UP}& \cref{list:objects_from_Phi} \\
	    $\cNPUP$ & image space normalization matrix & \cref{list:objects_from_Phi} \\
	    $\wChim, \wChimQuasi, \wChimVisual$ & {normals related to \acc{CHIM}:\\ standard, quasi- and visual normal} 
	        & \cref{list:objects_from_Phi} \\
 	    $\wChimQuasi_{\mathrm{c}}$ & {quasi-normal originating from\\ barycentric center of the \acc{CHIM}} 
 	        & \cref{list:objects_from_Phi} \\
	    $\tilde{e}_i$ & regularized unit vector & \eqref{eq:unit_vec_reg}\\
	    $\tilde{\Phi}$ & shifted pay-off matrix & \eqref{eq:Phi_shifted} \\
	    \bottomrule
	\end{tblr}
    \caption{Summary of important quantities.}%
    \label{tab:important_quantities}%
\end{table}

\subsection{Development of Performance Metrics}
\label{ssec:Performance_Metrics}
\paragraph*{Coverage}
First, the \acc{DM} method should be able to cover large parts of the \acc{PF}.
It therefore becomes reasonable to consider the ratio between the surface area that a \acc{DM} method can generate by sampling the \acc{PF} and the total surface area of the \acc{PF} as a coverage metric $m_{\mathrm{cov}}$.
Note that these surface area computations are also done in the \acc{NIS}.

The basis of the surface area computations are triangle meshes that are constructed by the ball pivoting algorithm \cite{bernardini_ball-pivoting_1999,digne_analysis_2014} using (dense) samples of the \acc{PF}.
The mesh reconstructed from a dense sampling $\imSetPareto[\tilde]$ of the \acc{PF} is denoted by $\imSetPareto[\tilde][,\mathrm{mesh}]$; details are given in \cref{sec:Approx_Karcher_Mean}.
The vertices of the mesh triangles are generated by the different \acc{DM} methods using $\baryCO$-values from a fine, uniform sampling of the unit simplex. 
This method of computing areas can be used for both convex and non-convex surfaces.

\paragraph*{Translation}
Second, the \acc{DM} method should yield a good translation of preferences to points on the \acc{PF}.
Such an assessment requires the establishment of a robust and geometrically meaningful mapping $\J_{\mathrm{des}}(\betaPref)$ from a preference vector $\betaPref \in \FunitSimplex$ to a corresponding point on the \acc{PF} $\imSetPareto$.
\\
To achieve this, the concept of Riemannian barycentric coordinates \cite{deylen_distortion_2016} is employed, which is defined via the Karcher mean \cite{karcher_riemannian_1977}.
The basis of the Karcher mean is the geodesic distance from a point $X_1$ to a point $X_2$ on a Riemannian manifold $S$
$d_{\mathrm{geod}}(S,X_1,X_2)$.
Details on the definition of the Karcher mean, its approximation using finitely many samples from $S$ and the associated implementation,
and the determination of the mapping $\J_{\mathrm{des}}(\betaPref)$
can be found in \cref{sec:Approx_Karcher_Mean}.

A meaningful metric for evaluating the translation quality can subsequently be formulated.
This metric is based on the geodesic distance between the expected point on the \acc{PF} $\J_{\mathrm{des}}(\betaPref)$ and the actual point realized by the $k$-th \acc{DM} method $\J_{\mathrm{DM},k}(\betaPref)$, yielding
\begin{equation}
	e_{\mathrm{transl},\, k}(\betaPref) := d_{\mathrm{geod}}
	\(
		\imSetPareto[\tilde][,\mathrm{mesh}], 	
		\J_{\mathrm{des}}(\betaPref),
		\J_{\mathrm{DM},k}(\betaPref)
	\).
\end{equation}
The translation quality metric is defined as
\begin{equation}
	m_{\mathrm{transl},\, k} :=
	1 -
	\frac{1}{s_{\mathrm{transl}}}
	\frac{1}{N_{\betaPref}}
	\sum_{i=1}^{N_{\betaPref}} 
	e_{\mathrm{transl},\, k}( \betaPref[i])
	\label{eq:translation_metric}
\end{equation}
where the scaling factor
\begin{equation}
	s_{\mathrm{transl}} =
	\frac{1}{n_J}
	\sum_{i=1}^{n_J} d_{\mathrm{geod}}
	\(
		\imSetPareto[\tilde][,\mathrm{mesh}], 	
		\J_{\mathrm{des}}(\beta_{\mathrm{c}}),
		\JMO(\xMO_{i}\opt)
	\)
	\label{eq:factor_translation_metric}
\end{equation}
is introduced to scale the geodesic distances suitably.
It is the mean geodesic distance between $\J_{\mathrm{des}}(\beta_{\mathrm{c}})$, the point assigned to the central preference $\beta_{\mathrm{c}} = 1/n_J \cdot \boldsymbol{1}$, and the $n_J$ \accp{IM} $\JMO(\xMO_{i}\opt)$.
Note that the value of \eqref{eq:translation_metric} can theoretically become negative since $s_{\mathrm{transl}}$ does not represent the maximum geodesic distance on $\imSetPareto[\tilde][,\mathrm{mesh}]$.

\subsection{Application to Exemplary Problems}
\label{ssec:Exemplary_Problems}
All six methods from \cref{list:DM_methods} are applied to two exemplary \accp{MOOP}.
The performance of the methods are evaluated and discussed based on two metrics from \cref{ssec:Performance_Metrics}.
\\
The two exemplary \accp{MOOP} are the following:

\begin{bThm}{eexample} 
\label{example1}
The first \acc{MOOP} with $n_J=3$ objectives is
\begin{equation}
	\begin{aligned}
		\min_{\xMO \in \Rset^3} \ &\J(\xMO) = \diag(\matInline{l_1,l_2,l_3}) \, \xMO
		\\
		\text{s.t.} \quad
		& \xMO_1^2 + \xMO_2^2 + \xMO_3^2 \leq 1.
	\end{aligned}
	\tag{$\mathsf{MOOP}_{\mathrm{ex,1}}$}
	\label{eq:MOOP_ex_1}
\end{equation}
The feasible image set $\imSetFeas$ is a non-rotated ellipsoid with the semi-axis lengths $l_1,l_2,l_3>0$ centered at the origin.
For $\matInline{l_1,l_2,l_3}=\matInline{1,1,1}$, its \acc{PF} is shown in \cref{fig:mainfig}.
In the following, only $\matInline{l_1,l_2,l_3}=\matInline{1,10,100}$ is considered.
The \acc{PF} of \tagref{eq:MOOP_ex_1} can be classified as an \emph{ideal example} of a \acc{PF} in terms of convexity and the location of the \accp{IM}, cf.~\cref{fig:DM_all_methods} (left).
Specifically, in \acc{NIS}, the normal vector of the \acc{CHIM} is $-1/n_J \cdot \boldsymbol{1}$ and the shifted pay-off matrix is
$\cNPUP(\Phi - \J_{\UP}) = 1_{n_J \times n_J} - \diag(\boldsymbol{1})$.
\end{bThm}

\begin{bThm}{eexample}
\label{example2}
The second \acc{MOOP} builds upon \tagref{eq:MOOP_ex_1} with $\matInline{l_1,l_2,l_3}=\matInline{1,10,100}$ by adding three inequality constraints to obtain
\begin{equation}
	\begin{aligned}
		\min_{\xMO \in \Rset^3} \ &\J(\xMO) = \diag(\matInline{l_1,l_2,l_3}) \, \xMO
		\\
		\text{s.t.} \quad
		& \xMO_1^2 + \xMO_2^2 + \xMO_3^2 \leq 1,
		\\
		& w_i\tr \xMO + b_i \geq 0,
		\qquad
		i=1,2,3
	\end{aligned}
	\tag{$\mathsf{MOOP}_{\mathrm{ex,2}}$}
	\label{eq:MOOP_ex_2}
\end{equation}
with
\begin{align}
	\matInline{w_1,w_2,w_3} &=
	\begin{mat}
		0.41003    & -0.63397 & \hphantom{-}0.1499 \\
		0.58997    & \hphantom{-}0.36603 &-0.8501 \\
		-1.254 \cdot 10^{-6} &        0 &      0
	\end{mat}
	,
	\\
	\matInline{b_1,b_2,b_3} &= \matInline{0.50625, 0, 0}
	.	
\end{align}
Each $w_i$ and $b_i$, respectively, contains the normal vector and the bias of a plane.
\\
The \acc{PF} of \tagref{eq:MOOP_ex_2} is depicted in \cref{fig:DM_all_methods} (right) and, together with the quantities derived from $\Phi$, in \cref{fig:quants_from_Phi}.
The location of its \accp{IM} in a \emph{\acc{NIS}}%
\footnote{Here, no specific notation is introduced to indicate the use of normalized objectives.}
results in the quantities
\begin{equation}
    \begin{split}
		\wChim &= \matInline{-0.487, 0.035, -0.477}\tr,
		\quad \text{\small(\acc{CHIM} normal)}
		\\
		\tilde{\Phi}
		&= 
		\begin{mat}
			0.000 &0.641 &1.000 \\
			0.716 &0.000 &1.000 \\
			1.000 &0.293 &0.000 \\
		\end{mat},
		\quad 
		\substack{
			\text{\small (shifted  pay-} 
			\\ 
			\text{\small off matrix)}
		}
		\\
		\theta &= 38.217^{\circ},
    \end{split}
\end{equation}
with $\theta = \min\{\measuredangle(\wChim,-1/n_J\boldsymbol{1}), \, 180^\circ - \measuredangle(\wChim,-1/n_J\boldsymbol{1})\}$,
where taking the minimum makes $\theta$ independent of the orientation of the normal vector.
The difference between the shown values and the \enquote{ideal} ones suggests that this \acc{PF} misses important properties of an ideal \acc{PF}.
The angle deviation is another clear indicator of this.
\end{bThm}

\paragraph*{Discussion}
For intuitive visual evaluation, each preference vector is assigned its own color.
This color mapping is shown in \cref{fig:two_unit_simplex}.
\begin{figure}[htb]
	\centering
		\adjincludegraphics[
			width=0.75\WidthBase, clip, 
			trim={{0.08\width} {0.02\height} {0\width} {0\height}}
		]{two_unit_simplex.\FigureFileExt}
	\caption{Assignment of each point on the $2$-unit simplex to a unique color.}
	\label{fig:two_unit_simplex}
\end{figure}
The results from \cref{example1} and \cref{example2} are presented in \cref{fig:DM_all_methods}
and quantitatively evaluated in \cref{tab:metrics_DM}.
\newlength{\DMFigWidth}%
\setlength{\DMFigWidth}{0.427\textwidth}%
\sbox{\tempbox}{%
    \includegraphics[width=\DMFigWidth]{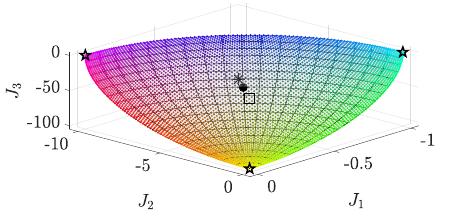}%
}%
\begin{figure*}[htbp]
    \centering
    
    \begin{subfigure}[t]{\DMFigWidth}
        \centering
        \includegraphics[width=\textwidth]{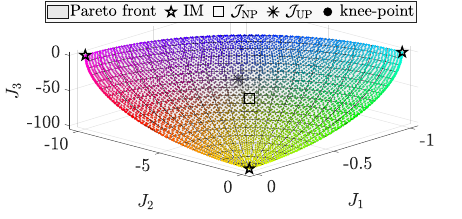}
    \end{subfigure}
    \InfoBox{
   		$\J_{\mathrm{des}}(\betaPref)$ obtained via the Karcher mean
    }
    \begin{subfigure}[t]{\DMFigWidth}
        \centering
        \includegraphics[width=\textwidth]{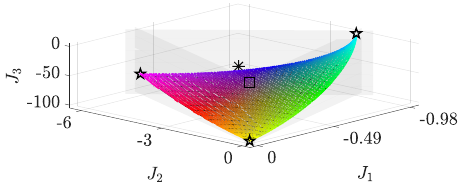}
    \end{subfigure}
    
    \begin{subfigure}[t]{\DMFigWidth}
        \centering
        \includegraphics[width=\textwidth]{WS_1__ex_1.\FigureFileExt}
    \end{subfigure}
    \InfoBox{
    	\ref{enum:WS_1}:\\
		\acc[s]{WS} \\
		with \acc{IM}\\
		scaling
		\\
		\&
		\\
       	\ref{enum:WS_2}:\\
		knee-point
    }
    \begin{subfigure}[t]{\DMFigWidth}
        \centering
        \includegraphics[width=\textwidth]{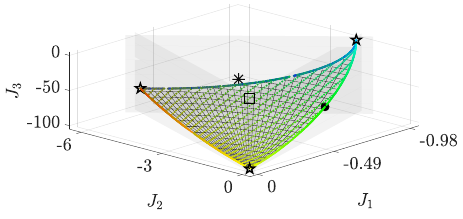}
    \end{subfigure}

    \begin{subfigure}[b]{\DMFigWidth}
        \centering
        \includegraphics[width=\textwidth]{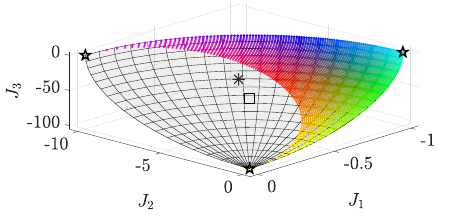}
    \end{subfigure}
    \InfoBox{
    	\ref{enum:PS_1}:\\
		\acc[s]{NBI} \\
		with \\
		$\realize{d} = \wChim$
    }
    \begin{subfigure}[b]{\DMFigWidth}
        \centering
        \includegraphics[width=\textwidth]{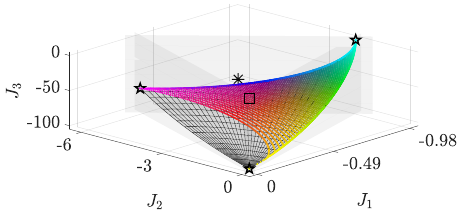}
    \end{subfigure}

    \begin{subfigure}[t]{\DMFigWidth}
        \centering
        \includegraphics[width=\textwidth]{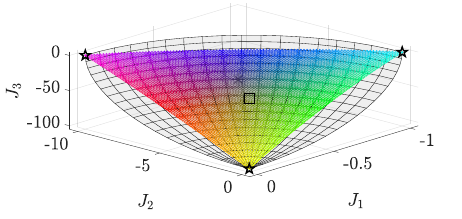}
    \end{subfigure}
    \InfoBox{
       	\ref{enum:PS_2}:\\
		\acc[s]{NBI} \\
		with \\
		$\realize{d} = \wChimQuasi_{\mathrm{c}}$
    }
    \begin{subfigure}[t]{\DMFigWidth}
        \centering
        \includegraphics[width=\textwidth]{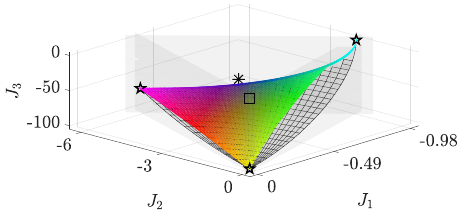}
    \end{subfigure}
    
    \begin{subfigure}[t]{\DMFigWidth}
        \centering
        \includegraphics[width=\textwidth]{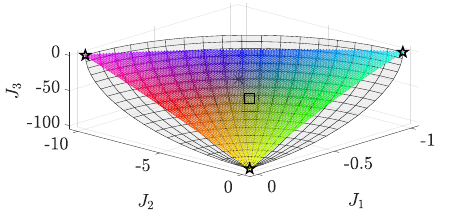}
    \end{subfigure}
    \InfoBox{
       	\textbf{\ref{enum:PS_3}}:\\
   		\acc{NBI}\\
   		with\\
   		$\realize{d}=\wChimVisual$
    }
    \begin{subfigure}[t]{\DMFigWidth}
        \centering
        \includegraphics[width=\textwidth]{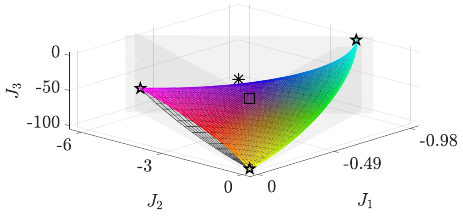}
    \end{subfigure}
    
    \begin{subfigure}[t]{\DMFigWidth}
        \centering
        \includegraphics[width=\textwidth]{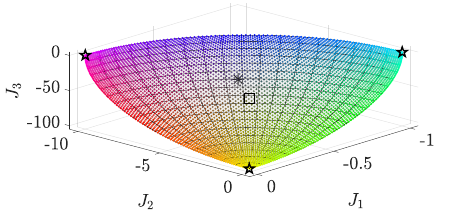}
    \end{subfigure}
    \InfoBox{
       	\textbf{\ref{enum:PS_4}}:\\
		\acc{NCHIM}
    }
    \begin{subfigure}[t]{\DMFigWidth}
        \centering
        \includegraphics[width=\textwidth]{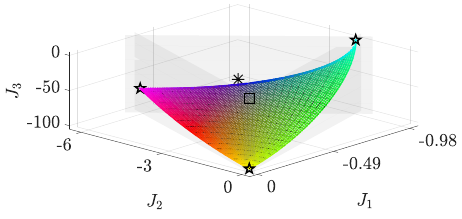}
    \end{subfigure}
    
    \caption[
    	{\protect\Accp[l]{PF}} of \cref{example1} (left) and of \cref{example2} (right).
    ]{
    	\Accp[l]{PF} of \cref{example1} (left) and of \cref{example2} (right).
    	The first row shows $\J_{\mathrm{des}}(\betaPref)$, whereas the remaining rows show $\J_{\DM}(\betaPref)$ for the various \acc{DM} methods.
    	Note that the points $\J_{\mathrm{des}}(\betaPref)$, by definition, coincide with the vertices of the underlying mesh.
    	Hence, their location is only an approximation for mesh triangles of non-vanishing size, which explains the imperfect positioning of the points in the upper figures.
    }
    \label{fig:DM_all_methods}
\end{figure*}

\definecolor{cA1}{rgb}{0.0000,0.7490,0.1294}
\definecolor{cA2}{rgb}{0.0000,0.7490,0.1294}
\definecolor{cA3}{rgb}{0.2125,0.7211,0.1305}
\definecolor{cA4}{rgb}{0.4967,0.5662,0.1343}
\definecolor{cB1}{rgb}{0.6115,0.4144,0.1350}
\definecolor{cB2}{rgb}{0.6856,0.2077,0.1338}
\definecolor{cB3}{rgb}{0.6899,0.1860,0.1337}
\definecolor{cB4}{rgb}{0.6522,0.3243,0.1346}
\definecolor{cC1}{rgb}{0.5913,0.4492,0.1351}
\definecolor{cC2}{rgb}{0.4967,0.5662,0.1343}
\definecolor{cC3}{rgb}{0.5355,0.5251,0.1348}
\definecolor{cC4}{rgb}{0.4098,0.6354,0.1331}
\definecolor{cD1}{rgb}{0.5070,0.5561,0.1345}
\definecolor{cD2}{rgb}{0.4747,0.5863,0.1341}
\definecolor{cD3}{rgb}{0.3420,0.6739,0.1321}
\definecolor{cD4}{rgb}{0.3420,0.6739,0.1321}
\definecolor{cE1}{rgb}{0.5070,0.5561,0.1345}
\definecolor{cE2}{rgb}{0.2742,0.7023,0.1312}
\definecolor{cE3}{rgb}{0.3420,0.6739,0.1321}
\definecolor{cE4}{rgb}{0.2125,0.7211,0.1305}
\definecolor{cF1}{rgb}{0.1152,0.7397,0.1298}
\definecolor{cF2}{rgb}{0.1152,0.7397,0.1298}
\definecolor{cF3}{rgb}{0.2125,0.7211,0.1305}
\definecolor{cF4}{rgb}{0.2742,0.7023,0.1312}

\newcommand{\CeCo}[1]{\cellcolor{#1}}

\DeclareRobustCommand{\rowA}{& \CeCo{cA1} 100.00 & \CeCo{cA2} 100.00 & \CeCo{cA3} 94.03 & \CeCo{cA4} 63.76}
\DeclareRobustCommand{\rowB}{& \CeCo{cB1} 35.10 & \CeCo{cB2} 6.49 & \CeCo{cB3} 5.91 & \CeCo{cB4} 21.17}
\DeclareRobustCommand{\rowC}{& \CeCo{cC1} 41.80 & \CeCo{cC2} 63.98 & \CeCo{cC3} 54.01 & \CeCo{cC4} 76.58}
\DeclareRobustCommand{\rowD}{& \CeCo{cD1} 61.08 & \CeCo{cD2} 66.12 & \CeCo{cD3} 85.86 & \CeCo{cD4} 85.37}
\DeclareRobustCommand{\rowE}{& \CeCo{cE1} 61.08 & \CeCo{cE2} 91.01 & \CeCo{cE3} 85.85 & \CeCo{cE4} 94.45}
\DeclareRobustCommand{\rowF}{& \CeCo{cF1} 100.00 & \CeCo{cF2} 99.99 & \CeCo{cF3} 94.02 & \CeCo{cF4} 91.48}

\definecolor{cmap0}{rgb}{0.6980,0.1333,0.1333}
\definecolor{cmap1}{rgb}{0.6940,0.1616,0.1335}
\definecolor{cmap2}{rgb}{0.6899,0.1860,0.1337}
\definecolor{cmap3}{rgb}{0.6856,0.2077,0.1338}
\definecolor{cmap4}{rgb}{0.6812,0.2275,0.1339}
\definecolor{cmap5}{rgb}{0.6767,0.2459,0.1341}
\definecolor{cmap6}{rgb}{0.6721,0.2631,0.1342}
\definecolor{cmap7}{rgb}{0.6673,0.2794,0.1343}
\definecolor{cmap8}{rgb}{0.6624,0.2950,0.1344}
\definecolor{cmap9}{rgb}{0.6574,0.3099,0.1345}
\definecolor{cmap10}{rgb}{0.6522,0.3243,0.1346}
\definecolor{cmap11}{rgb}{0.6469,0.3382,0.1347}
\definecolor{cmap12}{rgb}{0.6414,0.3517,0.1348}
\definecolor{cmap13}{rgb}{0.6358,0.3648,0.1349}
\definecolor{cmap14}{rgb}{0.6300,0.3776,0.1349}
\definecolor{cmap15}{rgb}{0.6240,0.3901,0.1350}
\definecolor{cmap16}{rgb}{0.6179,0.4024,0.1350}
\definecolor{cmap17}{rgb}{0.6115,0.4144,0.1350}
\definecolor{cmap18}{rgb}{0.6050,0.4262,0.1350}
\definecolor{cmap19}{rgb}{0.5982,0.4378,0.1351}
\definecolor{cmap20}{rgb}{0.5913,0.4492,0.1351}
\definecolor{cmap21}{rgb}{0.5841,0.4604,0.1350}
\definecolor{cmap22}{rgb}{0.5767,0.4715,0.1350}
\definecolor{cmap23}{rgb}{0.5690,0.4825,0.1350}
\definecolor{cmap24}{rgb}{0.5611,0.4933,0.1350}
\definecolor{cmap25}{rgb}{0.5528,0.5040,0.1349}
\definecolor{cmap26}{rgb}{0.5443,0.5146,0.1348}
\definecolor{cmap27}{rgb}{0.5355,0.5251,0.1348}
\definecolor{cmap28}{rgb}{0.5264,0.5355,0.1347}
\definecolor{cmap29}{rgb}{0.5169,0.5458,0.1346}
\definecolor{cmap30}{rgb}{0.5070,0.5561,0.1345}
\definecolor{cmap31}{rgb}{0.4967,0.5662,0.1343}
\definecolor{cmap32}{rgb}{0.4859,0.5763,0.1342}
\definecolor{cmap33}{rgb}{0.4747,0.5863,0.1341}
\definecolor{cmap34}{rgb}{0.4630,0.5962,0.1339}
\definecolor{cmap35}{rgb}{0.4507,0.6061,0.1337}
\definecolor{cmap36}{rgb}{0.4378,0.6159,0.1335}
\definecolor{cmap37}{rgb}{0.4242,0.6257,0.1333}
\definecolor{cmap38}{rgb}{0.4098,0.6354,0.1331}
\definecolor{cmap39}{rgb}{0.3945,0.6451,0.1329}
\definecolor{cmap40}{rgb}{0.3782,0.6547,0.1327}
\definecolor{cmap41}{rgb}{0.3608,0.6643,0.1324}
\definecolor{cmap42}{rgb}{0.3420,0.6739,0.1321}
\definecolor{cmap43}{rgb}{0.3216,0.6834,0.1318}
\definecolor{cmap44}{rgb}{0.2991,0.6928,0.1315}
\definecolor{cmap45}{rgb}{0.2742,0.7023,0.1312}
\definecolor{cmap46}{rgb}{0.2458,0.7117,0.1309}
\definecolor{cmap47}{rgb}{0.2125,0.7211,0.1305}
\definecolor{cmap48}{rgb}{0.1716,0.7304,0.1302}
\definecolor{cmap49}{rgb}{0.1152,0.7397,0.1298}
\definecolor{cmap50}{rgb}{0.0000,0.7490,0.1294}

\newcommand{\plotColorAxis}[3]{%
  \begin{tikzpicture}
    \foreach \y [count=\i from 0] in {0,1,...,50} {
      \pgfmathsetmacro{\ybottom}{( \y   / 51.0) * #2}
      \pgfmathsetmacro{\ytop}{  ((\y+1) / 51.0) * #2}
      \fill[cmap\i, draw=none] (0, \ybottom #3) rectangle (#1 #3, \ytop #3);
    }
    \pgfmathsetmacro{\ytick}{0.0000 * #2}
    \draw (#1 #3, \ytick #3) -- ++(0.1cm, 0) node[right, font=\small] {$\text{0.00}$};
    \pgfmathsetmacro{\ytick}{0.2500 * #2}
    \draw (#1 #3, \ytick #3) -- ++(0.1cm, 0) node[right, font=\small] {$\text{25.00}$};
    \pgfmathsetmacro{\ytick}{0.5000 * #2}
    \draw (#1 #3, \ytick #3) -- ++(0.1cm, 0) node[right, font=\small] {$\text{50.00}$};
    \pgfmathsetmacro{\ytick}{0.7500 * #2}
    \draw (#1 #3, \ytick #3) -- ++(0.1cm, 0) node[right, font=\small] {$\text{75.00}$};
    \pgfmathsetmacro{\ytick}{1.0000 * #2}
    \draw (#1 #3, \ytick #3) -- ++(0.1cm, 0) node[right, font=\small] {$\text{100.00}$};
  \end{tikzpicture}
}

\begin{table}[h]
	\centering
	\newcommand{\TableFontSize}{\small}
	\newsavebox{\myTblBox}
	\sbox{\myTblBox}{
		\renewcommand{\arraystretch}{1.3}
		\setlength{\tabcolsep}{3pt}
		\begin{tabular}[b]{ | *{5}{>{\TableFontSize}c |} }
			\hline
			
			&\multicolumn{2}{>{\TableFontSize}c|}{$m_{\mathrm{cov}}$ in \%}
			&\multicolumn{2}{>{\TableFontSize}c|}
				{$m_{\mathrm{transl}}$ in \%}
			\\ \hline
			method 
			& Ex. 1 & Ex. 2 
			& Ex. 1 & Ex. 2 
			\\ 
			\hline
			\ref{enum:WS_1} \rowA \\
			\ref{enum:PS_1} \rowC \\
			\ref{enum:PS_2} \rowD \\
			\textbf{\ref{enum:PS_3}} \rowE \\
			\textbf{\ref{enum:PS_4}} \rowF \\
			\hline
		\end{tabular}
	}
	\pgfmathsetmacro{\myTableHeightCM}{
		0.035277777777782*\dimexpr\the\ht\myTblBox\relax 
	}
	
		\begin{adjustbox}{valign=m}
			\usebox{\myTblBox}
		\end{adjustbox}%
		\hspace{0.5em}%
		\begin{adjustbox}{valign=m}
			\plotColorAxis{0.5}{\myTableHeightCM}{cm}
		\end{adjustbox}
	
	\vspace{-0.3\baselineskip}
	\caption{
		Quantitative evaluation of the \acc{DM} methods' performance when applied to \cref{example1} and \cref{example2}, as shown in \cref{fig:DM_all_methods}. 
		The novel methods are highlighted in bold.
		The cell colors visually support the values on an absolute scale from $0$ to $100\%$.
	}
	\label{tab:metrics_DM}
\end{table}

\begin{bThm}[Sensitivity to mesh resolution]{remark}
	\label{rem:sens_mesh_resolution}
	The following statements apply to both examples and all \acc{DM} methods:
	\\
	The mesh of the \acc{DM} method
	is computed for two different mesh resolutions.
	Increasing the number of mesh triangles by more than $37\%$ results in a decrease of the average triangle area by more than $24\%$.
	However, the respective \acc{PF} surface area values differ by at most 
	$0.013\%$.
	Hence, the finer mesh is considered a sufficiently accurate approximation and is therefore used for the surface area computation.
	\\
	The same computations are repeated for two further refined meshes. 
	The finest of these serves as the \emph{reference} for the \acc{PF} surface area and the mapping $\J_{\mathrm{des}}(\betaPref)$, which in turn determines the translation metric \eqref{eq:translation_metric}.
	The computations yield the percentage values $23\%$, $18\%$, and
	$0.0057\%$ (for the number of triangles, average triangle area, and total \acc{PF} surface area, respectively); 
	therefore, the same conclusion can be drawn as for the \acc{DM} mesh.	
	Furthermore, an analysis of the translation metric across the two different \emph{reference} mesh resolutions reveals a maximum difference of $0.17\%$.
	
	In summary, these results indicate that the \acc{PF} surface areas, and therefore also the coverage metrics, are sufficiently accurate.
	Moreover, the low sensitivity of the translation metric further validates the chosen \emph{reference} mesh resolution.
\end{bThm}
Since the \acc{KP} case \ref{enum:WS_2} yields only a single point, it has to be treated separately.
In \cref{example1} it yields a balanced trade-off between all objectives.
However, this is not the case in \cref{example2} where one rather obtains a trade-off between the first and third \acc{IM} on the boundary of the \acc{PF}.
This means that the preference was translated non-intuitively.

In terms of \emph{coverage}, 
one can see in \cref{tab:metrics_DM} that in \cref{example1} 
the methods \ref{enum:WS_1} and \ref{enum:PS_4} perform the best,
\ref{enum:PS_2} and \ref{enum:PS_3} perform acceptably, and
\ref{enum:PS_1} performs poorly.
This changes with \cref{example2}, where \ref{enum:PS_1}, \ref{enum:PS_3} perform better and the remaining methods perform similarly.
\\
In terms of \emph{translation}, in \cref{example1},
according to the presented metric \eqref{eq:translation_metric}, the methods \ref{enum:WS_1} and \ref{enum:PS_4} perform exceptionally, while \ref{enum:PS_2}, \ref{enum:PS_3} yield good results and the remaining strategies perform poorly.
When comparing the values from \cref{example2} with \cref{example1}, 
the results for \ref{enum:PS_1}, \ref{enum:PS_3} improve, 
whereas \ref{enum:WS_1} deteriorates
and \ref{enum:PS_2}, \ref{enum:PS_4} yield a similar performance.
Note that, in \cref{example1}, the metrics of \ref{enum:WS_1} and \ref{enum:PS_4} as well as those of \ref{enum:PS_2} and \ref{enum:PS_3} almost coincide.
For \ref{enum:PS_2} and \ref{enum:PS_3}, this is expected, since the quasi-normal and the visual normal coincide for the shifted pay-off matrix of this ideal \acc{PF}.
Hence, \cref{example1} does not discriminate between these pairs of methods.

A visual inspection of \cref{fig:DM_all_methods} reveals that the method \ref{enum:WS_1} applied to \cref{example2} yields a poor translation, compared to \ref{enum:PS_3} or \ref{enum:PS_4}.
Specifically, only some shades of green are in the interior of the \acc{PF} $\imSetPareto$, while most preferences yield points on the boundary of $\imSetPareto$.
However, the translation metric yields $63.76\%$, which suggests that the value of \eqref{eq:factor_translation_metric} is relatively large.
This insight makes it clear that values above $90\%$ signify an excellent translation, which is further confirmed by the visual assessment of the corresponding methods in \cref{fig:DM_all_methods}.

In summary, in the two examples, the established method \ref{enum:WS_1} performs well for the ideal \acc{PF} of \cref{example1} but deteriorates markedly for \cref{example2}, whereas \ref{enum:PS_1} already performs poorly for the ideal \acc{PF}.
\hyperanchor{anchor:NBI_robust}
If, instead, one uses $\wChimQuasi_{\mathrm{c}}$ or $\wChimVisual$ for \acc{NBI}, the results become much better and remain more consistent across \cref{example1,example2},
which indicates lower sensitivity to differing objective ranges;
see also the translation metrics in \cref{tab:metrics_DM}.
Altogether, the evaluation of the two examples suggests that the novel \acc{DM} methods \ref{enum:PS_3} and \ref{enum:PS_4} appear most promising.

\section{Individual Minima-Informed MOMPC}
\label{sec:IMI_MOMPC}
Now, \acc{MOMPC} problems are considered, which are transcribed to a \acc{MOOP} such that the scalarizations \eqref{eq:pOP_WS} and \eqref{eq:pOP_PS} can be applied.
It is proven that these methods stabilize a (controlled) fixed point.
The corresponding stability argument is independent of the specific geometry of the \accp{IM};
the latter only affects the selected point on the \acc{PF} through the \acc{DM} process.
Moreover, an easy-to-implement recipe for constructing stabilizing terminal ingredients is provided.

\subsection{Formulating the MOMPC problem}
\label{ssec:Formulating_MOMPC_problem}
An \acc{MOMPC} problem can now be formulated, together with a closed-loop algorithm based on the scalarization-based \acc{DM} methods derived above.
For the dynamical system, the \acc{DT} dynamics $\xSys_{j+1} = f_{\mathrm{dyn}}(\xSys_j,\uSys_j)$ are considered, with the system states $\xSys_j \in \Rset^{n_\xSys}$ and system inputs $\uSys_j \in \Rset^{n_\uSys}$ at time instances $j=0,1,\dots,N-1$ for $N\in \mathbb{N}$.
In total, the \emph{dynamical} constraints considered in the \acc{MOMPC} are
\begin{equation}
	\begin{alignedat}{2}
		\xSys_{j+1} &= f_{\mathrm{dyn}}(\xSys_j,\uSys_j),
		\quad &&j = 0,1,\ldots,N-1,	
		\\
		\xSys_j &\in \xSet, \ \uSys_j \in \uSet,
		&&j = 0,1,\ldots,N-1,
		\\
		\xSys_0 &= \param{\xSys}^0, \ \xSys_N \in \xSet\term
	\end{alignedat}
	\tag{$\mathcal{C}(\xMO,\param{\xSys}^0)$}
	\label{eq:dyn_constraints}
\end{equation}
with initial state $\param{\xSys}^0 \in \Rset^{n_\xSys}$ and terminal set $\xSet\term \subseteq \Rset^{n_\xSys}$.
The vector
$\xMO=[\xSys_0 \concSep \ldots \concSep \xSys_N \concSep \uSys_0 \concSep \ldots \concSep \uSys_{N-1}]$ 
shall become the optimization variable.
The set of optimization variables feasible w.r.t.\ the 
constraints \eqref{eq:dyn_constraints} is denoted by
\begin{equation}
	\decSetFeasMPC(\param{\xSys}^0):=\left\{\xMO \in \Rset^{n_{\xMO}} \ \middle| \ \text{\eqref{eq:dyn_constraints} can be satisfied}\right\},
\end{equation}
which is parameterized by the initial state $\param{\xSys}^0$.
Then, in \acc{MOMPC}, the objective vector
\begin{subequations}
	\label{eq:J_MOMPC}%
	\allowdisplaybreaks
	\begin{align}
		\JMO(\xMO) &= \left[J^N_1(\xMO) \concSep[\middle] \ldots \concSep[\middle] J^N_{n_J}(\xMO)\right],
		\numberthis
		\\
		\text{with}\quad
		J^N_i (\xMO) &=\sum_{j=0}^{N-1} \ell_i(\xSys_j,v_j) + F_i(\xSys_N) \, \forall i=1,\dots, n_J,
		\numberthis
	\end{align}
\end{subequations}
has to be minimized under the constraint that $\xMO \in \decSetFeasMPC(\param{\xSys}^0)$, on a moving horizon of length $N$ with changing current state $\param{\xSys}^0$.

The \acc{OP}
$
	\min_{\xMO \in \decSetFeasMPC(\param{\xSys}^0)}
	\JMO(\xMO)
$
has two terminal ingredients: the terminal costs $F_i$, $i\in\{1,\ldots,n_J\}$ in \eqref{eq:J_MOMPC} and the terminal set $\xSet\term$ in \eqref{eq:dyn_constraints}.
For appropriately defined terminal ingredients, the \acc{OP} belongs to the class of \acc{MO} quasi-infinite horizon \acc{MPC} (cf.~\cite{chen_quasi-infinite_1998}) problems.

In order to prove stability of the proposed \acc{MOMPC} algorithm, the admissible set is further constrained by an upper bound on the objective values.
The parameterized feasible set is defined by
\begin{equation}
	\decSetFeasMOMPC\big(\xSys^0,\J_{\ub}\big)
	:=
	\left\{
		\xMO \in \decSetFeasMPC\big(\xSys^0\big) 
	\ \middle| \ 
		\J(\xMO) \leq \J_{\ub}
	\right\}.
\end{equation}
Then, 
with an initial point $\param{\xSys}^0$ and the upper bound on the cost $\param{\J}_{\ub}$,
the \acc{DT} \acc{MOMPC} \acc{OP} with terminal ingredients is given by
\begin{equation}
	\begin{aligned}
			\min_{\xMO} \ &\JMO(\xMO)
			\\
			\text{s.t.} \quad
			&\xMO \in \decSetFeasMOMPC\big(\param{\xSys}^0,\param{\J}_{\ub}\big).
			\label{eq:MOMPC}
	\end{aligned}
	\tag{$\mathcal{P}(\param{\xSys}^0,\param{\J}_{\ub})$}
\end{equation}
Note that \ref{eq:MOMPC} possesses the structure of \eqref{eq:MOOP} as introduced in~\cref{ssec:Intro_MOO}.
Thus, it can be addressed using one of the scalarization methods derived in~\cref{ssec:IMG_Prob_Formulation}.
Two corresponding variants are proposed in \cref{alg:MOMPC}.

\newcommand{\WScount}{i)}%
\newcommand{\PScount}{ii)}%
\begin{algorithm}[htbp]
	\caption{%
		\acc{IM}-informed MOMPC for the scalarization approaches
		\WScount~WS and \PScount~PS
	}
	\label{alg:MOMPC}
	\vspace{0.2\baselineskip}
	\textit{Parameters and options:} 
		$\J_{\mathrm{conv}}$;
		$\boldsymbol{\delta}$ (cf.~\eqref{eq:bs_delta_restriction});
		initial descent condition approach: a) or b).
	\vspace{-0.5\baselineskip}
	\begin{steps}[label={\arabic*)},itemsep=0.5\baselineskip,labelindent=0.5ex]
		\setcounter{stepsi}{-1}
		\item 
		\label{step:MOMPC_init}
		Set $k=0$ and $t_0=t^0$.
		\\
		Apply the chosen initial descent condition approach:
		\begin{enumerate}[label={\alph*)},itemsep=0.3\baselineskip]
			\item Set $\J_{\DC}(t_0)=\boldsymbol{\infty}$.
			\item Set $\realize{\xSys}^0 = \xSys_{\meas}(t_0)$, $\realize{\J}_{\ub} = \boldsymbol{\infty}$.
			Set $\decSetFeas=\decSetFeasMOMPC\big(\realize{\xSys}^0,\realize{\J}_{\ub}\big)$
			and compute the regularized pay-off matrix $\Phi_{\mathrm{reg}}(\boldsymbol{\delta})$ (cf.~\eqref{eq:Phi_reg}).
			Then, set $\J_{\DC}(t_0)=\J_{\NP,\mathrm{reg}}(\boldsymbol{\delta})$ (row-wise maximum of $\Phi_{\mathrm{reg}}(\boldsymbol{\delta})$).
		\end{enumerate}	
		\item \label{step:algo_IM_all} 
		Set $\realize{\xSys}^0 = \xSys_{\meas}(t_k)$ and 
		$\realize{\J}_{\ub} = \J_{\DC}(t_k)$.
		\\[0.5\baselineskip]
		$
			\text{Then solve}
			\left\{
				\begin{alignedat}{2}
					&\text{\WScount} \ &&\text{\ref{eq:pOP_WS}} \\
					&\text{\PScount} \ &&\text{\ref{eq:pOP_PS}} \\
				\end{alignedat}
			\right.
			\text{ 
				with $p$ as
			 	in
		 	}
			\left\{
				\begin{alignedat}{2}
					&\text{\WScount} \ && \text{\eqref{eq:p_WS_IM}} \\
					&\text{\PScount} \ && \text{\eqref{eq:p_PS_IM}}
				\end{alignedat}
			\right.
		$
		\vspace{0.3\baselineskip}\\
		and $\decSetFeas=\decSetFeasMOMPC\big(\realize{\xSys}^0,\realize{\J}_{\ub}\big)$ 
		for all $i\in\{1,\ldots,n_J\}$ to obtain the \accp{IM} solutions.
		\item Derive the quantities from \cref{list:objects_from_Phi} based on $\Phi$. 
		\\[0.2\baselineskip]
		$
			\text{Set }
			\left\{
				\begin{alignedat}{2}
					&\text{\WScount} \ && \realize{w}_{\mathrm{WS}} \\
					&\text{\PScount} \ && \realize{\J}_{\SO}, \, \realize{d} \\
				\end{alignedat}
			\right.
			\
		$
		depending on the chosen method.
		\item 
		\label{step:algo_MOMPC_beta}
		Run an algorithm that computes a preference vector $\betaPref(t_k)$.
		\item 
		\label{step:algo_FHP_all}
		$
			\text{Solve }
			\left\{
				\begin{alignedat}{2}
					&\text{\WScount} \ &&\text{\ref{eq:pOP_WS}} \\
					&\text{\PScount} \ &&\text{\ref{eq:pOP_PS}}
				\end{alignedat}
			\right.
			\text{ 
				with $p$ as
			 	in
		 	}
			\left\{
				\begin{alignedat}{2}
					&\text{\WScount} \ &&\text{\eqref{eq:p_WS_FHP}} \\
					&\text{\PScount} \ &&\text{\eqref{eq:p_PS_FHP}}
				\end{alignedat}
			\right.
		$
		\vspace{0.3\baselineskip}\\
		and $\decSetFeas=\decSetFeasMOMPC\big(\realize{\xSys}^0,\realize{\J}_{\ub}\big)$
		to obtain the optimal solution $\xMO\opt$ and the corresponding optimal cost vector $\JStarTime:=\JMO(\xMO\opt)$.
		
		\item
		\label{step:MOMPC_terminate}
		If $\JStarTime \leq \J_{\mathrm{conv}}$, \textbf{terminate} the algorithm.
		
		\item 
		\label{step:feedback_and_update_all}
		Apply the feedback $\uSys_{0}\opt$ to the system to obtain
		$\xSys_{\meas}(t_{k+1}) = \xSys_{1}\opt$ (nominal case).
		Next, set
		$\J_{\DC}(t_{k+1}) = \JStarTime$.
		\\
		Then, perform the update $k\rightarrow k+1$, $t_k\rightarrow t_{k+1}=t_k+h$ and go to step \ref{step:algo_IM_all}.
	\end{steps}
	\vspace{-0.7\baselineskip}
\end{algorithm}

\Cref{alg:MOMPC} defines an \acc{MOMPC} scheme that makes real-time implementation more plausible in practice%
\footnote{
	If the solver used does not accept infinite bounds, $\boldsymbol{\infty}$ in \cref{step:MOMPC_init} of \cref{alg:MOMPC} can be replaced by a vector of sufficiently large numbers.
}%
.
In contrast to a standard single-objective \acc{MPC} scheme, \cref{alg:MOMPC} requires two sequential optimization stages.
The first stage consists of $n_J$ independent optimization problems for the computation of the \accp{IM}, which naturally provide the basis for the subsequent \acc{DM} step and can be solved in parallel.
Moreover, since these problems are structurally similar across sampling instants, they are also naturally amenable to warm-starting.
The second stage then consists of one additional optimization problem that realizes the \acc{DM} and can, in principle, be warm-started from the solution obtained at the previous sampling instant.

The \accp{IM}-informed methods proposed in this work do not primarily aim at accelerating the individual optimization problems themselves, for example through tailored problem transcriptions, hyperparameter choices, or specialized solvers.
Instead, their purpose is to generate, with as few sequential optimization stages as possible, a sufficiently informative basis for an automated and informed \acc{DM}.
At the same time, the size of the parallel workload should remain moderate, so that the resulting computational burden remains compatible with standard computing hardware.
Both requirements are addressed by the proposed scheme, which requires a total of only $n_J+1$ optimization problems organized into two sequential optimization stages.

\subsection{Stabilizing MOMPC Scheme}
\label{ssec:MOMPC_Scheme}
Having clarified the intended computational structure of \cref{alg:MOMPC}, the focus now turns to its theoretical properties.
To prove that \cref{alg:MOMPC} stabilizes a (controlled) fixed point $\xSys\cSS$, the following assumptions from \cite{grune_multiobjective_2019} are needed.
The notation is adapted accordingly.
\goodbreak
\begin{bThm}[Classical stage costs]{assumption}
	\label{asm:classical_stage_costs}
	\phantom{a}\\[-1\baselineskip]
	\begin{enumerate}[label={\arabic*.},labelindent=0pt]
		\item There exists an equilibrium pair $(\xSys\cSS, \uSys\cSS) \in \xSet \times \uSet$, i.e., $f(\xSys\cSS, \uSys\cSS)=\xSys\cSS$.
		\item There are comparison functions $\alpha_{\ell, i} \in \mathcal{K}$ (continuous, strictly increasing functions with $\alpha(0)=0$, cf.~\cite{kellett_compendium_2014})	such that all stage costs $\ell_i$, $i \in\{1, \ldots,\nObj\}$, satisfy 
		$\min_{\uSys \in \uSet} \, \ell_i(\xSys, \uSys) \geq \alpha_{\ell, i}\(\left\|\xSys-\xSys\cSS\right\|\) \forall \, \xSys \in \xSet$.
	\end{enumerate}
\end{bThm}

\begin{bThm}[Lyapunov function terminal cost]{assumption}
	\label{asm:Lyap_term_cost}
	Let $\xSys\cSS$ be an interior point of $\xSet\term \subseteq \xSet$, and let the terminal costs $F_i: \xSet\term \to \Rset_{\geq 0}$ satisfy $F_i(\xSys) \leq \alpha_{F}\(\left\|\xSys-\xSys\cSS\right\|\)$ for all $\xSys \in \xSet\term$, $i \in\{1, \ldots, \nObj\}$, and some $\alpha_{F} \in \mathcal{K}$.
	The existence of a local feedback controller $\termContr: \xSet\term \rightarrow \uSet$ is assumed, satisfying
	\begin{enumerate}[label={\arabic*.},labelindent=0pt]
		\item $f(\xSys, \termContr(\xSys)) \in \xSet\term$ for all $\xSys \in \xSet\term$ and
		\item%
		$\forall \, \xSys \in \xSet\term, i \in\{1, \ldots, \nObj\}\! :$
		\if\TwoCol1
			\\
			\quad
		\else
		\fi
		$F_i\(f(\xSys, \termContr(\xSys))\)+\ell_i\(\xSys, \termContr(\xSys)\) \leq F_i(\xSys)$.
	\end{enumerate}
\end{bThm}
Both \cref{asm:classical_stage_costs} and \cref{asm:Lyap_term_cost} are  straightforward extensions of those applied to stability proofs in the single-objective case.
The quadratic terminal ingredients constructed in \cref{ssec:terminal_ingredients} satisfy these requirements.

One more set needs to be defined due to its initial state dependency.
\goodbreak
\begin{bThm}[Feasible initial states]{definition}
	\label{def:feas_initial_states}
	The set of feasible initial states is defined by
	\vspace{-1ex}
	\begin{equation}
		\xSet\feasInitialStates:=
		\left\{
			\param{\xSys}^0 \in \xSet 
			\ \middle| \ 
			\decSetFeasMPC(\param{\xSys}^0) \neq \emptyset
		\right\}
		.
	\end{equation}
\end{bThm}
Using the set
\text{\hyperref[def:feas_initial_states]{$\xSet\feasInitialStates$}},
the following theorem can now be formulated.
\begin{bThm}[Stabilizing Scheme]{theorem}
	\label{thm:conv_proof}
	Consider the \acc{MOMPC} problems \ref{eq:pOP_WS} and \ref{eq:pOP_PS}.
	Let
	\begin{enumerate}[label={\roman*)},itemsep=0em,parsep=0em,labelindent=0pt]
		\item $N \geq 2$; 
		\item
		\cref{asm:classical_stage_costs} and \cref{asm:Lyap_term_cost} hold;
		\item the set $\decSetFeasMPC(\param{\xSys}^0)$ be compact for each feasible initial state
		$\param{\xSys}^0 \in \text{\hyperref[def:feas_initial_states]{$\xSet\feasInitialStates$}}$
		and $\ell_i$, $F_i$ be continuous;
		\item \cref{step:algo_FHP_all} of \cref{alg:MOMPC} yield a solution that is not strictly dominated by the candidate solution $\uSeq_{\xSys(\n+1)} := \big(\uSys_1\opt,\ldots,\uSys_{N-1}\opt,\termContr(\xSys_N\opt)\big)$ and, at $t_0$ with $\xSys_0 \in \xSet\term$, not by the input sequence generated by $\termContr$.
	\end{enumerate}
	Then, $\xSys\cSS$ is asymptotically stable for the closed loop defined by \cref{alg:MOMPC} without the termination criterion in \cref{step:MOMPC_terminate}, for both initial descent condition approaches and any sequence of preference vectors $\betaPref(t_k)$.
	That is, for each $\varepsilon > 0$ there is $\rho > 0$ such that $\|\xSys_0-\xSys\cSS\| < \rho$ implies $\|\xSys(\n)-\xSys\cSS\| < \varepsilon$ for all $\n \in \mathbb{N}_0$, and $\xSys(\n) \to \xSys\cSS$ for each $\xSys_0 \in \xSet\feasInitialStates$.
\end{bThm}
\newcommand{\PH}[1]{}
\begin{proof}
	The notation of \cite{grune_multiobjective_2019} is used, which is summarized in \cref{sec:conv_proof}, and $J_i^\star(\n) := J_i^{\PH{N}}\big(\xSys(\n), \uSeq_{\xSys(\n)}^{\star \PH{, N}}\big)$ abbreviates the cost of the solution applied at time $\n$.
	Due to the update in \cref{step:feedback_and_update_all} of \cref{alg:MOMPC}, the descent condition
	\begin{equation}
		J_i^\star(\n) \leq J_i^\star(\n-1)
		\quad \forall \, i \in \{1,\ldots,\nObj\}, \ \n \in \mathbb{N}
		\label{eq:descent_condition}
	\end{equation}
	holds.
	This descent condition is less restrictive than the one in \cite[Algorithm~2, step~(1)]{grune_multiobjective_2019}, where the upper bound is the cost $J_i^{\PH{N}}\big(\xSys(\n), \uSeq_{\xSys(\n)}^{\PH{N}}\big)$ of the candidate solution instead of $J_i^\star(\n-1)$.
	\\[0.3\baselineskip]
	\textbf{Feasibility:} 
	The sets $\decSetFeasMOMPC\big(\xSys(\n),\J_{\DC}(t_\n)\big)$ used in steps \ref{step:algo_IM_all} and \ref{step:algo_FHP_all} of \cref{alg:MOMPC} are nonempty at $t_0$ for $\xSys_0 \in \xSet\feasInitialStates$ and, as shown below, at all subsequent time steps.
	By iii), they are compact, since the upper bound on $\JMO$ defines a closed constraint.
	Hence, \ref{eq:pOP_WS} attains its minimum for both tasks, and so does \ref{eq:pOP_PS} for the computation of the \accp{IM}.
	In step \ref{step:algo_FHP_all}, the existence of a solution of \ref{eq:pOP_PS} is part of iv); it holds, e.g., for $\realize{d}<0$, because \ref{eq:pOP_PS} is then feasible and the ray-scaling variable $l$ is implicitly bounded from above due to $\J(\xMO) \geq \boldsymbol{0}$.
	For the candidate solution $\uSeq_{\xSys(\n+1)}^{\PH{N}}$ from iv), the resulting states are $\xSys_1\opt,\ldots,\xSys_N\opt,f\big(\xSys_N\opt,\termContr(\xSys_N\opt)\big)$.
	They satisfy \eqref{eq:dyn_constraints} with $\param{\xSys}^0=\xSys(\n+1)$ by item~1 of \cref{asm:Lyap_term_cost}.
	Hence, item~2 of \cref{asm:Lyap_term_cost} at $\xSys_N\opt \in \xSet\term$ and $\xSys_0\opt = \xSys(\n)$ yield (cf.~\cite[Lemma~3.2 and proof of Theorem~3.3]{grune_multiobjective_2019})
	\begin{equation}
		\begin{aligned}
			J_i^{\PH{N}}\(\xSys(\n+1), \uSeq_{\xSys(\n+1)}^{\PH{N}}\)
			&= \sum_{j=1}^{N-1} \ell_i\(\xSys_j\opt,\uSys_j\opt\)
			\proofColDepAlignBefore \proofColDepAfter
			+ \ell_i\(\xSys_N\opt,\termContr(\xSys_N\opt)\)
			+ F_i\(f\big(\xSys_N\opt,\termContr(\xSys_N\opt)\big)\)
			\\
			&\leq \sum_{j=1}^{N-1} \ell_i\(\xSys_j\opt,\uSys_j\opt\) + F_i\(\xSys_N\opt\)
			\\
			&= J_i^\star(\n) - \ell_i\(\xSys(\n), \uSys_{\xSys(\n)}^{\star \PH{, N}}(0)\)
			\proofColDepAlignBefore \proofColDepAfter
			\quad \forall \, i \in \{1,\ldots,\nObj\}.
		\end{aligned}
		\label{eq:candidate_decrease}
	\end{equation}
	Since $\ell_i \geq 0$, the candidate also satisfies \eqref{eq:descent_condition}, which implies recursive feasibility.
	\\[0.3\baselineskip]
	\textbf{Stability:} 
	Let $\alpha_\ell := \min_i \alpha_{\ell,i} \in \mathcal{K}$ \cite{kellett_compendium_2014} and $\xSys_0 \in \xSet\term$.
	By \cref{asm:Lyap_term_cost}, the input sequence generated by $\termContr$ is feasible and its cost is bounded by $F_i(\xSys_0)$, which follows from summing item 2 along the resulting trajectory.
	Thus, iv) implies the existence of an index $m$ with $J_m^\star(0) \leq F_m(\xSys_0)$.
	Together with the nonnegativity of $\ell_i$ and $F_i$ and \eqref{eq:descent_condition}, this yields
	\begin{equation}
		\begin{aligned}
			\alpha_\ell\(\|\xSys(\n)-\xSys\cSS\|\)
			&\leq
			\ell_m\(\xSys(\n), \uSys_{\xSys(\n)}^{\star \PH{, N}}(0)\)
			\leq
			J_m^\star(\n)
			\proofColDepAlignBefore \proofColDepAfter
			\leq
			J_m^\star(0)
			\leq
			\alpha_F\(\|\xSys_0-\xSys\cSS\|\)
			\quad \forall \, \n \in \mathbb{N}_0.
		\end{aligned}
	\end{equation}
	For $\varepsilon > 0$, choose $\rho > 0$ such that $\alpha_F(\rho) < \alpha_\ell(\varepsilon)$ and all $\xSys_0$ with $\|\xSys_0-\xSys\cSS\| < \rho$ lie in $\xSet\term$, which proves stability.
	\\[0.3\baselineskip]
	\textbf{Attractivity:} 
	Let $V(\n) := \sum_{i=1}^{\nObj} J_i^\star(\n) \geq 0$.
	By iv), there is an index $m=m(\n)$, which may change with $\n$, such that $J_m^\star(\n+1) \leq J_m^{\PH{N}}\big(\xSys(\n+1), \uSeq_{\xSys(\n+1)}^{\PH{N}}\big)$.
	Using \eqref{eq:candidate_decrease} for $m$ and \eqref{eq:descent_condition} for $i\neq m$ yields
	\begin{align}
		V(\n+1)
		&= J_m^\star(\n+1) + \sum_{i\neq m} J_i^\star(\n+1)
		\nonumber
		\\
		\overset{\text{iv)}}&{\leq} J_m^{\PH{N}}\(\xSys(\n+1), \uSeq_{\xSys(\n+1)}^{\PH{N}}\)
		\nonumber \proofColDepAlignBefore \proofColDepAfter
		+ \sum_{i\neq m} J_i^\star(\n+1)
		\nonumber
		\\
		\overset{\eqref{eq:candidate_decrease},\eqref{eq:descent_condition}}&{\leq} J_m^\star(\n) - \ell_m\(\xSys(\n), \uSys_{\xSys(\n)}^{\star \PH{, N}}(0)\)
		\nonumber \proofColDepAlignBefore \proofColDepAfter
		+ \sum_{i\neq m} J_i^\star(\n)
		\nonumber
		\\
		&= V(\n) - \ell_m\(\xSys(\n), \uSys_{\xSys(\n)}^{\star \PH{, N}}(0)\)
		\nonumber\\
		\overset{\text{Asm.~\ref{asm:classical_stage_costs}}}&{\leq} V(\n) - \alpha_\ell\(\|\xSys(\n)-\xSys\cSS\|\).
	\end{align}
	Summation gives 
	$
		\sum_{\n=0}^{\kMax-1} \alpha_\ell \big(\|\xSys(\n)-\xSys\cSS\|\big) 
		\leq V(0) - V(K) 
		\if\TwoCol1
			\linebreak
		\fi
		\leq V(0)
	$ 
	for all $\kMax \in \mathbb{N}$.
	Hence, $\alpha_\ell \big(\|\xSys(\n)-\xSys\cSS\|\big) \to 0$ and thus $\xSys(\n) \to \xSys\cSS$, cf.~\cite[Corollary~3.6]{grune_multiobjective_2019}.
\end{proof}
\begin{bThm}[Relation to \cite{grune_multiobjective_2019}]{remark}
	\label{rem:proof}
	\Cref{thm:conv_proof} reuses the terminal ingredients of \cite[Section~3]{grune_multiobjective_2019} and the candidate solution in iv), which corresponds to \cite[Algorithm~2, step~(2)]{grune_multiobjective_2019}, but employs a less restrictive descent condition.
	As a consequence, the individual performance bounds of \cite[Theorem~3.3]{grune_multiobjective_2019} are not available.
	Instead, a decrease of the sum of all objectives is exploited, which is sufficient for convergence, and the monotonicity of each objective additionally yields Lyapunov stability.
	Since a decrease by the stage cost is only guaranteed for one objective, which may change over time, \cref{asm:classical_stage_costs} requires each stage cost to be positive definite, which is stronger than \cite[Assumption~3.5]{grune_multiobjective_2019}.
	Condition iii) implies the external stability of the Pareto front assumed in \cite[Theorem~3.3]{grune_multiobjective_2019} (cf.~\cite[Theorem~2.4 and Lemma~2.5]{grune_multiobjective_2019}) and thus replaces this assumption.
	In addition, condition iii) ensures that the scalarized \accp{OP} attain their minima.
	\\
	Condition iv) is not part of \cite{grune_multiobjective_2019}, but is required here due to the less restrictive descent condition.
	Since the costs at $t_{k+1}$ are bounded by the optimal cost vector $\JStarTime$ of the previous time step (cf.~\cref{step:feedback_and_update_all} of \cref{alg:MOMPC}) instead of by the cost vector of the candidate solution, the computed solution is no longer guaranteed to be at least as good as the candidate in every objective.
	Apart from the existence of a solution, condition iv) requires the computed solution to be at least as good as the candidate in at least one objective, which is sufficient for both stability and attractivity.
	Condition iv) refers to the solution that is actually computed and cannot be guaranteed a priori if the \accp{OP} are solved by a local solver.
	However, condition iv) can be enforced for any of the \acc{DM} methods with one additional evaluation of the candidate's cost vector, by applying the candidate whenever \cref{step:algo_FHP_all} yields no solution or one that is strictly dominated by the candidate (at $t_0$, the input sequence generated by $\termContr$ or, if $\xSys_0 \notin \xSet\term$, an \acc{IM} solution from \cref{step:algo_IM_all}).
	Without this fallback, iv) holds for global solutions of \ref{eq:pOP_PS} and of \ref{eq:pOP_WS} with $\realize{w}_{\mathrm{WS}} \in \Rset_{\geq 0}^{\nObj} \setminus \{\boldsymbol{0}\}$, since these are weakly efficient \cite[Theorem~3.1.1]{miettinen_nonlinear_1998}, \cite[Theorem~2.1c)]{eichfelder_adaptive_2008}.
	This also holds with respect to the \acc{OP} without the upper bound $\J_{\DC}(t_k)$, because any solution that is strictly better in all objectives satisfies this upper bound as well.
	For global solutions, the \acc{KP} method \ref{enum:WS_2} is only covered if $-\wChim$ has no negative components (cf.~\cref{sec:KP}).
\end{bThm}
\begin{bThm}[Role of the \accp{IM}-informed \acc{DM}]{remark}
	\label{rem:IM_DM}
	Close to $\xSys\cSS$, the \accp{IM} may coincide, such that the quantities derived from them degenerate.
	In this case, any solution that satisfies iv), e.g., the candidate solution, can be applied without affecting \cref{thm:conv_proof}.
	In practice, the termination criterion in \cref{step:MOMPC_terminate} of \cref{alg:MOMPC} avoids this situation.
	\\
	The stability proof does not rely on geometric properties of the \accp{IM} or on the quality of the preference translation induced by the chosen \acc{DM} method.
	These aspects only affect which point on the \acc{PF} is selected within the feasible set.
	The less restrictive descent condition used here enlarges the feasible set and thereby preserves a broader set of \accp{ES} for the \acc{DM} process.
	Among all descent conditions that prevent any objective from increasing compared to the previous time step, \eqref{eq:descent_condition} is the least restrictive one.
\end{bThm}

\subsection{A Practical Approach to Determine Stabilizing Terminal Ingredients}
\label{ssec:terminal_ingredients}
Lastly, stabilizing terminal ingredients are explicitly constructed as required in \cref{thm:conv_proof} for a simplified setting.
As it has been proposed in
\cite{chen_quasi-infinite_1998,chen_terminal_2003,johansen_approximate_2004,gupta_application_2024} (a.o.), as well, 
let $\Delta \xSys = \xSys - \xSys\cSS$, $\Delta \uSys = \uSys - \uSys\cSS$, $i \in \{1,\ldots,n_J\}$, and assume
\goodbreak
\begin{enumerate}[label={(S\arabic*)}]
	\item \label{item:S1}
		quadratic stage costs 
		\\
		$\ell_i(\xSys,\uSys) = \Delta\xSys\tr \Q[][,i] \Delta\xSys + \Delta\uSys\tr \R[][,i] \Delta\uSys$ with $\Q[][,i],\R[][,i] \succ 0$, 
	\item \label{item:S2}
		quadratic terminal costs $\termCostFun_i(\xSys)=\Delta\xSys\tr P \Delta\xSys$ ($P \succ 0$),
	\item \label{item:S3}
		the terminal set to be defined as $\xSet\term = \{\xSys \in \Rset^{n_{\xSys}} \ | \ \Delta\xSys\tr P \Delta\xSys \leq \alpha\}$, 
		where $\alpha\geq0$ is the \emph{worst-case} terminal cost,
	\item \label{item:S4}
		the existence of a linear local feedback controller $\termContr\expkLoc(\xSys)=\uSys\cSS + K \Delta\xSys$.
\end{enumerate}

Based on the assumptions \ref{item:S1}-\ref{item:S4}, a straightforward procedure is proposed to determine the terminal cost matrix $P$, the controller gain $K$, and the terminal set parameter $\alpha$.
\\
The system matrices $A$ and $B$ are obtained by linearizing the system around $(\xSys\cSS, \uSys\cSS)$.
Let $\Q[\bar][]$ and $\R[\bar][]$ denote the arithmetic means of $\Q[][,i]$ and $\R[][,i]$, $i \in \{1,\ldots,n_J\}$, respectively.
Following the LQR approach shown in \cite{gupta_application_2024}, the gain $K := -K_{\mathrm{LQR}}$ and the matrix $P$ are obtained from the discrete-time LQR problem for $(A,B)$ with the weighting matrices $\rho_{\xSys} \Q[\bar][]$ and $\rho_{\uSys} \R[\bar][]$, where $\rho_{\xSys}, \rho_{\uSys} > 0$.
Hence, $P$ solves the discrete-time Lyapunov equation $A_K\tr P A_K - P = -\Omega$ with $A_K = A + BK$ and $\Omega := \rho_{\xSys} \Q[\bar][] + \rho_{\uSys} K\tr \R[\bar][] K$.
Let $\Qstar[][i] := \Q[][,i] \, + \, K\tr \R[][,i] K$, $i \in \{1,\ldots,n_J\}$.
By defining
\begin{equation}
	\begin{aligned}
			\DeltaQ_{i} &:= \Omega - \Qstar[][i],
			\quad i=1,\ldots,n_J,
	\end{aligned}
\end{equation}
it is ensured that the terminal cost matrix $P$ satisfies
$A_K\tr P A_K - P = -(\Qstar[][i] + \DeltaQ_{i})$
for all $i \in \{1,\ldots,n_J\}$.
With $\Qstar[\bar][] := \Q[\bar][] + K\tr \R[\bar][] K$, it holds that $\Qstar[][i] \preceq n_J \Qstar[\bar][]$ and $\Omega \succeq \min\{\rho_{\xSys},\rho_{\uSys}\} \Qstar[\bar][]$.
Hence, $\DeltaQ_{i} \succ 0$ is guaranteed for $\min\{\rho_{\xSys},\rho_{\uSys}\} > n_J$.

Now, let $\tilde{\alpha}_i$ be the largest possible realization of $\alpha_i$ such that 
\begin{equation}
    \begin{split}
		F_i\(f(\xSys, \termContr(\xSys))\) - F_i(\xSys) 
		&= -\Delta\xSys\tr(\Qstar[][i] + \DeltaQ_{i})\Delta\xSys + q(\xSys)
		\\
		&\leq -\ell_i\(\xSys, \termContr(\xSys)\)
		= -\Delta\xSys\tr \Qstar[][i] \Delta\xSys
    \end{split}
\end{equation}
with the nonlinear function $q(\xSys)$ depending on the dynamics, $K$ and $P$,
holds for all $\xSys \in \{\xSys \in \Rset^{n_{\xSys}} \ | \ \Delta\xSys\tr P \Delta\xSys \leq \alpha_i\}$, $i \in \{1,\ldots,n_J\}$.
Then, $\alpha = \min \{\tilde{\alpha}_1,\ldots,\tilde{\alpha}_{n_J}\}$ is a feasible choice (cf.~\cref{asm:Lyap_term_cost}) and all required terminal ingredients are determined.

\subsection{Adaptation of the Preference Vector}
\label{ssec:Adaptation_Preference_Vector}
In \cref{step:algo_MOMPC_beta} of \cref{alg:MOMPC}, a preference vector must be computed using all information available up to the current time step.
While the most straightforward approach is to apply a constant preference vector $\betaPref(t_k) = \beta_{\mathrm{fixed}}$ for all $k$, the proposed \acc{MO} framework allows for the dynamic, iterative adaptation of $\betaPref$.
This represents a significant advantage over classical single-objective or fixed-weight \acc{WS} formulations.

To illustrate this, \cref{alg:waterfilling_beta_update} is employed, which is presented and discussed in detail in \cref{sec:Water_Filling_Strategy}.
The algorithm implements a \emph{water-filling strategy} based on the classical projection method described in \cite{michelot_finite_1986}; see also \cite[Chapter 5.5.3]{boyd_convex_2004} and \cite[Chapter 10.4]{cover_elements_2001}.
Conceptually, this specific water-filling strategy can be interpreted as aiming to drive all cost components toward convergence as quickly as possible over time.
Enforcing synchronization across all cost components acts as an additional restrictive constraint.
Nevertheless, dynamically distributing the \enquote{preference budget} among the currently highest costs yields a fast overall decrease if the cost reductions scale proportionally with the assigned preferences.

\section{Numerical Example}
\label{sec:num_ex}
This section utilizes an exemplary room climate model to illustrate the proposed method.
The case study evaluates the complete closed-loop scheme, i.e., the \accp{IM}-informed \acc{DM} with online preference adaptation within the stabilizing \acc{MOMPC} scheme, in terms of convergence time. 
The model describes the coupled dynamics of temperature, humidity, and CO$_2$ concentration in a single room equipped with floor heating and a combined HVAC system (heating, ventilation, air conditioning) as well as humidity control.

The goal is to stabilize the controlled fixed point pair $(\xSys\cSS, \uSys\cSS)$
using quadratic stage costs as depicted in \ref{item:S1} from \cref{ssec:terminal_ingredients}.
The modeling of the system dynamics, the derivation of a financially optimal controlled fixed point pair $(\xSys\cSS, \uSys\cSS)$ 
and the computation of suitable terminal ingredients $\termContr(\xSys)$, $\termCostFun_i(\xSys)$, $i \in \{1,\ldots,n_J\}$, $\xSet\term$ can be found in \cref{sec:model_num_ex}.

The numerical evaluation considers eight test cases that result from
four different initial states (\circled{1}-\circled{4}) 
combined with 
two different sets of external influences (\circled{A},\circled{B}) on the system dynamics. 
Details can again be found in \cref{sec:model_num_ex}.

For \cref{alg:MOMPC}, the parameters $\boldsymbol{\delta}=\delta=0.002$, $\J_{\mathrm{conv}}= 10^{-2} \cdot \boldsymbol{1}$ are chosen, and both initial descent condition approaches are tested.
For \cref{alg:waterfilling_beta_update}, $\alpha=0.05$ is chosen.

\subsection{Challenges of the Problem}
\label{ssec:Challenges_of_Problem}
In order to get to know the special characteristics of the specific application, 
the test case \testSetup{A}{1} is examined.
Therefore, the first iteration of the \acc{MOMPC} scheme (\cref{alg:MOMPC}) is considered, for which $\realize{\J}_{\ub} = \J_{\DC}(t_0)=\boldsymbol{\infty}$ was set in \cref{step:MOMPC_init}, i.e.~no regularization was applied.

\acc{PF} samples, obtained via different scalarization methods, are given in \cref{fig:PF_first_iter}.
This investigation focuses on the properties of the \acc{PF} and of the \acc{CHIM} normal.

\begin{figure}[tb]
	\centering
		\adjincludegraphics[
			width=0.95\WidthBase, clip, 
			trim={{0\width} {0\height} {0.00\width} {0\height}}
		]{PF_first_iter_new.\FigureFileExt} 
	\caption{%
		\Acc{PF} samples of test case
		\testSetup{A}{1}
		in the first iteration of the \acc{MOMPC} scheme obtained via different scalarization methods.
		Each method computes $3500$ \acc[l]{PF} samples obtained via parameterizations that aim at an uniform sampling of the unit simplex (all methods except \acc{SRI}) and a hypersphere segment (\acc{SRI}, cf.~\cref{fig:SRI_rays_coverage}).
		The \acc[l]{PF} was reconstructed using the \acc{SRI} samples and the ball pivoting algorithm.
	}
	\label{fig:PF_first_iter}
\end{figure}

From first visual observation, one can conclude that the \acc{PF} is nonconvex. 
This is in agreement with the \acc{WS} method not being able to sample in large parts of the \acc{PF}.
However, similar behavior is observed for the \acc{NCHIM} method and the \acc{NBI} method.
Still, it is noted that the \acc{NBI} method with $\realize{d}=\wChim$ provides by far the worst sampling.
Thirdly, from the position of the \accp{IM}, the \acc{CHIM} normal can be computed as $\wChim \approx \matInline{0.024, -0.968, 0.008}\tr$,
or $\wChim \approx \matInline{-0.494, 0.427, -0.079}\tr$ in the \acc{NIS}.
\Cref{fig:DM_all_methods} (right) and \cref{fig:PF_first_iter} are compared with their corresponding \acc{CHIM} normal vectors (in a \acc{NIS}).
This comparison shows that the \acc{CHIM} normal vectors alone do not indicate how well the \acc{DM} methods perform in terms of \emph{coverage} and \emph{translation}.

\begin{figure}[tb]
	\centering
		\adjincludegraphics[
			width=0.95\WidthBase, clip, 
			trim={{0\width} {0\height} {0\width} {0\height}}
		]{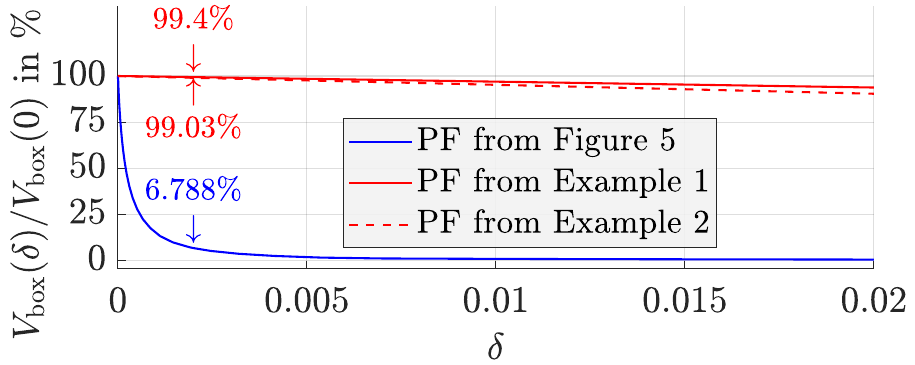} 
	\caption{Change of the hyperbox volume. For comparison, the curves corresponding to the \accp{PF} from \cref{example1} and \cref{example2} exhibit an approximately linear behavior. The percentage values are given for $\delta=0.002$.}
	\label{fig:pathological_PF_box_volume}
\end{figure}

In summary, for the test case shown in \cref{fig:PF_first_iter}, 
and some of the other test cases,
examples of \emph{pathological \accp{PF}} as defined in \cref{defi:pathological_PF} could be observed.
The approximately exponential decrease of the hyperbox size can be deduced from \cref{fig:pathological_PF_box_volume}.
\\
With such pathological \accp{PF}, the \accp{IM} seem to contain less information about the complete \acc{PF}.
It is noted that using automated \accp{IM}-informed \acc{DM} methods for \acc{MOMPC} on pathological \accp{PF} might lead to a deterioration in the translation of a high-level preference to a point on the \acc{PF}.
However, regardless of how poor the translation might be in such a potential case, the stability of the controller is not compromised.

\subsection{Analysis of the Numerical Results}
\label{ssec:Analysis_Numerical_Results}
\begin{table*}[htb]
	\centering
	\definecolor{cC1}{rgb}{0.4378,0.6159,0.1335}
\definecolor{cD1}{rgb}{0.6899,0.1860,0.1337}
\definecolor{cE1}{rgb}{0.0000,0.7490,0.1294}
\definecolor{cF1}{rgb}{0.0000,0.7490,0.1294}
\definecolor{cG1}{rgb}{0.3216,0.6834,0.1318}
\definecolor{cC2}{rgb}{0.6673,0.2794,0.1343}
\definecolor{cD2}{rgb}{0.6980,0.1333,0.1333}
\definecolor{cE2}{rgb}{0.3420,0.6739,0.1321}
\definecolor{cF2}{rgb}{0.0000,0.7490,0.1294}
\definecolor{cG2}{rgb}{0.4507,0.6061,0.1337}
\definecolor{cC3}{rgb}{0.5767,0.4715,0.1350}
\definecolor{cD3}{rgb}{0.4630,0.5962,0.1339}
\definecolor{cE3}{rgb}{0.0000,0.7490,0.1294}
\definecolor{cF3}{rgb}{0.2458,0.7117,0.1309}
\definecolor{cG3}{rgb}{0.0000,0.7490,0.1294}
\definecolor{cC4}{rgb}{0.5070,0.5561,0.1345}
\definecolor{cD4}{rgb}{0.4630,0.5962,0.1339}
\definecolor{cE4}{rgb}{0.2458,0.7117,0.1309}
\definecolor{cF4}{rgb}{0.0000,0.7490,0.1294}
\definecolor{cG4}{rgb}{0.0000,0.7490,0.1294}
\definecolor{cC5}{rgb}{0.2742,0.7023,0.1312}
\definecolor{cD5}{rgb}{0.4507,0.6061,0.1337}
\definecolor{cE5}{rgb}{0.0000,0.7490,0.1294}
\definecolor{cF5}{rgb}{0.2742,0.7023,0.1312}
\definecolor{cG5}{rgb}{0.0000,0.7490,0.1294}
\definecolor{cC6}{rgb}{0.3782,0.6547,0.1327}
\definecolor{cD6}{rgb}{0.4507,0.6061,0.1337}
\definecolor{cE6}{rgb}{0.0000,0.7490,0.1294}
\definecolor{cF6}{rgb}{0.0000,0.7490,0.1294}
\definecolor{cG6}{rgb}{0.0000,0.7490,0.1294}
\definecolor{cC7}{rgb}{0.3945,0.6451,0.1329}
\definecolor{cD7}{rgb}{0.2991,0.6928,0.1315}
\definecolor{cE7}{rgb}{0.0000,0.7490,0.1294}
\definecolor{cF7}{rgb}{0.0000,0.7490,0.1294}
\definecolor{cG7}{rgb}{0.2991,0.6928,0.1315}
\definecolor{cC8}{rgb}{0.5841,0.4604,0.1350}
\definecolor{cD8}{rgb}{0.5355,0.5251,0.1348}
\definecolor{cE8}{rgb}{0.0000,0.7490,0.1294}
\definecolor{cF8}{rgb}{0.2991,0.6928,0.1315}
\definecolor{cG8}{rgb}{0.0000,0.7490,0.1294}

\def\rowA{& 0.76 & 0.81 & 0.48 & 0.48 & 1.00 & 0.96 & 0.57 & 0.60}
\def\rowB{& 5.00 & 5.31 & 4.83 & 4.79 & 2.92 & 3.00 & 4.76 & 5.00}
\def\rowC{& \SetCell{bg=cC1} 1.12 & \SetCell{bg=cC2} 1.38 & \SetCell{bg=cC3} 1.24 & \SetCell{bg=cC4} 1.17 & \SetCell{bg=cC5} 1.04 & \SetCell{bg=cC6} 1.09 & \SetCell{bg=cC7} 1.10 & \SetCell{bg=cC8} 1.25 & 1.1725 & 0.1103}
\def\rowD{& \SetCell{bg=cD1} 1.41 & \SetCell{bg=cD2} 1.44 & \SetCell{bg=cD3} 1.14 & \SetCell{bg=cD4} 1.14 & \SetCell{bg=cD5} 1.12 & \SetCell{bg=cD6} 1.13 & \SetCell{bg=cD7} 1.05 & \SetCell{bg=cD8} 1.20 & 1.2035 & 0.1427}
\def\rowE{& \SetCell{bg=cE1} 1.00 & \SetCell{bg=cE2} 1.06 & \SetCell{bg=cE3} 1.00 & \SetCell{bg=cE4} 1.03 & \SetCell{bg=cE5} 1.00 & \SetCell{bg=cE6} 1.00 & \SetCell{bg=cE7} 1.00 & \SetCell{bg=cE8} 1.00 & 1.0121 & 0.0237}
\def\rowF{& \SetCell{bg=cF1} 1.00 & \SetCell{bg=cF2} 1.00 & \SetCell{bg=cF3} 1.03 & \SetCell{bg=cF4} 1.00 & \SetCell{bg=cF5} 1.04 & \SetCell{bg=cF6} 1.00 & \SetCell{bg=cF7} 1.00 & \SetCell{bg=cF8} 1.05 & 1.0158 & 0.0222}
\def\rowG{& \SetCell{bg=cG1} 1.06 & \SetCell{bg=cG2} 1.12 & \SetCell{bg=cG3} 1.00 & \SetCell{bg=cG4} 1.00 & \SetCell{bg=cG5} 1.00 & \SetCell{bg=cG6} 1.00 & \SetCell{bg=cG7} 1.05 & \SetCell{bg=cG8} 1.00 & 1.0289 & 0.0458}
\def\rowN{& 17 & 16 & 29 & 29 & 24 & 23 & 21 & 20}

\definecolor{cmap0}{rgb}{0.6980,0.1333,0.1333}
\definecolor{cmap1}{rgb}{0.6940,0.1616,0.1335}
\definecolor{cmap2}{rgb}{0.6899,0.1860,0.1337}
\definecolor{cmap3}{rgb}{0.6856,0.2077,0.1338}
\definecolor{cmap4}{rgb}{0.6812,0.2275,0.1339}
\definecolor{cmap5}{rgb}{0.6767,0.2459,0.1341}
\definecolor{cmap6}{rgb}{0.6721,0.2631,0.1342}
\definecolor{cmap7}{rgb}{0.6673,0.2794,0.1343}
\definecolor{cmap8}{rgb}{0.6624,0.2950,0.1344}
\definecolor{cmap9}{rgb}{0.6574,0.3099,0.1345}
\definecolor{cmap10}{rgb}{0.6522,0.3243,0.1346}
\definecolor{cmap11}{rgb}{0.6469,0.3382,0.1347}
\definecolor{cmap12}{rgb}{0.6414,0.3517,0.1348}
\definecolor{cmap13}{rgb}{0.6358,0.3648,0.1349}
\definecolor{cmap14}{rgb}{0.6300,0.3776,0.1349}
\definecolor{cmap15}{rgb}{0.6240,0.3901,0.1350}
\definecolor{cmap16}{rgb}{0.6179,0.4024,0.1350}
\definecolor{cmap17}{rgb}{0.6115,0.4144,0.1350}
\definecolor{cmap18}{rgb}{0.6050,0.4262,0.1350}
\definecolor{cmap19}{rgb}{0.5982,0.4378,0.1351}
\definecolor{cmap20}{rgb}{0.5913,0.4492,0.1351}
\definecolor{cmap21}{rgb}{0.5841,0.4604,0.1350}
\definecolor{cmap22}{rgb}{0.5767,0.4715,0.1350}
\definecolor{cmap23}{rgb}{0.5690,0.4825,0.1350}
\definecolor{cmap24}{rgb}{0.5611,0.4933,0.1350}
\definecolor{cmap25}{rgb}{0.5528,0.5040,0.1349}
\definecolor{cmap26}{rgb}{0.5443,0.5146,0.1348}
\definecolor{cmap27}{rgb}{0.5355,0.5251,0.1348}
\definecolor{cmap28}{rgb}{0.5264,0.5355,0.1347}
\definecolor{cmap29}{rgb}{0.5169,0.5458,0.1346}
\definecolor{cmap30}{rgb}{0.5070,0.5561,0.1345}
\definecolor{cmap31}{rgb}{0.4967,0.5662,0.1343}
\definecolor{cmap32}{rgb}{0.4859,0.5763,0.1342}
\definecolor{cmap33}{rgb}{0.4747,0.5863,0.1341}
\definecolor{cmap34}{rgb}{0.4630,0.5962,0.1339}
\definecolor{cmap35}{rgb}{0.4507,0.6061,0.1337}
\definecolor{cmap36}{rgb}{0.4378,0.6159,0.1335}
\definecolor{cmap37}{rgb}{0.4242,0.6257,0.1333}
\definecolor{cmap38}{rgb}{0.4098,0.6354,0.1331}
\definecolor{cmap39}{rgb}{0.3945,0.6451,0.1329}
\definecolor{cmap40}{rgb}{0.3782,0.6547,0.1327}
\definecolor{cmap41}{rgb}{0.3608,0.6643,0.1324}
\definecolor{cmap42}{rgb}{0.3420,0.6739,0.1321}
\definecolor{cmap43}{rgb}{0.3216,0.6834,0.1318}
\definecolor{cmap44}{rgb}{0.2991,0.6928,0.1315}
\definecolor{cmap45}{rgb}{0.2742,0.7023,0.1312}
\definecolor{cmap46}{rgb}{0.2458,0.7117,0.1309}
\definecolor{cmap47}{rgb}{0.2125,0.7211,0.1305}
\definecolor{cmap48}{rgb}{0.1716,0.7304,0.1302}
\definecolor{cmap49}{rgb}{0.1152,0.7397,0.1298}
\definecolor{cmap50}{rgb}{0.0000,0.7490,0.1294}

\def\plotColorAxis#1#2#3{%
  \begin{tikzpicture}
    \foreach \y [count=\i from 0] in {0,1,...,50} {
      \pgfmathsetmacro{\ybottom}{( \y   / 51.0) * #2}
      \pgfmathsetmacro{\ytop}{  ((\y+1) / 51.0) * #2}
      \fill[cmap\i, draw=none] (0, \ybottom #3) rectangle (#1 #3, \ytop #3);
    }
    \draw (#1 #3, 0.0000 * #2 #3) -- ++(0.1cm, 0) node[right, font=\tableFontSizeWiley] {1.44};
    \draw (#1 #3, 0.2500 * #2 #3) -- ++(0.1cm, 0) node[right, font=\tableFontSizeWiley] {1.33};
    \draw (#1 #3, 0.5000 * #2 #3) -- ++(0.1cm, 0) node[right, font=\tableFontSizeWiley] {1.22};
    \draw (#1 #3, 0.7500 * #2 #3) -- ++(0.1cm, 0) node[right, font=\tableFontSizeWiley] {1.11};
    \draw (#1 #3, 1.0000 * #2 #3) -- ++(0.1cm, 0) node[right, font=\tableFontSizeWiley] {1.00};
  \end{tikzpicture}
}

	\tableFontSizeWiley
	\newsavebox{\myTblBoxBeta}
	\sbox{\myTblBoxBeta}{
		\begin{tblr}[
			expand = \rowN\rowA\rowB\rowC\rowD\rowE\rowF\rowG
		]{
			colspec = {c r *{10}{c}},
			colsep = 3pt,
			vline{2}  = {4-Z}{0.7pt, dashed},
			vline{3}  = {1.2pt, solid},
			vline{11} = {1.2pt, solid},
			hline{2}  = {3-10}{0.7pt, dashed},
			hline{3}  = {2-10}{0.7pt, dashed},
			hline{4}  = {1-10}{1.2pt, solid},
			hline{6}  = {3-12}{0.7pt, solid},
		}
			& & \SetCell[c=8]{c} test case \\
			& & \testSetup{A}{1} & \testSetup{A}{2} & \testSetup{A}{3} & \testSetup{A}{4} & \testSetup{B}{1} & \testSetup{B}{2} & \testSetup{B}{3} & \testSetup{B}{4} & & \\
			& normalization in minutes \rowN \\
			\SetCell[r=7]{c,m} method & FWWS best \rowA \\
			& FWWS worst \rowB & mean & std. dev.\\
			& \ref{enum:WS_1} \rowC \\
			& \ref{enum:PS_1} \rowD \\
			& \ref{enum:PS_2} \rowE \\
			& \textbf{\ref{enum:PS_3}} \rowF \\
			& \textbf{\ref{enum:PS_4}} \rowG \\
		\end{tblr}
	}

	\pgfmathsetmacro{\myTableHeightCM}{
		0.035277777777782*\dimexpr\the\ht\myTblBoxBeta+\the\dp\myTblBoxBeta\relax 
	}
	
	\begin{adjustbox}{valign=m}
		\usebox{\myTblBoxBeta}
	\end{adjustbox}%
	\hspace{1.5em}%
	\begin{adjustbox}{valign=m}
		\plotColorAxis{0.5}{\myTableHeightCM}{cm}
	\end{adjustbox}
	
	\vspace{-0.3\baselineskip}
	\caption{
		The numerical values in the test case columns report normalized convergence times.
		For each test case, the convergence times are divided by the minimum convergence time achieved by any \acc{DM} method, which is given in the row \enquote{normalization in minutes}.
		The acronym \enquote{FWWS} stands for \enquote{fixed-weight \acc[l]{WS}}.
		For the \acc{MO} \acc{DM} methods, the cell colors visually support the normalized convergence times, with the color scale ranging from $1$ to $1.44$, the largest value of these methods.
		The last two columns report the mean and standard deviation of the normalized convergence times for each \acc{DM} method.
	}
	\label{tab:metrics_convergence}
\end{table*}

\begin{figure*}[htb]
	\centering
		\adjincludegraphics[
			width=1\linewidth, clip, 
			trim={{0.0\width} {0.06\height} {0.02\width} {0.08\height}}
		]{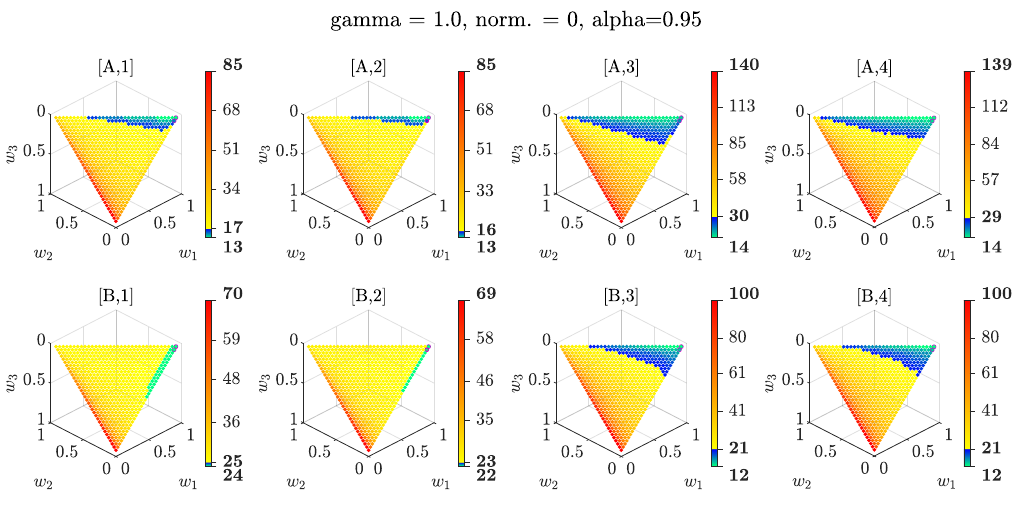}
	\caption{
		Convergence times for all fixed-weight \acc{WS} \accp{MPC}.
		For each test case, convergence times are color-coded using the corresponding color bar.
		Each color bar consists of two segments separated by bold tick labels.
		The bold labels at the bottom and top denote the best and worst convergence times of the fixed-weight \acc{WS} approach, whereas the bold labels in the middle are set to the convergence times of method \ref{enum:PS_3}, i.e., the normalized values in the corresponding row of \cref{tab:metrics_convergence} multiplied by the normalization values.
		The number of cases in which the fixed-weight \acc{WS} approach converges faster than \ref{enum:PS_3} (lower segments, green to blue) is consistently much smaller than the number of cases in which \ref{enum:PS_3} converges faster (upper segments, yellow to red).
		Moreover, \ref{enum:PS_3} is outperformed in all test cases simultaneously by only two weight vectors, marked by dots with a magenta outline.
	}
	\label{fig:convergence_time_beta_adapt}
\end{figure*}
\begin{figure*}[htb]
	\centering
		\adjincludegraphics[
			width=1\linewidth, clip, 
			trim={{0.0\width} {0.0\height} {0.005\width} {0.08\height}}
		]{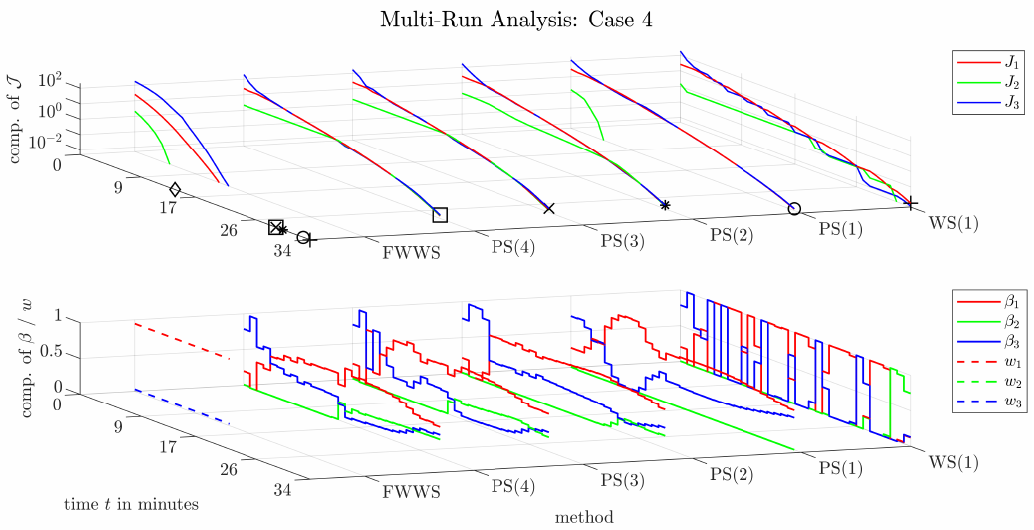}
	\caption{
		Both plots show results for the test case \testSetup{A}{4}. 
		It was chosen to compactly display the trajectories in three-dimensional plots where different $y$-components spatially separate the methods from one another.
		\textbf{Upper figure:} Cost trajectories.
		Each method is assigned a different marker. 
		The convergence time is also displayed on the timeline using these markers.
		\textbf{Lower figure:} Preference vector or weight vector trajectories.
	}
	\label{fig:results_sim_case_beta_adaptation_4}
\end{figure*}

In the process of formulating an optimization problem, it is common to identify several competing objectives. 
In such cases, the \acc{WS} method provides a straightforward and fundamental strategy to combine these objectives, particularly in the absence of more advanced methodologies.
A classic example is the infinite-time LQR problem where either the (diagonal) entries of the weighting matrices $Q$ and $R$ can be interpreted as the weights or the matrices themselves are scaled/weighted ($\param{w}_{\mathrm{WS},1} Q$, $\param{w}_{\mathrm{WS},2} R$).
\\
Hence, as comparative solutions, a \acc{WS} scalarization of \eqref{eq:MOMPC} is run with 595 realizations of fixed weights and $\realize{\J}_{\ub} = \boldsymbol{\infty}$, i.e.~the descent condition needed for \acc{MOMPC} is practically \emph{turned off}.
All realizations of the weights are obtained by uniformly sampling the standard simplex (excluding the boundary).
Since the weights are fixed, this results in a simple single-objective \acc{MPC} scheme.
Using these comparative solutions it is easier to put the solution of the \acc{MOMPC} methods using \accp{IM}-informed \acc{DM} into context.

The basis for assessing the numerical results is the time required until \cref{alg:MOMPC} terminates.
This time is determined by the algorithm iteration index $k$ and the step size $h_{\mathrm{macro}}$ (cf.~\cref{ssec:Setup_Numerical_Tests}), that is, $t_{\mathrm{conv}} = k h_{\mathrm{macro}}$.
In the remainder of this section, this quantity is referred to as the \enquote{convergence time} and is measured in minutes.

The numerical results with the \emph{initial descent condition approach b)} (cf.~\cref{step:MOMPC_init} of \cref{alg:MOMPC}) are summarized in \cref{tab:metrics_convergence} and \cref{fig:convergence_time_beta_adapt}.
Moreover, \cref{fig:results_sim_case_beta_adaptation_4} depicts the trajectories of the costs and the preference vector for test case \testSetup{A}{4}.

From the last two columns of \cref{tab:metrics_convergence}, it can be seen that among all \accp{IM}-informed \acc{DM} methods, \ref{enum:PS_2}, \ref{enum:PS_3}, and \ref{enum:PS_4} yield the best convergence times overall.
Among these three methods, \ref{enum:PS_2} performs best on average, while \ref{enum:PS_2} and \ref{enum:PS_3} show the smallest standard deviations.
At the same time, the mean differences between \ref{enum:PS_2}, \ref{enum:PS_3}, and \ref{enum:PS_4} amount to less than one sampling step $h_{\mathrm{macro}} = 1\,\unit{min}$, so that no single method emerges as a clear winner across all test cases.

One hypothesis for the slight advantage of \ref{enum:PS_2} is a higher robustness with respect to unfavorably oriented \acc{CHIM} hyperplanes%
\footnote{%
	Recall that, for \ref{enum:PS_2}, the quasi-normal is given by the vector connecting the barycentric center of the \acc{CHIM} and the \acc{UP} (cf.~\cref{item:quants_from_Phi_6} in \cref{list:objects_from_Phi}).
	Hence, the resulting shooting direction is independent of the orientation of the \acc{CHIM}.
}.
By contrast, the poorer performance of \ref{enum:WS_1} and \ref{enum:PS_1} may be attributable to a lower translation quality or coverage; cf.~\cref{ssec:Exemplary_Problems}.
In addition, \ref{enum:WS_1} exhibits pronounced chattering in the preference vector trajectories, as can be seen in \cref{fig:results_sim_case_beta_adaptation_4}.

A further classification of the results is provided by \cref{fig:convergence_time_beta_adapt}.
For each individual test case, there exist fixed weights such that the corresponding single-objective fixed-weight \acc{WS} \acc{MPC} yields a shorter convergence time than the \acc[l]{MO} approach (here represented by \ref{enum:PS_3}, cf.~\cref{fig:convergence_time_beta_adapt}).
However, when considering those fixed weights that are faster across all test cases simultaneously, only two weight vectors remain.
Interestingly, the favorable fixed weights tend to be located in regions with comparatively large $w_1$ and small $w_2,w_3$.
This is consistent with the hypothesis that, for the present problem, the first objective is coupled more strongly to the remaining ones, so that reducing $J_1$ may also accelerate the decay of $J_2$ and $J_3$.
If such structural effects are known a priori, they may be exploited to tune fixed weights more effectively.
However, obtaining this type of problem knowledge typically requires broader simulation studies or other computationally intensive analyses, and the resulting conclusions may still be scenario-dependent.
Hence, while fixed weights can be highly effective for a given operating scenario, their successful tuning may depend on prior problem-specific knowledge that is not readily available in practice.

The strength of the \acc{MO} approach does not primarily lie in outperforming the best tuned fixed-weight controller in every single test case.
Rather, its key advantage is the possibility to adapt the preference vector online.
This leads to a more broadly applicable controller and avoids the need for experienced engineers to tune weight parameters manually or to identify them through computationally expensive brute-force procedures.

\emph{Initial descent condition approach a)}, i.e., without regularization in the initial step, was also tested.
Since this directly affects the first decision, it also influences all subsequent decisions through the closed-loop evolution.
Although pathological \accp{PF} were observed in several test cases, the convergence times of approaches a) and b) differ by at most four minutes in each test case, and neither approach is systematically faster.
This suggests that the overall closed-loop performance is reasonably robust with respect to the initial decision in the present setup.
A more systematic investigation of how strongly an initially poor decision can affect the overall performance of \accp{IM}-informed \acc{MOMPC} remains an interesting topic for future work and is beyond the scope of this paper.

\section{Conclusion}
\label{sec:conclusion}
This paper developed a unified framework for \accp[l]{IM}-informed \acc[l]{MOMPC} that guarantees stabilization of a fixed point.
Six scalarization-based \acc[l]{DM} methods informed by the \accp[l]{IM} were considered, two of which are novel.
All of these methods aim to facilitate a real-time implementation of \acc[l]{MOMPC} by requiring only $n_J+1$ optimizations in two sequential stages.
In the two considered examples, certain strategies,
for example using the \emph{visual} normal for the \acc{NBI} method, 
provided a more robust translation of user preferences than the established \acc[l]{WS} and \acc[l]{NBI} methods with the \acc{CHIM} normal.
\\
Nevertheless, future research should include the evaluation of all methods for a large batch of benchmark \accp{MOOP}, both convex and concave.
Furthermore, there might exist additional \acc[l]{DM} methods closely related to the novel \acc{NCHIM} method that exhibit similarly promising numerical behavior.
One possible variation here could be to choose a different constant \acc[l]{SO} that does not correspond to the \acc[l]{NP}.

Next, the \acc[l]{DM} methods were successfully integrated into a quasi-infinite horizon \acc{MOMPC} scheme. 
A proof of asymptotic stability for the resulting closed-loop system was provided,
which is valid for all presented \acc[l]{DM} methods as long as the computed solution is not strictly dominated by the candidate solution.
This can be ensured by one additional cost evaluation.
Moreover, a practical method for constructing the necessary terminal ingredients was offered.

The paper concluded with a numerical case study that illustrated the closed-loop convergence of the proposed framework in all considered nominal test cases.
Furthermore, it was shown that, even in the presence of pathological \accp{PF} that may limit the information content of the \accp{IM}, the convergence of the resulting closed loop is not compromised in the considered setup.
The results showed that the best convergence times are achieved by an established and two novel \acc{DM} methods, without a clear overall winner.
They also indicated that suitably tuned fixed-weight \acc{WS} controllers can yield faster convergence times in individual test cases, but only a single fixed weight vector, identified by a brute-force search over all test cases, was faster in all of them.
In this sense, the main advantage of the \acc{MO} approach lies in the online adaptation of the preference vector, which avoids manual weight tuning and yields a more broadly applicable controller.
While this paper considered one particular way of adapting $\beta$ online, future research could develop and investigate alternative algorithms that exploit the information accumulated up to the current time in order to improve performance criteria other than convergence time.

\appendix
\crefalias{section}{appendix}

\section{Comparison of NBI and SRI} 
\label{sec:Comp_NBI_SRI}
The key difference between the two subforms is the parameterization of $(\param{\J}_{\SO}, \param{d})$.
	\Acc{NBI} samples via $\baryCO \in \FunitSimplex$, while
	\Acc{SRI} utilizes $\tau \in \FunitBox :=  [0,1]^{n_J-1}$.
The numerical example depicted in \cref{fig:mainfig} shows that the \acc{SRI} method provides better coverage of the entire \acc{PF}.

The unsampled parts of the \acc{PF} in \cref{fig:NBI_rays_coverage}, could be found by the following modifications:
i) scaling the \acc{CHIM} simplex, i.e.
$\param{\J}_{\SO}(\baryCO,c)=\Phi(1+c)\baryCO + c/n_{J} (\J_{\NP}-\J_{\UP})$, $c \geq 0$ 
or
ii) choosing suitable points on the hyperplane that contains the \acc{CHIM}.
However, in general, it is not clear how to choose $c$ or where these \enquote{suitable points} lie.
Thus, \acc{SRI} is a promising approach to be considered further as a scalarization method, especially if the sampling of the entire \acc{PF} is desired.

\section{Knee-Point}
\label{sec:KP}
In the following, the concept of \accp{FHP} is introduced, of which the \acc{KP} presented later constitutes a special case.
An \acc{FHP} is a point on the boundary of the feasible objective set $\imSetFeas$ from \eqref{eq:feasImage}, geometrically defined via the distance to a hyperplane.
Recall that all points $\chi \in \Rset^{n_J}$ belonging to an $(n_J-1)$-dimensional hyperplane have to satisfy the equation 
$\wFHP\tr \chi + \bFHP=0$, with normal vector $\wFHP\in \Rset^{n_J}$ and the scalar bias $\bFHP\in \Rset$.

The \emph{signed distance} between a hyperplane, parameterized by the normal vector $\wFHP$ and the bias $\bFHP$, and a point $\hat{\chi}$ along a vector $\dFHP$ is%
\begin{equation}
	\signD(\hat{\chi},\wFHP,\bFHP,\dFHP) = (\wFHP\tr \hat{\chi} + \bFHP) / \|\wFHP\| \cdot 1 / \cos(\theta),
	\label{eq:signedDist}
\end{equation}
where 
$\theta=\measuredangle(\wFHP,\dFHP)$ and $\cos(\theta) = \wFHP\tr\dFHP / (\|\wFHP\|\|\dFHP\|)$%
.%

\begin{bThm}[]{proposition}
	\label{prop:FHP}
	Let $(\wFHP,\bFHP)  \in \Rset^{n_J} \times \Rset$ define a hyperplane with $\wFHP$ scaled according to \cref{defi:scaling_normal} and let $\dFHP \in \Rset^{n_J}$ satisfy $|\measuredangle(\wFHP,\dFHP)|<90^\circ$ (i.e.~$\wFHP\tr\dFHP>0$). 
	\\
	Then, the \accp{OP}
	\begin{equation}
		\underset{\J \in \imSetFeas}{\max} \, \signD(\J,\wFHP,\bFHP,\dFHP)
		\quad \text{and} \quad
		\underset{\J \in \imSetFeas}{\min} \, J_{\FHP}(\J,\wFHP)
		\label{eq:FHP_equiv_form}
	\end{equation}
	with 
	$
		J_{\FHP}(\hat{\chi},\wFHP) = -\wFHP\tr \hat{\chi},
		\label{eq:J_FHP}
	$
	are equivalent.
\end{bThm}
\begin{proof}
	Observe that $\wFHP$, $\bFHP$ and $\theta$ are constant and $\cos(\theta)>0$ in the definition of the signed distance $\signD$ in \eqref{eq:signedDist}.
	Thus, the objective can be simplified and transformed from maximization to minimization by a simple sign change.
	In particular, the position of the hyperplane, as determined by $b$, has no influence on the solutions of the \accp{OP} in \eqref{eq:FHP_equiv_form}. 
\end{proof}

\begin{figure}[tb]
	\centering
		\adjincludegraphics[
			width=0.6\WidthBase, clip, 
			trim={{0\width} {0\height} {0\width} {0\height}}
		]{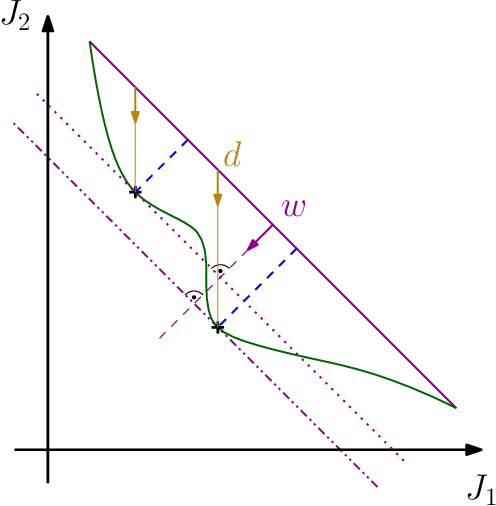}
	\caption{
		Example of a nonconvex \acc{PF} with two \accp{FHP}. 
		Both can be found by 
		either of the two \accp{OP} in \protect\eqref{eq:FHP_equiv_form}.
	}
	\label{fig:FHP_along_d}
\end{figure}

\Acc{FHP} can now be formally defined as (local) minima of $J_{\FHP}$, as is also illustrated in \cref{fig:FHP_along_d}.
In the implementation, the right problem in \eqref{eq:FHP_equiv_form} is replaced by $\min_{\xMO \in \decSetFeas} J_{\FHP}(\J(\xMO),\wFHP)$.

\begin{bThm}[Furthest-to-hyperplane point]{definition}
	\label{defi:char_FHP}%
	Let the assumptions of \cref{prop:FHP} hold.
	Then, $\hpAcc{\J} \in \imSetFeas$ is called a local \acc{FHP} if there exists an $\epsilon > 0$ such that 
	no $\hat{\J} \in \mathcal{B}_{\epsilon}(\hpAcc{\J}) \cap \imSetFeas$ satisfies 
	$J_{\FHP}(\hat{\J},\wFHP) < J_{\FHP}(\hpAcc{\J},\wFHP)$.
	Here, $\mathcal{B}_{\epsilon}(\hpAcc{\J})$ denotes an $\epsilon$-ball centered at $\hpAcc{\J}$.
	\\
	A local \acc{FHP} $\hpAcc{\J}$ is also a global \acc{FHP}, if there exists no $\hat{\J} \in \imSetFeas$ such that
	$J_{\FHP}(\hat{\J},\wFHP) < J_{\FHP}(\hpAcc{\J},\wFHP)$.

\end{bThm}

\begin{bThm}{remark}
	Readers familiar with the \acc{WS} method\footnote{Note that \acc{WS} scalarization for \acc{MOO} is formally introduced in \cref{sssec:WS}.} might have noticed the structural equivalence to the \acc{FHP} concept.
	In fact, if a hyperplane has a scaled normal $w$ with non-positive components, the \acc{WS} method can be used to find the point that (locally) has the largest signed distance to that hyperplane.
	This interpretation of the equivalence in \cref{prop:FHP} was also hinted at in \cite{das_characterizing_1999}, but not concretely stated.
	Note that a \acc{FHP} is not necessarily a point on the \acc{PF}, possibly only a dominated boundary point of $\imSetFeas$.
	This is because $\wFHP$ is allowed to have positive components, while \acc{WS} is restricted to non-negative weights by definition.
\end{bThm}

A special case of a \acc{FHP} is the so-called \acc[l]{KP}.
In line with \cite{das_characterizing_1999}, the \acc{KP} is defined to be the point in feasible image set that is furthest from the \acc{CHIM} hyperplane.
Hence, the \acc{KP} is a \acc{FHP} where $\wFHP=\wChim$.
\\
In \cite{das_normal-boundary_1998,das_characterizing_1999} it was proved that the \acc{KP} can be determined using either the \acc{WS} scalarization with suitable weights or the \acc{NBI} scalarization. 
The proof made use of multiplier theory and has the critical assumption that there exists a point on the \acc{CHIM} that is the start point of a vector that has the direction of the (quasi-)normal vector and ends on the \acc{KP}.
\\
In contrast, \cref{prop:FHP} and \cref{defi:char_FHP} demonstrate a geometric derivation of the weight that is needed to compute a \acc{FHP}.
In the case of the \acc{KP}, the weight resulting from this approach is a positive multiple of the weight obtained by the approach of \cite{das_characterizing_1999}.
Note that this result applies to any hyperplane, not just the one that corresponds to the \acc{CHIM}.
Furthermore, the geometric approach arguably enhances interpretability.

Depending on the geometric structure of the \acc{PF}, the \acc{CHIM} normal $\wChim$, that is scaled as described in \cref{defi:scaling_normal}, can have positive components.
As a result $J_{\FHP}(\J,\wChim)$ can be interpreted as the cost function of a \acc[l]{WS} scalarization with negative weights, which can yield dominated points.
The authors in \cite{das_normal-boundary_1998} state that such a weight vector \enquote{does not qualify}.

\section{Approximation of the Karcher Mean and Its Application}
\label{sec:Approx_Karcher_Mean}
The Karcher mean $P(w,S) \in \Rset^{n}$ associated with a weight vector $w\in\FunitSimplex$ is defined as the point on a $(n-1)$-dimensional Riemannian manifold $S$ that minimizes the weighted sum of squared geodesic distances $d_{\mathrm{geod}}(\cdot,\cdot,\cdot)$ to the anchor points $A_i$:
\begin{equation}
	P(w) = 
	\argmin_{X \in S} 
	\sum_{i=1}^{n} w_{i} \cdot d_{\mathrm{geod}}(S,X, A_i)^2.
	\label{eq:karcher_mean_cont}
\end{equation}
This definition (\cref{eq:karcher_mean_cont}) interprets the weight vector components $w_{i}$ as masses placed at the anchor points $A_i$, and the resulting point as the system's center of mass (or center of inertia) on the curved manifold.
Note that the Riemannian barycentric coordinates using the Karcher mean is a true generalization of barycentric coordinates in flat (Euclidean) space.
It can thus be shown that, if the manifold $S$ is a flat Euclidean space (i.e., $S = \Rset^{n-1}$), the geodesic distance $d_{\mathrm{geod}}(\cdot,\cdot,\cdot)$ becomes the standard Euclidean norm $\| \cdot \|_2$ and the minimization problem \eqref{eq:karcher_mean_cont} has a unique, closed-form solution that is mathematically identical to the standard linear barycentric coordinates $P(w) = \sum_{i=1}^{n} w_i A_i$.

The exact solution of \eqref{eq:karcher_mean_cont} cannot be determined, however, if the closed-form representation of $S$ is not available. 
Then, the solution of \eqref{eq:karcher_mean_cont} can be approximated through
\begin{equation}
	\tilde{P}(w) = 
	\argmin_{X \in \tilde{S}} 
	\sum_{i=1}^{n} w_{i} \cdot d_{\mathrm{geod}}(\tilde{S}_{\mathrm{mesh}},X, A_i)^2,
	\label{eq:karcher_mean_discr}
\end{equation}
where $\tilde{S} = \{X_1,\ldots,X_N\}$ is a fine sampling of $S$ and $\tilde{S}_{\mathrm{mesh}}$ is a surface mesh reconstructed from $\tilde{S}$.
The geodesic distance is obtained using the \emph{Fast Marching} method \cite{sethian_level_1999,sethian_fast_2000}. 
This method adapts the Dijkstra-like logic to continuous surfaces by allowing paths across the mesh triangles, rather than restricting them to the discrete mesh edges.
\\
To compute this mapping efficiently, the $N \times n$ geodesic distance matrix $D$ is first precomputed, where $D_{ji} = d_{\mathrm{geod}}(\tilde{S}_{\mathrm{mesh}},X_j, A_i)$.
Given $M$ weight vectors stored in an $n \times M$ matrix $W$, the $N \times M$ cost matrix $C$ for all points and all preferences is computed via a single matrix multiplication
$
C = (D \circ D) \, W,
$
where $\circ$ denotes the element-wise (Hadamard) product. 
The optimal point index $k_m$ for each $m$-th weight vector $w_{m}$ is then found by a \enquote{min} operation on the $m$-th column of $C$. 
This discretized approach provides a computationally feasible and highly efficient method to map any weight vector $w$ to its corresponding approximated Karcher mean.
In the present work, this construction is applied to \accp{PF} in order to obtain reference mappings from preference vectors to points on the corresponding \accp{PF}.
More specifically, the general approximation in \cref{eq:karcher_mean_discr} is instantiated for sampled \accp{PF} and their reconstructed meshes, as described next.

In order to obtain the mapping $\J_{\mathrm{des}}(\betaPref)$ for finitely many preference vectors $\betaPref$, the \acc{PF} is finely sampled to obtain $\imSetPareto[\tilde]$, and a mesh $\imSetPareto[\tilde][,\mathrm{mesh}]$ is reconstructed from the $N_{\mathrm{samples}}$ samples using the ball pivoting algorithm \cite{bernardini_ball-pivoting_1999,digne_analysis_2014}.
The mapping $\J_{\mathrm{des}}(\betaPref)$ is then computed via \cref{eq:karcher_mean_discr} by setting
\begin{align}
	\tilde{S} = \imSetPareto[\tilde], \quad
	\tilde{S}_{\mathrm{mesh}} = \imSetPareto[\tilde][,\mathrm{mesh}], \quad
	A_i = \JMO(\xMO_{i}\opt)
	\text{ with } i=1,\ldots,n_J
\end{align}
that is, the \accp{IM} are chosen as the anchor points.
This choice yields the natural correspondence $\J_{\mathrm{des}}(e_i) = \JMO(\xMO_{i}\opt)$.

\section{Notation Used in the Stability Proof}
\label{sec:conv_proof}
The proof of \cref{thm:conv_proof} uses the following notation adopted from \cite{grune_multiobjective_2019}.
\begin{itemize}
	\item $\xSys(\n)$ denotes the closed-loop state at time $t_\n$ with $\xSys_0 \equiv \xSys(0)$, whereas the subscript $j$ in $\xSys_j$ and $\uSys_j$ refers to the prediction steps in \eqref{eq:dyn_constraints}.
	\item $\uSeq = \big(\uSys(0),\ldots,\uSys(N-1)\big) \in \uSet^N$ denotes a sequence of $N$ control values
	and	$\uSeq(j)=\uSys(j)$ denotes the \enquote{extraction} of the $j$-th component.
	\item $\xSys^{\uSeq}(\xSys(\n))$ denotes the state sequence resulting from applying $\uSeq$ to the considered dynamics 
	and with the initial state $\param{\xSys}^0 \equiv \xSys(\n)$; 
	$\xSys^{\uSeq}(j, \xSys(\n))$ denotes the $j$-th element of that sequence.
	\item The input sequences 
		i) $\uSeq^{\PH{N}}_{\xSys(\n)}$ and 
		ii) $\uSeq^{\star \PH{, N}}_{\xSys(\n)}$  
		with their corresponding state sequences are a 
		i) feasible solution or the
		ii) applied solution
		to \eqref{eq:MOMPC} with the initial state $\param{\xSys}^0\equiv \xSys(\n)$.
		In particular, $\uSeq^{\star \PH{, N}}_{\xSys(\n)} = \big(\uSys_0\opt,\ldots,\uSys_{N-1}\opt\big)$ consists of the input components of the solution $\xMO\opt$ obtained in \cref{step:algo_FHP_all} of \cref{alg:MOMPC}, or of the candidate if it is applied instead (cf.~\cref{rem:proof}).
	\item Let 
	\begin{equation}
		\begin{aligned}
			J_i^{\PH{N}} \(\xSys(\n),\uSeq\) =
			&\sum_{j=0}^{N-1} \ell_i\(\xSys^{\uSeq}\(j,\xSys(\n)\) \! , \uSys(j)\) 
			\proofColDepAlignBefore
			\proofColDepAfter + F_i\(\xSys^{\uSeq}\(N,\xSys(\n)\)\)
		\end{aligned}
	\end{equation}
	define the cost from applying the input sequence $\uSeq$ to the initial state $\xSys(\n)$. 
\end{itemize}

\section{Water-Filling Strategy}
\label{sec:Water_Filling_Strategy}
In the water-filling strategy presented in \cref{alg:waterfilling_beta_update} local sensitivity information are leveraged to continuously adapt the allocation of resources, namely the individual components of $\betaPref$.
Note that, unlike the standard unweighted projection in \cite{michelot_finite_1986}, \cref{step:algo_wf_loop} computes a weighted projection onto the simplex.

As previously discussed, the algorithm could be interpreted as adapting $\betaPref$ in such a way as to achieve rapid, synchronized convergence.
Because \cref{alg:waterfilling_beta_update} operates as a local heuristic that allocates resources based on element-wise inverse sensitivities $W_i$, achieving this heavily relies on the underlying \acc{PF} topologies.
Specifically, the algorithm implicitly assumes that cost reductions scale independently ($\Delta \J_i \propto \betaPref[i](t_{k-1})$), thereby ignoring formal cross-dependencies.
However, despite this isolated calculation of $W_i$, the algorithm can implicitly exploit constructive cross-couplings (such as a high preference $\betaPref[1]$ driving down both $\J_1$ and a coupled cost $\J_2$) provided this simultaneous descent is immediately observable in the subsequent time step.

It must be emphasized that deriving an optimal algorithm for adapting $\betaPref$ is beyond the scope of this paper.
Instead, this specific strategy serves primarily as a proof of concept to highlight the general flexibility and the potential of iterative adaptations enabled by the \acc{MO} framework.

\begin{algorithm}[htbp]
	\caption{%
		Compute $\beta(t_k)$ via water-filling strategy
	}
	\label{alg:waterfilling_beta_update}
	\vspace{0.2\baselineskip}
	\emph{Input:} 
		dampening factor $\alpha \in [0,1]$; 
		current \acc{UP} $\J_{\UP}(t_k)$;
		previous cost values $\J\opt(t_{k-2})$, $\J\opt(t_{k-1})$; 
		previous preference vector $\beta(t_{k-1})$.
	\vspace{-0.5\baselineskip}
	\begin{steps}[label={\arabic*)},itemsep=0.5\baselineskip,labelindent=0.5ex]
		\item Coldstart: If $k<2$ set $\beta(t_k) = \J_{\UP}(t_k) / \operatorname{sum}(\J_{\UP}(t_k))$ and \textbf{terminate} the algorithm.
		\item 
			\label{step:algo_wf_metrics} 
			Calculate cost descent $\Delta \J = \J\opt(t_{k-2}) - \J\opt(t_{k-1})$.
			Apply lower bounds for numerical stability: $\Delta \tilde{\J}_i = \max(\Delta \J_i, \epsilon)$ and $\tilde{\beta}_{i} = \max(\beta_i(t_{k-1}), \epsilon)$
			with $\epsilon = 10^{-12}$.
			Compute dampened inverse sensitivity $W$:
			$
				W_{\mathrm{raw}, i} = \tilde{\beta}_{i} \, / \,  \Delta \tilde{\J}_i
			$
			and
			$
				W_i = (1-\alpha) W_{\mathrm{raw}, i} + \alpha \, / \, n_J \sum_{j=1}^{n_J} W_{\mathrm{raw}, j}
			$%
			,
			$i=1,\ldots,n_J$.
			
		\item 
			\label{step:algo_wf_init}
			Initialize active set $\mathcal{A} = \{1, \ldots, n_J\}$ and tentative preference vector $\beta_{\mathrm{tent}} = \boldsymbol{0}$.
			
		\item 
			\label{step:algo_wf_loop}
			Repeat to find the valid active set:
			\begin{enumerate}[label={\alph*)},itemsep=0.3\baselineskip]
				\item Compute common target cost $\J^{\mathrm{tar}}$ for $\mathcal{A}$:
				\\[0.2\baselineskip]
				$
					\J^{\mathrm{tar}} = 
					\displaystyle
					\frac{\sum_{i \in \mathcal{A}} \big(W_i \J_i\opt(t_{k-1})\big) - 1}{\sum_{i \in \mathcal{A}} W_i}
				$
				\item Compute tentative preferences for $i \in \mathcal{A}$:
				\\[0.2\baselineskip]
				$
					\beta_{\mathrm{tent}, i} = W_i \big(\J_i\opt(t_{k-1}) - \J^{\mathrm{tar}}\big)
				$
				\item 
				If $\beta_{\mathrm{tent}, i} \ge 0 \ \forall\, i \in \mathcal{A}$, go to \cref{step:algo_wf_finalize}.
				Else, remove violating indices: $\mathcal{A} \leftarrow \mathcal{A} \setminus \{i \mid \beta_{\mathrm{tent}, i} < 0\}$ and repeat.
			\end{enumerate}

		\item 
			\label{step:algo_wf_finalize}
			For all $i=1,\ldots,n_J$
			set $\beta_i(t_k)=\beta_{\mathrm{tent}, i}$ if $i \in \mathcal{A}$, else set $\beta_i(t_k)=0$.
			Normalize to ensure exact summation to 1: $\beta(t_k) \leftarrow \beta(t_k) / \operatorname{sum}(\beta(t_k))$.
	\end{steps}
	\vspace{-0.5\baselineskip}
	\emph{Output:} Preference vector $\beta(t_k)$.
\end{algorithm}


\section{Modeling of the Room Climate System}
\label{sec:model_num_ex}

\subsection{Modeling the Dynamics and the Financial Cost Rate}

\paragraph*{Relevant Quantities}
The state vector and its subvector are denoted as
\begin{equation}
    \CS{\xSys} = [ T \concSep h \concSep c ], 
    \quad \PS{\xSys} = [ T \concSep h ],
\end{equation}
where
\begin{itemize}
    \item $T$ is the room temperature in $^{\circ}\mathrm{C}$,
    \item $h$ is the absolute humidity in kg water per kg dry air ($\unit{kg/kg}$),
    \item $c$ is the CO$_2$ concentration in ppm.
\end{itemize}

Furthermore, the control vector 
and 
the vector of (piecewise constant) parameters or disturbances
is denoted as
\begin{equation}
	\begin{aligned}
		\uSys &= [ P_{\mathrm{floor}} \concSep \alpha_{\mathrm{mix}} \concSep Q_{\mathrm{vent}} \concSep \dT \concSep g_{\mathrm{humid}} ],		
		\\
		\text{and} \quad
		\param{d} &= [ N_{\mathrm{occ}} \concSep T_{\mathrm{out}} \concSep S_{\mathrm{solar}} \concSep \dot{m}_{h,\mathrm{dist}} ],
	\end{aligned}
\end{equation}
where
\begin{itemize}
    \item $P_{\mathrm{floor}}$ is the floor heating power input in $\unit{W=J/s}$ 
    \item $\alpha_{\mathrm{mix}}$ is the fraction of recirculated indoor air (0: all outside air, 1: all inside air),
    \item $Q_{\mathrm{vent}}$ is the ventilation air flow rate in $\unit{m^3/s}$,
    \item $\dT$ is the temperature adjustment of supply air by the HVAC system in $\unit{^{\circ}C}$,
    \item $g_{\mathrm{humid}}$ is the central humidity control input in kg water per second ($\unit{kg/s}$)
\end{itemize}
and
\begin{itemize}
    \item $N_{\mathrm{occ}}$ is the number of occupants,
    \item $T_{\mathrm{out}}$ is the outside air temperature in $\unit{^{\circ}C}$,
    \item $S_{\mathrm{solar}}$ is the solar irradiance in $\unit{W/m}^{2}$,
    \item $\dot{m}_{h,\mathrm{dist}}$ is an additional humidity disturbance in kg water per second ($\unit{kg/s}$).
\end{itemize}

\paragraph*{State and Input Bounds}
To ensure physical realism and to specify certain preferences, box constraints are imposed in the \acc{MPC} controller on the states 
\begin{equation}
	\begin{aligned}
		T &\in [17, 25] \  \unit{^{\circ}C},
		& h &\in 
		[4 \cdot 10^{-3}, 15 \cdot 10^{-3}] \ \unit{kg/kg},
		\\
		c &\in 
		[0, 1000] \ \unit{ppm}
	\end{aligned}
	\label{eq:state_bounds}
\end{equation}
and on the inputs
\begin{equation}
	\begin{aligned}
		P_{\mathrm{floor}} &\in 
		[0, 2000] \ \unit{W},
		& \airInside &\in [0,1],
		\\
		Q_{\mathrm{vent}} &\in [0,0.5] \ \unit{m^3/s},
		&\Delta T_{\mathrm{cond}} &\in [-5,5] \ \unit{^{\circ}C},
		\\
		g_{\mathrm{humid}} &\in 
		\scaleTime \cdot [-0.4,0.4] \ \unit{kg/h}
		.
	\end{aligned}
	\label{eq:input_bounds}
\end{equation}
with $\scaleTime = 1/3600 \, \unit{h/s}$. The lower and upper bounds are summarized in $\CS{\xSys}_{\lb}$, $\uSys_{\lb}$ and $\CS{\xSys}_{\ub}$, $\uSys_{\ub}$, respectively. 
Then $\CS{\xSet} = [\CS{\xSys}_{\lb}, \CS{\xSys}_{\ub}]$ and $\uSet = [\uSys_{\lb}, \uSys_{\ub}]$.

\paragraph*{Continuous-Time Dynamics}
Assuming that the disturbances are (piecewise) constant, the continuous-time dynamics are written as
\begin{equation}
	\dot{\CS{\xSys}}(t) = \fCT[\CS](\CS{\xSys}(t),\uSys(t),\param{d}),
\end{equation}
where $\fCT[\CS]$ is obtained by combining thermal, moisture, and mass-balance equations with 
\begin{subequations}
	\begin{align}
		C_{\mathrm{r}} \dot T &= 
		\begin{aligned}[t]
				&\frac{T_{\mathrm{out}} - T}{R_{\mathrm{rw}}}
			+ P_{\mathrm{floor}} + N_{\mathrm{occ}} P_{\mathrm{occ}}
			\\
			+ \, &\alpha_{\mathrm{solar}} A_{\mathrm{solar}} 	S_{\mathrm{solar}}
			+ \rho_{\mathrm{air}} c_{p,\mathrm{air}} Q_{\mathrm{vent}} 	(T_{\mathrm{sup}} - T),
		\end{aligned}
		\\
		m_{\mathrm{air}} \dot h &= 
		\begin{aligned}[t]
			\rho_{\mathrm{air}} Q_{\mathrm{vent}} (h_{\mathrm{sup}} - h)
			&+ g_{\mathrm{humid}} + N_{\mathrm{occ}} \dot m_{h,\mathrm{occ}}
			\\
			&+ \dot{m}_{h,\mathrm{dist}},
		\end{aligned}
		\\	
		m_{\mathrm{air}} \dot c &= 
			\rho_{\mathrm{air}} Q_{\mathrm{vent}} (c_{\mathrm{sup}} - c)
			+ N_{\mathrm{occ}} \dot m_{c,\mathrm{occ}}.
	\end{align}
	\label{eq:sys_dynamics_CT}%
\end{subequations}
with $m_{\mathrm{air}}=V_{\mathrm{r}}\rho_{\mathrm{air}}$.
The supply-air variables are defined by mixing and instantaneous conditioning:
\begin{subequations}
	\begin{align}
		T_{\mathrm{sup}} &= \airInside T + (1-\airInside) T_{\mathrm{out}} + \dT,\\
		h_{\mathrm{sup}} &= \airInside h + (1-\airInside) h_{\mathrm{out}},\\
		c_{\mathrm{sup}} &= \airInside c + (1-\airInside) c_{\mathrm{out}}.
	\end{align}
\end{subequations}

\paragraph*{Modeling the Power Consumption}
The instantaneous electrical power for the different HVAC actions are given by
\begin{equation}
	\begin{alignedat}{3}
		P_{\mathrm{vent}}(Q) &= p_{\mathrm{SFP}} Q, 
		&P_{\mathrm{cond}}(v,\Delta T) &= v \rho_{\mathrm{air}} c_{p,\mathrm{air}} |\Delta T|_{\pSmoothAbs},
		\\
		P_{\mathrm{humid}}(g) &= \frac{|g|_{\pSmoothAbs} h_{\mathrm{vap}}}{\eta_{\mathrm{hum}}},
		\quad
		&P_{\mathrm{heat}}(Q) &= \frac{P_{\mathrm{floor}}}{p_{\mathrm{COP}}} + P_{\mathrm{pump}},
	\end{alignedat}
	\label{eq:power_cons_funs}
\end{equation}
where a constant pump power of $P_{\mathrm{pump}}=30\unit{W}$ is assumed.
To ensure smoothness, the smooth absolute-value approximation $|\cdot|_{\pSmoothAbs} \equiv \sqrt{\cdot^2 + 1/\pSmoothAbs}$ with $L=50$ is employed.
The total cost rate
\begin{equation}
	l_{\mathrm{money}}(\uSys) = 
		C_{\mathrm{elec}}
		\begin{aligned}[t]
			\big( \hspace{1em}
				&P_{\mathrm{vent}}(Q_{\mathrm{vent}}) 
				+ P_{\mathrm{cond}}(Q_{\mathrm{vent}},\dT) 
				\\
				+ \, &P_{\mathrm{humid}}(
					g_{\mathrm{humid}}
				) 
				+ P_{\mathrm{heat}}(P_{\mathrm{floor}})
			\hspace{0.5em} \big)
		\end{aligned}
	\label{eq:l_money}\end{equation}
is obtained by summing up the electrical powers and multiplying it with $C_{\mathrm{elec}}$.
\\
Numerical values of the constant parameters are summarized in \cref{tab:params}.
\begin{table}[htbp]
	\centering
	\setlength{\tabcolsep}{4pt}
	\begin{tabular}{l c c l}
		\toprule
		Parameter & Symbol & Value & Unit \\
		\midrule
		Air density & $\rho_{\mathrm{air}}$ & 1.2 & $\unit{kg/m^3}$ \\
		Specific heat of air & $c_{p,\mathrm{air}}$ & 1005 & $\unit{J/(kg\,K)}$ \\
		Room volume & $V_{\mathrm{r}}$ & 40 & $\unit{m^3}$ \\
		Thermal capacity & $C_{\mathrm{r}}$ & $0.9 \cdot 10^{6}$ & $\unit{J/K}$ \\
		Thermal resistance & $R_{\mathrm{rw}}$ & 0.012 & $\unit{K/W}$ \\
		Occupant heat gain & $P_{\mathrm{occ}}$ & 120 & $\unit{W}$ \\
		CO$_2$ generation rate & $\dot{m}_{c,\mathrm{occ}}$ & $6.2$ & $\unit{ppm \cdot kg/s}$ \\
		Moisture generation rate & $\dot{m}_{h,\mathrm{occ}}$ & $1.39 \cdot 10^{-5}$ & $\unit{kg/s}$ \\
		Effective solar absorption coeff. & $\alpha_{\mathrm{solar}}$ & 0.6 & -- \\
		Effective area & $A_{\mathrm{solar}}$ & 1 & $\unit{m^2}$ \\
		Outdoor CO$_2$ & $c_{\mathrm{out}}$ & $500$ & $\unit{ppm}$ \\
		Outdoor humidity & $h_{\mathrm{out}}$ & $5 \cdot 10^{-3}$ & $\unit{kg/kg}$ 
		\\[1ex]
		Specific fan power param. & $p_{\mathrm{SFP}}$ & 1000 &  $\unit{W/(m^3/s)}$ \\
		Enthalpy of vaporization & $h_{\mathrm{vap}}$ & $2.45 \cdot 10^{6}$ & $\unit{J/kg}$ \\
		Humidifying device efficiency & $\eta_{\mathrm{hum}}$ & 0.7 & -- \\
		Heat pump coeff. of performance & $p_{\mathrm{COP}}$ & 2.5 & -- \\
		Electrical cost (non-base units) & $C_{\mathrm{elec}}$ & 0.30 & $\unit{\text{\euro}/(kWh)}$ \\
		\bottomrule
	\end{tabular}
	\caption{Physical parameters of the exemplary room model.}
	\label{tab:params}
\end{table}

\subsection{Preparations for the MPC Controller}
\paragraph*{Scaling}
The state $h$, and subsequently the input $g_{\mathrm{humid}}$, the disturbance $\dot{m}_{h,\mathrm{dist}}$ and the parameter $\dot{m}_{h,\mathrm{occ}}$, are scaled with $\scale_h = 10^{3}$.
Furthermore, the input $P_{\mathrm{floor}}$ is scaled with $\scale_{P_{\mathrm{floor}}} = 10^{-3}$.
For simplicity, no new notation is introduced for scaled quantities%
\footnote{From this point on the quantities are used in their scaled version.}%
.
Note that the system dynamics \eqref{eq:sys_dynamics_CT}, the state/input bounds \eqref{eq:state_bounds}/\eqref{eq:input_bounds} and the power consumption functions in \eqref{eq:power_cons_funs} are modified accordingly.

\paragraph*{Discretization of the Continuous-Time Dynamics}
The (scaled) continuous-time dynamics \eqref{eq:sys_dynamics_CT} are discretized via the second-order Runge-Kutta method of Heun \cite[chapter 4.2]{atkinson_numerical_2009} with the step-size $\hMicro=\hMacro/M$.
The composition of $M$ Runge-Kutta steps then results in a map $\fDT[\CS](\cdot,\cdot,\cdot)$ that generates the discrete-time dynamic constraints of the form
\begin{equation}
	\CS{\xSys}_{k+1} = \fDT[\CS]\(\CS{\xSys}_k,\uSys_k,\param{d}\) \, \text{ for } k \in \mathbb{N}_0
\end{equation}
given an initial state $\CS{\xSys}_0 $. 
The subsystem $\PS{\xSys}_{k+1} = \fDT[\PS](\PS{\xSys}_k,\uSys_k,\param{d})$ for $k \in \mathbb{N}_0$, which is decoupled from the state $c$, can be obtained by using the first two components of the complete system.

\paragraph*{Financially Optimal Controlled Steady State}
The aim is to find an optimal steady state (OSS) pair $(\CS{\xSys}\cSS, \uSys\cSS)$, cf.~\cref{asm:classical_stage_costs}, that minimizes the financial cost modeled by \eqref{eq:l_money} under some additional constraints on state and input. The optimal (controlled) steady state problem reads
\begin{equation}
	\begin{aligned}
		\min_{\CS{\xSys}, \uSys} \ &
		l_{\mathrm{money}}(\uSys)
		\\
		\text{s.t.} \quad 
		\CS{\xSys} &= \fDT[\CS](\CS{\xSys}, \uSys, \param{d}),
		\\
		\CS{\xSys} &\in [\CS{\xSys}_{\mathrm{OSS},\lb}, \CS{\xSys}_{\mathrm{OSS},\ub}] ,
		\\
		\uSys &\in [\uSys_{\mathrm{OSS},\lb}, \uSys_{\mathrm{OSS},\ub}] .
	\end{aligned}
\end{equation}
The solution to this optimization problem then provides the sought-after pair $\(\CS{\xSys}\cSS, \uSys\cSS\)$.
In the following, the units are omitted for most of the quantities.
To enforce pleasant ambient conditions, the bounds are set to
$\CS{\xSys}_{\mathrm{OSS},\lb} = \matInline{19.8, 5.35, 0}\tr$,
$\CS{\xSys}_{\mathrm{OSS},\ub} = \matInline{22.2, 8.65, 700}\tr$, 
resulting from a maximal deviation of 15\% of the respective state range in \eqref{eq:state_bounds} from $21\unit{^\circ C}$ and $7 \cdot 10^{-3} \, \unit{kg/kg}$\footnote{Note that this value refers to the unscaled humidity state.}
and a maximum $\text{CO}_2$ concentration that is acceptable in terms of health. 
Furthermore,
$\uSys_{\mathrm{OSS},\lb}=\uSys_{\lb} + \matInline{0, 0, 0.05, 0, 0}\tr$ and 
$\uSys_{\mathrm{OSS},\ub}=\uSys_{\ub}$.

As it is not reasonable to penalize a $\text{CO}_2$ concentration lower than $c\cSS = \CS{\xSys}\cSS[,3]$, the MPC controller penalizes only the state $\xSys$ with the dynamics
	$\PS{\xSys}_{k+1} = \fDT[\PS](\PS{\xSys}_k,\uSys_k,\param{d})$
	which is partially decoupled from $c$.

\subsection{Setup of the Numerical Tests}
\label{ssec:Setup_Numerical_Tests}
The parameters relevant for discretizing the dynamics are $h_{\mathrm{macro}}= 60 \unit{s}$ and $M=1$. 
The prediction horizon is $T=2\unit{h}$, which results in $N = \ceil{T/h_{\mathrm{macro}}}$.
The regularization term (cf.~\eqref{eq:w_reg_i}) is implemented with $\delta=0.002$.
\\
The weighting matrices are chosen as
\begin{equation}
	 \begin{alignedat}{3}
	 	\Q[][,1] &= \diag(\matInline{1,10^{-6}}), & \R[][,1] &= 10^{-6} \cdot \id{n_{\uSys}}, \\
	 	\Q[][,2] &= \diag(\matInline{10^{-6},1}), & \R[][,2] &= 10^{-6} \cdot \id{n_{\uSys}}, \\
	 	\Q[][,3] &= 10^{-6} \cdot \id{n_{\xSys}}, &\quad \R[][,3] &= \id{n_{\uSys}},
	 \end{alignedat}
\end{equation}
and the parameters of the LQR approach in \cite{gupta_application_2024} as $\rho_{\xSys}=10$, $\rho_{\uSys}=10$.
This choice results in a three-objective optimization problem. 
Suitable terminal ingredients are determined as described in \cref{ssec:terminal_ingredients}.
\\
The following eight test cases are considered:
Either
\circled{A} $\realize{d} = \matInline{2, 16, 200, 0}\tr$ or 
\circled{B} $\realize{d} = \matInline{4, 28, 800, 0.5}\tr$ 
with one of the four initial states
\circled{1} $\realize{\xSys}^0 = \matInline{24, 7.0}\tr$,\ColDepLinebreak
\circled{2} $\realize{\xSys}^0 = \matInline{24, 4.5}\tr$,
\circled{3} $\realize{\xSys}^0 = \matInline{17, 7.0}\tr$,
\circled{4} $\realize{\xSys}^0 = \matInline{17, 4.5}\tr$.

\if\AddAckn1
\section*{Acknowledgments}
A preprint version of this manuscript is available on arXiv at \url{https://arxiv.org/abs/2510.23454}.

The authors would like to thank the anonymous reviewers for their insightful comments and constructive feedback, which significantly improved the quality of this manuscript.

The authors used Claude Opus 5.5 (Anthropic, October 2026) for language polishing and to obtain feedback on the cohesion of the manuscript and the argument structure of the proof of \cref{thm:conv_proof}.
All suggestions were critically reviewed by the authors, who take full responsibility for the content of this article.
\fi

\addtolength{\textheight}{-\ShortenLastPage} 
	 
\bibliography{26_IMI_sMOMPC.bib}

\end{document}